\documentclass[12pt]{amsart}
\usepackage{amssymb}
\usepackage{amsmath}
\usepackage{amscd}
\usepackage[all]{xy}
\theoremstyle{plain}
\newtheorem{theorem}{Theorem}[section]
\newtheorem{corollary}[theorem]{Corollary}
\newtheorem{proposition}[theorem]{Proposition}
\newtheorem{lemma}[theorem]{Lemma}
{\theoremstyle{remark}

\newtheorem{remark}[theorem]{Remark}}
{\theoremstyle{definition}
\newtheorem{definition}[theorem]{Definition}

\newtheorem{example}[theorem]{Example}}

\newcommand{\benu}{\begin{enumerate}\renewcommand{\labelenumi}{{\rm (\roman{enumi})}}\renewcommand{\itemsep}{0pt}}
\newcommand{\eenu}{\end{enumerate}}

\newcommand{\N}{\mathbb{N}}
\newcommand{\Z}{\mathbb{Z}}

\newcommand{\C}{\mathbb{C}}
\newcommand{\T}{\mathbb{T}}
\newcommand{\e}{\varepsilon}
\newcommand{\tT}{\widetilde{T}}
\newcommand{\cK}{{\mathcal K}}
\newcommand{\cL}{{\mathcal L}}
\newcommand{\cM}{{\mathcal M}}
\newcommand{\K}{\mathbb{K}}
\newcommand{\M}{\mathbb{M}}

\newcommand{\cO}{{\mathcal O}}
\newcommand{\cT}{{\mathcal T}}
\newcommand{\ip}[2]{\langle{#1},{#2}\rangle}

\newcommand{\s}[3]{{{#1}^{#2}_{\textnormal{#3}}}}
\newcommand{\rs}[1]{{\textnormal #1}}

\newcommand{\F}{{\mathcal F}}
\newcommand{\G}{{\mathcal G}}
\DeclareMathOperator{\id}{id}
\DeclareMathOperator{\Aut}{Aut}

\DeclareMathOperator{\diag}{diag}
\DeclareMathOperator{\dom}{dom}
\DeclareMathOperator{\Ad}{Ad}

\begin{document}
\title[A class of $C^*$-algebras II, examples]
{A class of \boldmath{$C^*$}-algebras generalizing 
both graph algebras and homeomorphism \boldmath{$C^*$}-algebras II, 
examples}
\author[Takeshi KATSURA]{Takeshi KATSURA}
\address{Department of Mathematics, 
Hokkaido University, Kita 10, Nishi 8, 
Kita-Ku, Sapporo, 060-0810, JAPAN}
\email{katsura@math.sci.hokudai.ac.jp}
\thanks{The author is partially supported by Research Fellowship 
for Young Scientists of the Japan Society for the Promotion of Science.} 
\date{}

\begin{abstract}
We show that the method to construct $C^*$-algebras 
from topological graphs, introduced in our previous paper, 
generalizes many known constructions. 
We give many ways to make 
new topological graphs from old ones, 
and study the relation of $C^*$-algebras 
constructed from them. 
We also give a characterization
of our $C^*$-algebras 
in terms of their representation theory. 
\end{abstract}

\maketitle

\setcounter{section}{-1}

\section{Introduction}

In a previous paper \cite{Ka1}, 
we introduced the notion of topological graphs and 
a method to construct $C^*$-algebras from them. 
Topological graphs generalize ordinary graphs 
and homeomorphisms on locally compact spaces, 
and our method for constructing $C^*$-algebras from them 
generalizes the constructions of graph algebras 
and homeomorphism algebras (see \cite{Ka1} for detail). 
In this paper, 
we give many ways to make new topological graphs 
from old ones, 
and study $C^*$-algebras constructed from them. 
We also see that the way of constructing $C^*$-algebras 
from topological graphs and the class of such $C^*$-algebras 
relates many known constructions and classes 
besides graph algebras and homeomorphism algebras. 
In \cite{Ka1}, 
we show that our $C^*$-algebras are always nuclear 
and satisfy the Universal Coefficient Theorem. 
So far, we know of no examples which satisfy these two conditions, 
but are not in our class. 
Almost all ``classifiable'' $C^*$-algebras 
can be obtained as $C^*$-algebras of topological graphs. 
Thus our $C^*$-algebras are useful to study the structure 
of classifiable $C^*$-algebras. 

In Section \ref{SecPre}, 
we recall definitions and results 
in our previous paper \cite{Ka1}. 
In Section \ref{SecFactor}, 
we define factor maps between two topological graphs and 
show that these give $*$-ho\-mo\-mor\-phisms 
between $C^*$-algebras associated with them. 
In Section \ref{SecToepfam}, 
we investigate $C^*$-algebras $C^*(T)$ 
generated by Toeplitz pairs $T=(T^0,T^1)$. 
Thanks to this investigation, 
we get a characterization of our $C^*$-algebra $\cO(E)$ 
without using the space $\s{E}{0}{rg}$ 
(Proposition \ref{OEminimum}). 
We use the results here in the next paper \cite{Ka2}. 
In Section \ref{SecProj}, 
we define projective systems of topological graphs 
and their projective limits, 
and study how these relate to $C^*$-algebras $\cT(E)$ and $\cO(E)$. 
In Section \ref{Secsubalg} and Section \ref{SecMorita}, 
we give methods to make a new topological graph 
from given one so that the $C^*$-algebras they define 
are strongly Morita equivalent. 
Section \ref{SecOther} is devoted 
to give other operations to make new topological graphs. 
In the final three sections, 
we discuss examples which show how 
a number of $C^*$-algebras that appear in the literature 
may be realized as topological graph algebras.

\section{Preliminaries}\label{SecPre}

\begin{definition}
A {\em topological graph} $E=(E^0,E^1,d,r)$ consists of 
two locally compact spaces $E^0$ and $E^1$, 
and two maps $d,r\colon E^1\to E^0$, 
where $d$ is locally homeomorphic 
and $r$ is continuous.
\end{definition}

We regard an element $v$ of $E^0$ as a vertex, 
and an element $e$ of $E^1$ as a directed edge 
pointing from its domain $d(e)\in E^0$ 
to its range $r(e)\in E^0$. 
For a topological graph $E=(E^0,E^1,d,r)$, 
the triple $(E^1,d,r)$ is called 
a {\em topological correspondence} on $E^0$, 
which is considered as a generalization of a continuous map. 
By the local homeomorphism $d\colon E^1\to E^0$, 
$E^1$ is ``locally'' isomorphic to $E^0$, 
and the pair $(E^1,d)$ defines a ``domain'' of a continuous map $r$
which is ``locally'' a continuous map from $E^0$ to $E^0$. 

Let us denote by $C_d(E^1)$ 
the set of continuous functions $\xi$ on $E^1$ such that 
$\ip{\xi}{\xi}(v)=\sum_{e\in d^{-1}(v)}|\xi(e)|^2<\infty$ for any $v\in E^0$ 
and $\ip{\xi}{\xi}\in C_0(E^0)$.
For $\xi,\eta\in C_d(E^1)$ and $f\in C_0(E^0)$, 
we define $\xi f\in C_d(E^1)$ and $\ip{\xi}{\eta}\in C_0(E^0)$ by 
$$(\xi f)(e)=\xi (e)f(d(e)) \mbox{ for } e\in E^1$$
$$\ip{\xi}{\eta}(v)=\sum_{e\in d^{-1}(v)}\overline{\xi(e)}\eta(e)
\mbox{ for } v\in E^0.$$
With these operations,  
$C_d(E^1)$ is a (right) Hilbert $C_0(E^0)$-module 
(\cite[Proposition 1.10]{Ka1}).
We define a map $\pi\colon C_b(E^1)\to\cL(C_d(E^1))$ 
by $(\pi(f)\xi)(e)=f(e)\xi(e)$ 
for $f\in C_b(E^1)$, $\xi\in C_d(E^1)$ and $e\in E^1$, 
where $C_b(E^1)$ is the set of all bounded continuous functions on $E^1$. 
We have $\pi(f)\in \cK(C_d(E^1))$ if and only if $f\in C_0(E^1)$ 
(\cite[Proposition 1.17]{Ka1}). 
We define a left action $\pi_r$ of $C_0(E^0)$ on $C_d(E^1)$ 
by $\pi_r(f)=\pi(f\circ r)\in \cL(C_d(E^1))$ for $f\in C_0(E^0)$. 
Thus we get a $C^*$-correspondence $C_d(E^1)$ over $C_0(E^0)$. 

We set $d^0=r^0=\id_{E^0}$ and $d^1=d, r^1=r$.
For $n=2,3,\ldots$, 
we recursively define a space $E^n$ of paths with length $n$ 
and domain and range maps $d^n,r^n\colon E^n\to E^0$ by
$$E^n=\{(e',e)\in E^1\times E^{n-1}\mid d^1(e')=r^{n-1}(e)\},$$ 
$d^n((e',e))=d^{n-1}(e)$ and $r^n((e',e))=r^1(e')$. 
For each $n\in\N$, 
$d^n$ is a local homeomorphism of $E^n$ to $E^0$ 
and of course $r^n$ is continuous 
(the triple $(E^n,d^n,r^n)$ is the $n$-times composition of 
the topological correspondence $(E^1,d,r)$ on $E^0$, 
see \cite[Section 1]{Ka1}). 
Thus we may define a $C^*$-correspondence $C_{d^n}(E^n)$ over $C_0(E^0)$ 
in a fashion similar to the definition of $C_d(E^1)$. 
We have that 
$C_{d^{n+m}}(E^{n+m})\cong C_{d^n}(E^n)\otimes C_{d^{m}}(E^{m})$ 
as $C^*$-correspondences over $C_0(E^0)$ 
for any $n,m\in\N=\{0,1,2,\ldots\}$.
As long as no confusion arises, 
we omit the superscript $n$ 
and simply write $d,r$ for $d^n,r^n$. 

\begin{definition}\label{DefTpl}
Let $E=(E^0,E^1,d,r)$ be a topological graph.
A {\em Toeplitz $E$-pair} in a $C^*$-algebra $A$ 
is a pair of maps $T=(T^0,T^1)$ consisting of 
a $*$-ho\-mo\-mor\-phism $T^0\colon C_0(E^0)\to A$ and 
a linear map $T^1\colon C_d(E^1)\to A$ satisfying 
\benu
\item $T^1(\xi)^*T^1(\eta)=T^0(\ip{\xi}{\eta})$ for $\xi,\eta\in C_d(E^1)$, 
\item $T^0(f)T^1(\xi)=T^1(\pi_r(f)\xi)$ 
for $f\in C_0(E^0)$ and $\xi\in C_d(E^1)$.
\eenu

We denote by $\cT(E)$ the universal $C^*$-algebra
generated by a Toeplitz $E$-pair.
\end{definition}

For a Toeplitz $E$-pair $T=(T^0,T^1)$, 
the equation $T^1(\xi)T^0(f)=T^1(\xi f)$ 
holds automatically 
from the condition (i). 
We write $C^*(T)$ to denote the $C^*$-algebra 
generated by the images of the maps $T^0$ and $T^1$. 
We define a $*$-ho\-mo\-mor\-phism 
$\varPhi\colon \cK(C_d(E^1))\to C^*(T)$
by $\varPhi(\theta_{\xi,\eta})=T^1(\xi)T^1(\eta)^*$ 
for $\xi,\eta\in C_d(E^1)$. 
We say that a Toeplitz $E$-pair $T=(T^0,T^1)$ is {\em injective} 
if $T^0$ is injective. 
If a Toeplitz $E$-pair $T=(T^0,T^1)$ is injective,
then $T^1$ and $\varPhi$ are isometric.

\begin{definition}
Let $E=(E^0,E^1,d,r)$ be a topological graph.
We define three open subsets 
$\s{E}{0}{sce},\s{E}{0}{fin}$ and $\s{E}{0}{rg}$ of $E^0$ by
$\s{E}{0}{sce}=E^0\setminus\overline{r(E^1)}$,
\begin{align*}
\s{E}{0}{fin}=\{v\in E^0\mid\mbox{ there exists}&\mbox{ a neighborhood } 
V \mbox{ of } v\\
&\mbox{ such that }
r^{-1}(V)\subset E^1 \mbox{ is compact}\},
\end{align*}
and $\s{E}{0}{rg}=\s{E}{0}{fin}\setminus\overline{\s{E}{0}{sce}}$.
We define two closed subsets $\s{E}{0}{inf}$ and $\s{E}{0}{sg}$ of $E^0$ by 
$\s{E}{0}{inf}=E^0\setminus \s{E}{0}{fin}$ and 
$\s{E}{0}{sg}=E^0\setminus \s{E}{0}{rg}$. 
\end{definition}

We have $\s{E}{0}{sg}=\s{E}{0}{inf}\cup\overline{\s{E}{0}{sce}}$. 
A vertex in $\s{E}{0}{sce}$ is called a {\em source}.
Vertices in $\s{E}{0}{rg}$ are said to be {\em regular},
and those in $\s{E}{0}{sg}$ are said to be {\em singular}. 
The following is a characterization of regular vertices. 

\begin{lemma}[{\cite[Proposition 2.8]{Ka1}}]\label{E0r}
For $v\in E^0$, we have $v\in\s{E}{0}{rg}$ 
if and only if there exists a neighborhood $V$ of $v$ such that 
$r^{-1}(V)$ is compact and $r(r^{-1}(V))=V$. 
\end{lemma}

Note that if it exists, such a neighborhood $V$ is compact, 
and every compact neighborhood $V'$ of $v$ contained in $V$ 
satisfies the same conditions. 
We have that $\ker\pi_r=C_0(\s{E}{0}{sce})$ and 
$\pi_r^{-1}(\cK(C_d(E^1)))=C_0(\s{E}{0}{fin})$ 
(\cite[Proposition 1.24]{Ka1}).
Hence the restriction of $\pi_r$ to $C_0(\s{E}{0}{rg})$ is 
an injection into $\cK(C_d(E^1))$.

\begin{definition}
A Toeplitz $E$-pair $T=(T^0,T^1)$ is called 
a {\em Cuntz-Krieger $E$-pair} 
if $T^0(f)=\varPhi(\pi_r(f))$ for any $f\in C_0(\s{E}{0}{rg})$. 

The universal $C^*$-algebra 
generated by a Cuntz-Krieger $E$-pair $t=(t^0,t^1)$ 
is denoted by $\cO(E)$. 
\end{definition}

In \cite{Ka6}, 
the author suggests the way to associate a $C^*$-algebra $\cO_X$ 
for each $C^*$-correspondence $X$, 
which is a modification of the construction of Pimsner algebras 
in \cite{P}. 
The $C^*$-algebra $\cO(E)$ is nothing but the $C^*$-algebra $\cO_{C_d(E^1)}$ 
associated with the $C^*$-correspondence $C_d(E^1)$. 

Since $t=(t^0,t^1)$ is injective (\cite[Proposition 3.7]{Ka1}), 
$\varphi\colon \cK(C_d(E^1))\to \cO(E)$ is injective. 
By the universality of $\cO(E)$, 
there exists an action $\beta\colon \T\curvearrowright \cO(E)$ defined by 
$\beta_z(t^0(f))=t^0(f)$ and $\beta_z(t^1(\xi))=zt^1(\xi)$ 
for $f\in C_0(E^0)$, $\xi\in C_d(E^1)$ and $z\in\T$. 
The action $\beta$ is called the {\em gauge action}. 
We say that a Toeplitz $E$-pair $T$ {\em admits a gauge action} 
if there exists an automorphism $\beta'_z$ on $C^*(T)$ 
with $\beta'_z(T^0(f))=T^0(f)$ and $\beta'_z(T^1(\xi))=zT^1(\xi)$
for every $z\in\T$. 
Note that if such automorphisms $\beta'_z$ exist, 
then $\beta'\colon \T\ni z\mapsto \beta'_z\in\Aut(C^*(T))$ 
becomes automatically a strongly continuous homomorphism. 
The following proposition is called the gauge-invariant uniqueness theorem.

\begin{proposition}[{\cite[Theorem 4.5]{Ka1}}]\label{GIUT}
For a topological graph $E$ and 
a Cuntz-Krieger $E$-pair $T$, 
the natural surjection $\cO(E)\to C^*(T)$ is an isomorphism 
if and only if $T$ is injective and admits a gauge action. 
\end{proposition}

\section{Factor maps}\label{SecFactor}

In this section, 
we define factor maps between two topological graphs and 
show that these give $*$-ho\-mo\-mor\-phisms 
between $C^*$-algebras associated with them.
Let us take two topological graphs 
$E=(E^0,E^1,d_E,r_E)$ and $F=(F^0,F^1,d_F,r_F)$, 
and fix them. 
For a locally compact space $X$, 
we denote by $\widetilde{X}=X\cup\{\infty\}$ 
the one-point compactification of $X$. 
We consider elements of $C_0(X)$ 
as continuous functions on $\widetilde{X}$ 
vanishing at $\infty\in \widetilde{X}$. 

\begin{definition}\label{DefFactor}
A {\em factor map} from $F$ to $E$ is 
a pair $m=(m^0,m^1)$ consisting of continuous maps 
$m^0\colon \widetilde{F}^0\to \widetilde{E}^0$ 
and $m^1\colon \widetilde{F}^1\to \widetilde{E}^1$ 
which send $\infty$ to $\infty$, such that 
\benu
\item For every $e\in F^1$ with $m^1(e)\in E^1$, 
we have $r_E(m^1(e))=m^0(r_F(e))$ and $d_E(m^1(e))=m^0(d_F(e))$. 
\item If $e'\in E^1$ and $v\in F^0$ satisfies $d_E(e')=m^0(v)$, 
then there exists a unique element $e\in F^1$ 
such that $m^1(e)=e'$ and $d_F(e)=v$. 
\eenu
\end{definition}

The term `factor map' comes from 
cosidering topological graphs as dynamical systems. 
Note that the domain and range maps $d_E,r_E\colon E^1\to E^0$ 
of a topological graph $E=(E^0,E^1,d_E,r_E)$ 
may not extend to 
continuous maps from $\widetilde{E}^1$ to $\widetilde{E}^0$ in general. 

\begin{lemma}
Let $m=(m^0,m^1)$ be a factor map from $F$ to $E$. 
Then $m$ induces a pair of maps $\mu=(\mu^0,\mu^1)$ by 
\begin{align*}
\mu^0\colon C_0(E^0)\ni f&\mapsto f\circ m^0\in C_0(F^0)\\
\mu^1\colon C_{d_E}(E^1)\ni \xi&\mapsto \xi\circ m^1\in C_{d_F}(F^1),
\end{align*}
which is a morphism in the sense of \cite[Definition 2.3]{Ka7}, 
that is, 
$\mu^0$ is a $*$-ho\-mo\-mor\-phism, 
$\ip{\mu^1(\xi)}{\mu^1(\eta)}=\mu^0(\ip{\xi}{\eta})$ 
and $\pi_{r_F}(\mu^0(f))\mu^1(\xi)=\mu^1(\pi_{r_E}(f)\xi)$
for $\xi,\eta\in C_{d_E}(E^1)$ and $f\in C_0(E^0)$. 
\end{lemma}

\begin{proof}
Clearly $f\mapsto f\circ m^0$ defines a $*$-ho\-mo\-mor\-phism $\mu^0$. 
Take $\xi,\eta\in C_{d_E}(E^1)$, 
and we will show the equality 
$\ip{\mu^1(\xi)}{\mu^1(\eta)}=\mu^0(\ip{\xi}{\eta})$. 
Let us take $v\in F^0$. 
For $e\in (d_F)^{-1}(v)\cap (m^1)^{-1}(E^1)$, 
we have $m^1(e)\in (d_E)^{-1}(m^0(v))$ 
by the condition (i) in Definition \ref{DefFactor}.
Conversely for $e'\in (d_E)^{-1}(m^0(v))$, 
there exists unique $e\in (d_F)^{-1}(v)$ with $m^1(e)=e'$. 
Hence the map 
$$(d_F)^{-1}(v)\cap (m^1)^{-1}(E^1)\ni e \mapsto 
m^1(e)\in (d_E)^{-1}(m^0(v))$$
is bijective. 
Therefore we have 
\begin{align*}
\ip{\mu^1(\xi)}{\mu^1(\eta)}(v)
&=\sum_{e\in (d_F)^{-1}(v)}\overline{\mu^1(\xi)(e)}\,\mu^1(\eta)(e)\\
&=\!\!\sum_{e\in (d_F)^{-1}(v)\cap (m^1)^{-1}(E^1)}\!\!
  \overline{\xi(m^1(e))}\,\eta(m^1(e))\\
&=\!\!\sum_{e'\in (d_E)^{-1}(m^0(v))}\!\!\overline{\xi(e')}\eta(e')\\
&=\ip{\xi}{\eta}(m^0(v))\\
&=\mu^0(\ip{\xi}{\eta})(v)
\end{align*}
(note that when $m^0(v)=\infty$ 
we have $(d_E)^{-1}(m^0(v))=\emptyset$, 
hence in this case the both hands of the fourth equality are zero).
Thus we get 
$\ip{\mu^1(\xi)}{\mu^1(\eta)}=\mu^0(\ip{\xi}{\eta})$. 
By taking $\xi=\eta$ in the above equality, 
we see that $\mu^1$ is well-defined. 
Finally it is easy to see 
$\mu^1(\pi_{r_E}(f)\xi)=\pi_{r_F}(\mu^0(f))\mu^1(\xi)$ 
for $f\in C_0(E^0)$ and $\xi\in C_{d_E}(E^1)$. 
We are done. 
\end{proof}

By this lemma, 
we get the following proposition.

\begin{proposition}\label{functT}
Let $\mu^0\colon C_0(E^0)\to C_0(F^0)$, 
$\mu^1\colon C_{d_E}(E^1)\to C_{d_F}(F^1)$ be maps 
defined from a factor map $m=(m^0,m^1)$ from $F$ to $E$ as above. 
Then there exists a unique $*$-ho\-mo\-mor\-phism $\mu\colon \cT(E)\to\cT(F)$ 
such that $\mu\circ T_E^i=T_F^i\circ\mu^i$ for $i=0,1$, 
where $T_E=(T_E^0,T_E^1)$ and $T_F=(T_F^0,T_F^1)$ 
are the universal Toeplitz $E$-pair in $\cT(E)$ 
and the universal Toeplitz $F$-pair in $\cT(F)$, respectively. 
\end{proposition}

\begin{proposition}\label{comp}
Let $E,F,G$ be topological graphs, 
and $m,n$ be factor maps from $F$ to $E$ and $G$ to $F$ respectively. 
We set $(m\circ n)^i=m^i\circ n^i\colon \widetilde{G}^i\to 
\widetilde{E}^i$ for $i=0,1$. 
Then $m\circ n=((m\circ n)^0,(m\circ n)^1)$ is a factor map from $G$ to $E$, 
and the $*$-ho\-mo\-mor\-phism 
$\omega\colon \cT(E)\to\cT(G)$ defined from $m\circ n$ 
is the composition of $\mu\colon \cT(E)\to\cT(F)$ 
and $\nu\colon \cT(F)\to\cT(G)$ 
which are defined from $m$ and $n$ respectively. 
\end{proposition}

\begin{proof}
Take $e\in G^1$ with $(m^1\circ n^1)(e)\in E^1$. 
Then $n^1(e)\in F^1$. 
Hence we have $d_F(n^1(e))=n^0(d_G(e))$ and $r_F(n^1(e))=n^0(r_G(e))$. 
Since $m^1(n^1(e))\in E^1$, we have 
$d_E(m^1(n^1(e)))=m^0(d_F(n^1(e)))$ and $r_E(m^1(n^1(e)))=m^0(r_F(n^1(e)))$. 
Therefore we get 
$$d_E((m^1\circ n^1)(e))=(m^0\circ n^0)(d_G(e))\,,\quad
r_E((m^1\circ n^1)(e))=(m^0\circ n^0)(r_G(e)).$$ 
Take $e\in E^1$ and $v\in G^0$ with $d_E(e)=(m^0\circ n^0)(v)$.
Since $m$ is a factor map, 
there exists unique $e'\in F^1$ 
with $m^1(e')=e$ and $d_F(e')=n^0(v)$.
Since $n$ is a factor map,
there exists unique $e''\in G^1$ 
with $n^1(e'')=e'$ and $d_G(e'')=v$.
Therefore $e''\in G^1$ is the unique element satisfying
$m^1(n^1(e''))=e$ and $d_G(e'')=v$.
Thus $m\circ n$ is a factor map. 
By Proposition \ref{functT}, 
we have 
\begin{align*}
\omega\circ T_E^i
&=T_G^i\circ\omega^i
=T_G^i\circ\nu^i\circ\mu^i
=\nu\circ T_F^i\circ\mu^i
=\nu\circ\mu\circ T_E^i
\end{align*}
for $i=0,1$. 
Since $\cT(E)$ is generated by the images of $T_E^0$ and $T_E^1$, 
we have $\omega=\nu\circ\mu$. 
\end{proof}

The factor map $m\circ n$ defined in Proposition \ref{comp}
is called the {\em composition} of factor maps $m$ and $n$.
Thus we get a contravariant functor $E\mapsto\cT(E)$ 
from the category of topological graphs 
with factor maps as morphisms 
to the one of $C^*$-algebras with $*$-ho\-mo\-mor\-phisms as morphisms.
We study which factor maps from $F$ to $E$ give 
$*$-ho\-mo\-mor\-phisms from $\cO(E)$ to $\cO(F)$. 
For a factor map $m=(m^0,m^1)$ from $F$ to $E$, 
we can define a $*$-ho\-mo\-mor\-phism 
$\psi\colon \cK(C_{d_E}(E^1))\to\cK(C_{d_F}(F^1))$ 
by $\psi(\theta_{\xi,\eta})=\theta_{\mu^1(\xi),\mu^1(\eta)}$
for $\xi,\eta\in C_{d_E}(E^1)$ 
where $\mu^1\colon C_{d_E}(E^1)\to C_{d_F}(F^1)$ 
is defined from $m^1$ as above. 
To get a $*$-ho\-mo\-mor\-phism from $\cO(E)$ to $\cO(F)$, 
we need to know whether 
we have $\psi(\pi_{r_E}(f))=\pi_{r_F}(\mu^0(f))$ 
for $f\in C_0(\s{E}{0}{rg})$. 
As defined in Section \ref{SecPre}, 
the map $\pi_{r_E}\colon C_0(E^0)\to\cL(C_d(E^1))$ 
is the composition of the map 
$C_0(E^0)\ni f\mapsto f\circ r_E\in C_b(E^1)$ 
and $\pi_E\colon C_b(E^1)\to \cL(C_{d_E}(E^1))$. 
The map $\pi_{r_F}\colon C_0(F^0)\to\cL(C_d(F^1))$ 
is also the composition of the map 
$C_0(F^0)\ni f\mapsto f\circ r_F\in C_b(F^1)$ 
and $\pi_F\colon C_b(F^1)\to \cL(C_{d_F}(F^1))$. 
In order to get the equality 
$\psi(\pi_{r_E}(f))=\pi_{r_F}(\mu^0(f))$ 
for $f\in C_0(\s{E}{0}{rg})$, 
it suffices to see that 
$\psi(\pi_E(g))=\pi_F(\mu^1(g))$ for $g=f\circ r_E\in C_0(E^1)$
and that $\mu^1(g)=\mu^0(f)\circ r_F$. 
The former equality is valid for arbitrary factor maps. 

\begin{proposition}\label{psi}
Let $\psi\colon \cK(C_{d_E}(E^1))\to\cK(C_{d_F}(F^1))$ be 
the $*$-ho\-mo\-mor\-phism defined 
from a $*$-ho\-mo\-mor\-phism $\mu^1\colon C_{d_E}(E^1)\to C_{d_F}(F^1)$ 
as above. 
Then we have $\psi(\pi_E(g))=\pi_F(\mu^1(g))$ 
for all $g\in C_0(E^1)$. 
\end{proposition}

\begin{proof}
Since $\psi\circ\pi_E$ and $\pi_F\circ\mu^1$ are continuous, 
it suffices to show the equality 
only for elements in $C_c(E^1)$. 
Take $g\in C_c(E^1)$. 
By \cite[Lemma 1.16]{Ka1}, 
there exist $\xi_k,\eta_k\in C_c(E^1)$ for $k=1,\ldots,m$ 
such that $g=\sum_{k=1}^m\xi_k\overline{\eta_k}$
and that $\xi_k(e)\,\overline{\eta_k(e')}=0$ for any $k$ and any 
$e,e'\in E^1$ with $d_E(e)=d_E(e')$ and $e\neq e'$. 
By \cite[Lemma 1.15]{Ka1}, 
we have $\pi_E(g)=\sum_{k=1}^m\theta_{\xi_k,\eta_k}$. 
Thus we have 
$\psi(\pi_E(g))=\sum_{k=1}^m\theta_{\mu^1(\xi_k),\mu^1(\eta_k)}$. 
We also have $\mu^1(g)=\sum_{k=1}^m\mu^1(\xi_k)\overline{\mu^1(\eta_k)}$. 
To prove $\pi_F(\mu^1(g))=\sum_{k=1}^m\theta_{\mu^1(\xi_k),\mu^1(\eta_k)}$ 
by \cite[Lemma 1.15]{Ka1}, 
we need to check that $\mu^1(\xi_k)(e)\,\overline{\mu^1(\eta_k)(e')}=0$ 
for each $k$ and $e,e'\in F^1$ with $d_F(e)=d_F(e')$ and $e\neq e'$. 
If either $e$ or $e'$ is not in $(m^1)^{-1}(E^1)$, 
then clearly $\mu^1(\xi_k)(e)\,\overline{\mu^1(\eta_k)(e')}=0$.
When both $e$ and $e'$ are in $(m^1)^{-1}(E^1)$,
$d_F(e)=d_F(e')$ implies $d_E(m^1(e))=d_E(m^1(e'))$ 
by the condition (i) in Definition \ref{DefFactor}, 
and $e\neq e'$ implies $m^1(e)\neq m^1(e')$ 
by the condition (ii). 
Hence 
$$\mu^1(\xi_k)(e)\,\overline{\mu^1(\eta_k)(e')}
=\xi_k(m^1(e))\,\overline{\eta_k(m^1(e'))}=0.$$
Thus \cite[Lemma 1.15]{Ka1} implies 
$$\pi_F(\mu^1(g))=\sum_{k=1}^m\theta_{\mu^1(\xi_k),\mu^1(\eta_k)}
=\psi(\pi_E(g)).$$ 
We are done. 
\end{proof}

The equality $\mu^1(f\circ r_E)=\mu^0(f)\circ r_F$ 
for $f\in C_0(\s{E}{0}{rg})$ is not true 
for a general factor map $m=(m^0,m^1)$. 
We need the following notion. 

\begin{definition}\label{DefRegular}
A factor map $m=(m^0,m^1)$ from $F$ to $E$ is called {\em regular} 
if $(r_F)^{-1}(v)$ is non-empty and contained in $(m^1)^{-1}(E^1)$ 
for every $v\in F^0$ with $m^0(v)\in\s{E}{0}{rg}$. 
\end{definition}

\begin{lemma}\label{regular}
For a regular factor map $m$ from $F$ to $E$, 
we have $(m^0)^{-1}(\s{E}{0}{rg})\subset \s{F}{0}{rg}$. 
\end{lemma}

\begin{proof}
Let us take $v\in (m^0)^{-1}(\s{E}{0}{rg})\subset F^0$. 
Take a compact neighborhood $V$ of $m^0(v)\in \s{E}{0}{rg}$ 
such that $V\subset \s{E}{0}{rg}$, 
and set $U=(r_E)^{-1}(V)\subset E^1$. 
Then $U$ is compact and $r_E(U)=V$ 
(see Lemma \ref{E0r}). 
Set $V'=(m^0)^{-1}(V)\subset F^0$ and $U'=(r_F)^{-1}(V')\subset F^1$. 
Then $V'$ is a neighborhood of $v$. 
By the regularity of $m$, 
we have $r_F(U')=V'$ and $U'\subset (m^1)^{-1}(E^1)$. 
The condition (i) of Definition \ref{DefFactor} 
tells us that $U'=(m^1)^{-1}(U)$. 
Since $U$ is compact and $m^1$ is proper, 
we see that $U'$ is compact. 
Thus we have found a neighborhood $V'$ of $v$ 
such that $(r_F)^{-1}(V')=U'$ is compact 
and $r_F((r_F)^{-1}(V'))= r_F(U')=V'$. 
By Lemma \ref{E0r}, 
we see that $v\in\s{F}{0}{rg}$. 
Thus we have $(m^0)^{-1}(\s{E}{0}{rg})\subset \s{F}{0}{rg}$. 
\end{proof}

\begin{lemma}\label{funct0}
Let $m=(m^0,m^1)$ be a regular factor map from $F$ to $E$, 
and define a $*$-ho\-mo\-mor\-phism $\mu^0\colon C_0(E^0)\to C_0(F^0)$, 
a linear map $\mu^1\colon C_{d_E}(E^1)\to C_{d_F}(F^1)$ 
and a $*$-ho\-mo\-mor\-phism 
$\psi\colon \cK(C_{d_E}(E^1))\to\cK(C_{d_F}(F^1))$ as before. 
Then we have $\psi(\pi_{r_E}(f))=\pi_{r_F}(\mu^0(f))$ 
for $f\in C_0(\s{E}{0}{rg})$. 
\end{lemma}

\begin{proof}
By Proposition \ref{psi}, 
it suffices to show that $\mu^1(f\circ r_E)=\mu^0(f)\circ r_F$ 
for $f\in C_0(\s{E}{0}{rg})$. 
For $e\in (m^1)^{-1}(E^1)$, 
we have 
$$\mu^1(f\circ r_E)(e)=f(r_E(m^1(e)))=f(m^0(r_F(e)))=\mu^0(f)(r_F(e)).$$
For $e\notin (m^1)^{-1}(E^1)$, 
we have $\mu^1(f\circ r_E)(e)=0$.
If $r_F(e)\notin (m^0)^{-1}(E^0)$ then $\mu^0(f)(r_F(e))=0$.
If $r_F(e)\in (m^0)^{-1}(E^0)$ then we have $m^0(r_F(e))\notin\s{E}{0}{rg}$ 
by the regularity of the factor map $m$. 
Hence in this case, $\mu^0(f)(r_F(e))=f(m^0(r_F(e)))=0$.
Therefore we have $\mu^1(f\circ r_E)=\mu^0(f)\circ r_F$. 
We are done. 
\end{proof}

\begin{proposition}\label{funct1}
Let $\mu^0\colon C_0(E^0)\to C_0(F^0)$ 
and $\mu^1\colon C_{d_E}(E^1)\to C_{d_F}(F^1)$ be maps 
defined from a regular factor map $m$ from $F$ to $E$. 
Then there exists a unique $*$-ho\-mo\-mor\-phism $\mu\colon \cO(E)\to\cO(F)$ 
such that $\mu\circ t_E^i=t_F^i\circ\mu^i$ for $i=0,1$.

The $*$-ho\-mo\-mor\-phism $\mu$ is injective if and only if $m^0$ is surjective.
\end{proposition}

\begin{proof}
To define a $*$-ho\-mo\-mor\-phism $\mu\colon \cO(E)\to\cO(F)$ 
such that $\mu\circ t_E^i=t_F^i\circ\mu^i$ for $i=0,1$, 
it suffices to check that the pair of maps 
$T^0=t_F^0\circ\mu^0\colon C_0(E^0)\to \cO(F)$ and 
$T^1=t_F^1\circ\mu^1\colon C_{d_E}(E^1) \to \cO(F)$ 
is a Cuntz-Krieger $E$-pair.
We already saw that 
$T=(T^0,T^1)$ is a Toeplitz $E$-pair in Proposition \ref{functT}. 
The map $\varPhi\colon \cK(C_{d_E}(E^1))\to\cO(F)$ defined by 
$\varPhi(\theta_{\xi,\eta})=T^1(\xi)T^1(\eta)^*$ 
satisfies the equation $\varPhi=\varphi_F\circ\psi$. 
For $f\in C_0(\s{E}{0}{rg})$, 
we have $\mu^0(f)\in C_0(\s{F}{0}{rg})$ by Lemma \ref{regular}. 
Hence we have 
\begin{align*}
T^0(f)
&=t_F^0(\mu^0(f))
 =\varphi_F(\pi_{r_F}(\mu^0(f)))
 =\varphi_F(\psi(\pi_{r_E}(f)))
 =\varPhi(\pi_{r_E}(f)) 
\end{align*}
by Lemma \ref{funct0}. 
This implies that $T$ is a Cuntz-Krieger $E$-pair. 
Therefore there exists a $*$-ho\-mo\-mor\-phism $\mu\colon \cO(E)\to \cO(F)$ 
such that $\mu\circ t_E^i=t_F^i\circ\mu^i$ for $i=0,1$.
The uniqueness is easily verified.

The $C^*$-algebra $\cO(F)$ has the gauge action $\beta$ 
and we see that $\beta_z(T^0(f))=T^0(f)$ 
and $\beta_z(T^1(\xi))=zT^1(\xi)$ 
for $f\in C_0(E^0)$, $\xi\in C_d(E^1)$ and $z\in\T$.
Hence by Proposition \ref{GIUT}
the $*$-ho\-mo\-mor\-phism $\mu$ is injective 
if and only if $T^0=t_F^0\circ\mu^0$ is injective.
Since $t_F^0$ is injective, 
$t_F^0\circ\mu^0$ is injective
if and only if so is $\mu^0$.
It is easy to see that $\mu^0$ is injective 
exactly when $m^0$ is surjective.
Thus $\mu$ is injective if and only if $m^0$ is surjective.
\end{proof}

Note that if $m^0$ is surjective then so is $m^1$ 
by the condition (ii) in Definition \ref{DefFactor}. 

\begin{proposition}\label{compreg}
Let $E,F,G$ be topological graphs,
and $m,n$ be regular factor maps from $F$ to $E$ and $G$ to $F$ respectively.
Then the composition $m\circ n$ of $m$ and $n$ is regular 
and the $*$-ho\-mo\-mor\-phism $\cO(E)\to\cO(G)$ defined from $m\circ n$ 
is the composition of the two maps $\cO(E)\to\cO(F)$ and $\cO(F)\to\cO(G)$ 
defined from $m$ and $n$ respectively. 
\end{proposition}

\begin{proof}
Take $v\in G^0$ with $(m^0\circ n^0)(v)\in\s{E}{0}{rg}$. 
We have $v\in \s{G}{0}{rg}$ 
because Lemma \ref{regular} implies 
$$(m^0\circ n^0)^{-1}(\s{E}{0}{rg})\subset (n^0)^{-1}(\s{F}{0}{rg})\subset 
\s{G}{0}{rg}.$$
Hence the set $(r_G)^{-1}(v)$ is not empty. 
Take $e\in (r_G)^{-1}(v)$. 
Since $n^0(r_G(e))\in (m^0)^{-1}(\s{E}{0}{rg})\subset \s{F}{0}{rg}$, 
we have $n^1(e)\in F^1$. 
Since $m^0\big(r_F(n^1(e))\big)=m^0(n^0(r_G(e)))\in \s{E}{0}{rg}$, 
we have $m^1(n^1(e))\in E^1$. 
Thus $m\circ n$ is regular. 
The remainder of the proof is similar to that of Proposition \ref{comp}. 
\end{proof}

\section{$C^*$-algebras generated by Toeplitz pairs}\label{SecToepfam}

In this section, 
we investigate $C^*$-algebras $C^*(T)$ 
generated by Toeplitz pairs $T=(T^0,T^1)$. 
To this end, 
we introduce a construction of new topological graph $E_Y$ 
from an original topological graph $E$ and 
a closed subset $Y$ of $\s{E}{0}{rg}$, 
and see that Cuntz-Krieger $E_Y$-pairs 
are useful to study Toeplitz $E$-pairs. 
We will use the rsults in this section 
for analysing ideal structures of $\cO(E)$ in \cite{Ka2}.

\begin{definition}
Let $E$ be a topological graph. 
For a Toeplitz $E$-pair $T=(T^0,T^1)$, 
we define a closed subset $Y_T$ of $\s{E}{0}{rg}$ by 
$$C_0(\s{E}{0}{rg}\setminus Y_T)
=\{f\in C_0(\s{E}{0}{rg})\mid T^0(f)=\varPhi(\pi_r(f))\}.$$
\end{definition}

It is not difficult to see that 
the right hand side of the equation above 
is an ideal of $C_0(\s{E}{0}{rg})$. 
We have $Y_T=\emptyset$ if and only if 
$T$ is a Cuntz-Krieger $E$-pair. 
Thus $Y_T$ measures how far $T$ is 
from being a Cuntz-Krieger $E$-pair. 

\begin{lemma}\label{YT}
Let $E$ be a topological graph 
and $T$ be an injective Toeplitz $E$-pair. 
For $f\in C_0(E^0)$, 
we have $T^0(f)\in \varPhi\big(\cK(C_d(E^1))\big)$ 
if and only if $f\in C_0(\s{E}{0}{rg}\setminus Y_T)$. 
\end{lemma}

\begin{proof}
For $f\in C_0(\s{E}{0}{rg}\setminus Y_T)$, 
we have $T^0(f)=\varPhi(\pi_r(f))\in \varPhi\big(\cK(C_d(E^1))\big)$. 
Conversely take $f\in C_0(E^0)$ with 
$T^0(f)\in \varPhi\big(\cK(C_d(E^1))\big)$. 
By \cite[Proposition 2.11]{Ka1}, 
we have $f\in C_0(\s{E}{0}{rg})$ 
and $T^0(f)=\varPhi(\pi_r(f))$. 
Hence we get $f\in C_0(\s{E}{0}{rg}\setminus Y_T)$. 
\end{proof}

We will construct a topological graph $E_Y$ 
from a topological graph $E$ and 
a closed subset $Y$ of $\s{E}{0}{rg}$, 
and prove that there exists a one-to-one correspondence 
between Toeplitz $E$-pairs $T$ with $Y_T\subset Y$ 
and Cuntz-Krieger $E_Y$-pairs. 
This enables us to use results on Cuntz-Krieger pairs 
for analyzing $C^*$-algebras $C^*(T)$ 
generated by the Toeplitz pairs $T$. 
This fact also gives a new definition of $\cO(E)$ 
which does not use the space $\s{E}{0}{rg}$ 
or the notion of Cuntz-Krieger pairs (Proposition \ref{OEminimum}). 

Let us take a topological graph $E=(E^0,E^1,d,r)$ 
and a closed subset $Y$ of $\s{E}{0}{rg}$. 
We define a topological graph $E_Y=(E_Y^0,E_Y^1,d_Y,r_Y)$ as follows. 
Set $\partial Y=\overline{Y}\setminus Y$ 
where $\overline{Y}$ is taken in $E^0$. 
Since $Y$ is a closed subset of an open subset $\s{E}{0}{rg}$, 
we see that $Y$ is open in $\overline{Y}$. 
Hence $\partial Y$ is closed in $\overline{Y}$. 
A locally compact space 
$E_Y^0=E^0\mathop{\amalg}_{\partial Y}\overline{Y}$ 
is defined to be a topological space obtained from 
the disjoint union $E^0\amalg\overline{Y}$ 
by identifying the common closed subset $\partial Y$. 
Similarly we define 
$E_Y^1=E^1\mathop{\amalg}_{d^{-1}(\partial Y)} d^{-1}(\overline{Y})$. 
Note that we have $d^{-1}(\overline{Y})=\overline{d^{-1}(Y)}$ 
and $d^{-1}(\partial Y)=\partial (d^{-1}(Y))$ 
because $d\colon E^1\to E^0$ is locally homeomorphic. 
We consider $E^0$ and $E^1$ as subsets of $E_Y^0$ and $E_Y^1$, 
respectively. 
Both inclusions $Y\to\overline{Y}\subset E_Y^0$ and 
$d^{-1}(Y)\to d^{-1}(\overline{Y})\subset E_Y^1$ are 
denoted by $\omega$. 
Thus $\omega(Y)$ and $\omega(d^{-1}(Y))$ are 
the complements of the closed subsets $E^0\subset E_Y^0$ 
and $E^1\subset E_Y^1$, respectively. 
We may extend $d$ and $r$ to maps $d_Y$ and $r_Y$ on $E_Y^1$ 
by setting 
$$d_Y(\omega(e))=\omega(d(e))\in \omega(Y)\subset E_Y^0,
\quad\text{and}\quad
r_Y(\omega(e))=r(e)\in E^0\subset E_Y^0$$
for $e\in d^{-1}(Y)$. 
It is not difficult to see that $d_Y\colon E_Y^1\to E_Y^0$ is a local homeomorphism 
and $r_Y\colon E_Y^1\to E_Y^0$ is a continuous map. 
Thus we get a topological graph $E_Y=(E_Y^0,E_Y^1,d_Y,r_Y)$. 
Note that $E_Y$ is obtained from the topological graph $E$ 
by attaching extra vertices $\omega(Y)$ 
and extra edges $\omega(d^{-1}(Y))$ whose domains are in $\omega(Y)$ 
and ranges are in $E^0$. 

\begin{lemma}\label{F0a}
We have $(E_Y^0)_{\rs{fin}}=\s{E}{0}{fin}\cup\omega(Y)$, 
$(E_Y^0)_{\rs{sce}}=\s{E}{0}{sce}\cup\omega(Y)$,
and $(E_Y^0)_{\rs{rg}}=\s{E}{0}{rg}$.
\end{lemma}

\begin{proof}
For any $v\in Y$, 
we have $(r_Y)^{-1}(\omega(v))=\emptyset$.
Hence the open subset $\omega(Y)$ of $E_Y^0$ 
is contained in $(E_Y^0)_{\rs{sce}}$. 
It is easy to see that for $v\in E^0\setminus \partial Y$ 
we have 
$v\in(E_Y^0)_{\rs{fin}}$ if and only if $v\in\s{E}{0}{fin}$, 
and $v\in(E_Y^0)_{\rs{sce}}$ if and only if $v\in\s{E}{0}{sce}$.
Let us take $v\in \partial Y$. 
For a neighborhood $V$ of $v\in E^0$, 
the set $V'=V\cup \omega(V\cap Y)\subset E^0$ 
is a neighborhood of $v\in E_Y^0$, 
and we have $r_Y^{-1}(V')=r^{-1}(V)$.
Since we can find a neighborhood of the above form 
in every neighborhood of $v\in \partial Y\subset E_Y^0$, 
we see that 
$v\in(E_Y^0)_{\rs{fin}}$ if and only if $v\in\s{E}{0}{fin}$, 
and that $v\in(E_Y^0)_{\rs{sce}}$ if and only if $v\in\s{E}{0}{sce}$.
Therefore we get $(E_Y^0)_{\rs{fin}}=\s{E}{0}{fin}\cup\omega(Y)$ 
and $(E_Y^0)_{\rs{sce}}=\s{E}{0}{sce}\cup\omega(Y)$.
Finally we have 
$(E_Y^0)_{\rs{rg}}=(E_Y^0)_{\rs{fin}}\setminus\overline{(E_Y^0)_{\rs{sce}}}
=\s{E}{0}{rg}\setminus \partial Y$. 
Since $Y$ is a closed subset of $\s{E}{0}{rg}$, 
we have $\s{E}{0}{rg}\cap\overline{Y}=Y$. 
Hence $\s{E}{0}{rg}\cap \partial Y=\emptyset$. 
Thus we get 
$(E_Y^0)_{\rs{rg}}=\s{E}{0}{rg}\setminus \partial Y=\s{E}{0}{rg}$. 
\end{proof}

We define a map $m_Y^0\colon E_Y^0\to E^0$ and $m_Y^1\colon E_Y^1\to E^1$ by 
the identities on $E^0$ and $E^1$, 
and $m_Y^0(\omega(v))=v$, $m_Y^1(\omega(e))=e$ 
for $v\in V$ and $e\in d^{-1}(V)$. 
Both $m_Y^0$ and $m_Y^1$ are 
proper continuous surjections. 
Hence these extend the continuous maps 
$\widetilde{E}_Y^0\to \widetilde{E}^0$ 
and $\widetilde{E}_Y^1\to \widetilde{E}^1$, 
which are still denoted by $m_Y^0$ and $m_Y^1$. 
It is not difficult to see that 
the pair $m_Y=(m_Y^0,m_Y^1)$ 
is a factor map from $E_Y$ to $E$. 
The factor map $m_Y$ is not regular when $Y\neq\emptyset$ 
and $m_Y$ is identity when $Y=\emptyset$. 
We define a $*$-ho\-mo\-mor\-phism $\mu_Y^0\colon C_0(E^0)\to C_0(E_Y^0)$ 
and a linear map $\mu_Y^1\colon C_d(E^1)\to C_{d_Y}(E_Y^1)$ 
from $m_Y=(m_Y^0,m_Y^1)$. 
The $*$-ho\-mo\-mor\-phism $\mu_Y^0$ is an isomorphism onto 
the subalgebra 
$$\{h\in C_0(E_Y^0)\mid h(v)=h(\omega(v))\mbox{ for }v\in Y\},$$
and the linear map $\mu_Y^1$ is an isometric map onto 
$$\{\zeta\in C_{d_Y}(E_Y^1)\mid \zeta(e)=\zeta(\omega(e))
\mbox{ for }e\in d^{-1}(Y)\}.$$
Let $t_Y=(t_Y^0,t_Y^1)$ be the universal Cuntz-Krieger $E_Y$-pair 
on $\cO(E_Y)$ and define 
a $*$-ho\-mo\-mor\-phism $T_Y^0\colon C_0(E^0)\to \cO(E_Y)$ 
and a linear map $T_Y^1\colon C_d(E^1)\to \cO(E_Y)$ by
$T_Y^i=t_Y^i\circ\mu_Y^i$ for $i=0,1$. 
Let $\varphi_Y\colon \cK(C_{d_Y}(E_Y^1))\to \cO(E_Y)$ 
and $\varPhi_Y\colon \cK(C_d(E^1))\to \cO(E_Y)$ be 
$*$-ho\-mo\-mor\-phisms determined by $t_Y$ and $T_Y$, 
respectively. 
Note that we have $\varPhi_Y=\varphi_Y\circ\psi$ 
where $\psi\colon \cK(C_d(E^1))\to \cK(C_{d_Y}(E_Y^1))$ 
is defined 
by $\psi(\theta_{\xi,\eta})=\theta_{\mu_Y^1(\xi),\mu_Y^1(\eta)}$ 
for $\xi,\eta\in C_d(E^1)$. 

\begin{lemma}\label{lemEY}
For an element $f\in C_0(\s{E}{0}{rg})$, 
we have $\pi_{r_Y}(\mu^0(f))=\psi(\pi_r(f))\in\cK(C_{d_Y}(E_Y^1))$. 
\end{lemma}

\begin{proof}
Take $f\in C_0(\s{E}{0}{rg})$. 
By Proposition \ref{psi}, 
it suffices to see that $\mu^0(f)\circ r_Y=\mu^1(f\circ r)$. 
This is clear because we have $m^0(r_Y(e))=r(m^1(e))$ 
for all $e\in E^1$. 
The proof is completed. 
\end{proof}

\begin{proposition}\label{EY}
The pair $T_Y=(T_Y^0,T_Y^1)$ 
is an injective Toeplitz $E$-pair in $\cO(E_Y)$ 
such that $Y_{T_Y}=Y$. 
\end{proposition}

\begin{proof}
By Proposition \ref{functT}, 
$T_Y$ is a Toeplitz $E$-pair. 
Clearly the pair $T_Y$ is injective. 
We will show that $Y_{T_Y}=Y$. 
Take $f\in C_0(\s{E}{0}{rg}\setminus Y)$. 
We have $\mu_Y^0(f)\in C_0((E_Y^0)_{\rs{rg}})$ 
by Lemma \ref{F0a}. 
We have 
$$T_Y^0(f)=t_Y^0(\mu_Y^0(f))
=\varphi_Y\big(\pi_{r_Y}(\mu_Y^0(f))\big)
=\varphi_Y\big(\psi(\pi_r(f))\big)
=\varPhi_Y(\pi_r(f))$$
by Lemma \ref{lemEY}. 
Conversely take $f\in C_0(\s{E}{0}{rg})$ 
with $T_Y^0(f)=\varPhi_Y(\pi_r(f))$. 
We see 
$$t_Y^0(\mu_Y^0(f))=
T_Y^0(f)=\varPhi_Y(\pi_r(f))\in \varPhi_Y\big(\cK(C_d(E^1))\big)
\subset \varphi_Y\big(\cK(C_d(E_Y^1))\big).$$
Hence we get $\mu_Y^0(f)\in C_0((E_Y^0)_{\rs{rg}})$ 
by Lemma \ref{YT}. 
This implies that $f\in C_0(\s{E}{0}{rg}\setminus Y)$ 
by Lemma \ref{F0a}. 
Thus we have 
$$\{f\in C_0(\s{E}{0}{rg})\mid T_Y^0(f)=\varPhi_Y(\pi_r(f))\}
=C_0(\s{E}{0}{rg}\setminus Y).$$
This means that $Y_{T_Y}=Y$. 
\end{proof}

By Proposition \ref{EY}, 
for any $*$-ho\-mo\-mor\-phism $\rho\colon \cO(E_Y)\to B$, 
the pair $T=(T^0,T^1)$ defined by $T^i=\rho\circ T_Y^i$ for $i=0,1$ 
is a Toeplitz $E$-pair in the $C^*$-algebra $B$ 
satisfying $Y_T\subset Y$. 
We will prove the converse. 
Take a topological graph $E$, 
a closed subset $Y$ of $\s{E}{0}{rg}$ 
and a Toeplitz $E$-pair $T$ with $Y_T\subset Y$. 
To get a $*$-ho\-mo\-mor\-phism $\rho\colon \cO(E_Y)\to C^*(T)$ 
such that $T^i=\rho\circ T_Y^i$ for $i=0,1$, 
it suffices to construct a Cuntz-Krieger $E_Y$-pair 
$\tT=(\tT^0,\tT^1)$ on $C^*(T)$ 
such that $T^i=\tT^i\circ \mu_Y^i$ for $i=0,1$. 

Define a $*$-ho\-mo\-mor\-phism 
$T^0_{\rs{rg}}\colon C_0(\s{E}{0}{rg})\to C^*(T)$
by $T^0_{\rs{rg}}=\varPhi\circ\pi_r$. 
The condition $Y_T\subset Y$ implies that 
$T^0(f)=T^0_{\rs{rg}}(f)$ 
for $f\in C_0(\s{E}{0}{rg}\setminus Y)
\subset C_0(\s{E}{0}{rg}\setminus Y_T)$.
For $f\in C_0(E^0)$ and $g\in C_0(\s{E}{0}{rg})$,
we have 
$$T^0(f)T_{\rs{rg}}^0(g)=T^0(f)\varPhi(\pi_r(g))
=\varPhi(\pi_r(f)\pi_r(g))
=T_{\rs{rg}}^0(fg).$$
We also have $T_{\rs{rg}}^0(g)T^0(f)=T_{\rs{rg}}^0(gf)$.

\begin{lemma}\label{DeftT0}
For $h\in C_0(E_Y^0)$, 
there exist $f\in C_0(E^0)$ 
and $g\in C_0(\s{E}{0}{rg})\subset C_0(E^0)$ 
such that 
$h(v)=f(v)+g(v)$ for $v\in E^0$ 
and $h(\omega(v))=f(v)$ for $v\in Y$. 
The element $T^0(f)+T^0_{\rs{rg}}(g)\in C^*(T)$ 
does not depend on the choices of $f$ and $g$ 
satisfying the above two equations. 
\end{lemma}

\begin{proof}
Take $h\in C_0(E_Y^0)$. 
Define a function $g_0$ on $Y$ 
by $g_0(v)=h(v)-h(\omega(v))$ for $v\in Y$. 
For a net $\{v_i\}$ in $Y$ converges 
to an element in $\partial Y$, 
we have $\lim g_0(v_i)=\lim \big(h(v_i)-h(\omega(v_i))\big)=0$. 
Hence we get $g_0\in C_0(Y)$. 
Since $Y$ is closed in $\s{E}{0}{rg}$, 
$g_0$ extends to a function 
$g\in C_0(\s{E}{0}{rg})$ 
such that $g(v)=h(v)-h(\omega(v))$ for $v\in Y$. 
Define $f\in C_0(E^0)$ by $f(v)=h(v)-g(v)$ for $v\in E^0$. 
Now it is easy to see that $f$ and $g$ satisfy the equations 
$h(v)=f(v)+g(v)$ for $v\in E^0$ 
and $h(\omega(v))=f(v)$ for $v\in Y$. 
Let us take other $f'\in C_0(E^0)$ and $g'\in C_0(\s{E}{0}{rg})$ 
satisfying the two conditions. 
Then we have $g-g'=-(f-f')\in C_0(\s{E}{0}{rg}\setminus Y)$. 
Hence we see 
\begin{align*}
\big(T^0(f)+T^0_{\rs{rg}}(g)\big)
-\big(T^0(f')+T^0_{\rs{rg}}(g')\big)
&=T^0(f-f')+T^0_{\rs{rg}}(g-g')\\
&=T^0(f-f')-T^0_{\rs{rg}}(f-f')=0.
\end{align*}
Thus the element $T^0(f)+T^0_{\rs{rg}}(g)$ 
does not depend on the choices of $f$ and $g$ 
satisfying the two conditions. 
\end{proof}

For $h\in C_0(E_Y^0)$, 
we define $\tT^0(h)\in C^*(T)$ 
by $\tT^0(h)=T^0(f)+T^0_{\rs{rg}}(g)$ 
where $f\in C_0(E^0)$ and $g\in C_0(\s{E}{0}{rg})$ 
are elements satisfying the two equations in Lemma \ref{DeftT0}. 

\begin{proposition}
The map $\tT^0\colon C_0(E_Y^0)\to C^*(T)$ 
is a $*$-ho\-mo\-mor\-phism satisfying 
the equation $T^0=\tT^0\circ \mu_Y^0$. 
\end{proposition}

\begin{proof}
It is clear that $\tT^0$ is linear and $*$-preserving. 
We will show that it is multiplicative. 
Take $h_1,h_2\in C_0(E_Y^0)$. 
For $i=1,2$, 
choose $f_i\in C_0(E^0)$ 
and $g_i\in C_0(\s{E}{0}{rg})$ 
satisfying the two conditions in Lemma \ref{DeftT0} 
for $h_i\in C_0(E_Y^0)$.
Set $f=f_1f_2\in C_0(E^0)$ and 
$$g=(f_1+g_1)(f_2+g_2)-f_1f_2
=f_1g_2+g_1f_2+g_1g_2\in C_0(\s{E}{0}{rg}).$$ 
Then we have $(h_1h_2)(v)=f(v)+g(v)$ for $v\in E^0$ 
and $(h_1h_2)(\omega(v))=f(v)$ for $v\in Y$. 
Hence we get $\tT^0(h_1h_2)=T^0(f)+T^0_{\rs{rg}}(g)$. 
On the other hand, 
we have 
\begin{align*}
\tT^0(h_1)\tT^0(h_2)&=\big(T^0(f_1)+T^0_{\rs{rg}}(g_1)\big)
\big(T^0(f_2)+T^0_{\rs{rg}}(g_2)\big)\\
&=T^0(f_1f_2)+T^0_{\rs{rg}}(f_1g_2+g_1f_2+g_1g_2)\\
&=T^0(f)+T^0_{\rs{rg}}(g)
\end{align*}
Thus $\tT^0$ is multiplicative. 

Take $f_0\in C_0(E^0)$ and set $h=\mu_Y^0(f_0)\in C_0(E_Y^0)$. 
We can take $f=f_0$ and $g=0$ 
in the definition of $\tT^0(h)$. 
Hence we get $\tT^0(h)=T^0(f_0)$. 
This proves $T^0=\tT^0\circ \mu_Y^0$. 
\end{proof}

\begin{proposition}\label{tTinj}
The $*$-ho\-mo\-mor\-phism $\tT^0\colon C_0(E_Y^0)\to C^*(T)$ 
is injective if and only if 
$T$ is an injective Toeplitz $E$-pair such that $Y_T=Y$. 
\end{proposition}

\begin{proof}
Suppose that $\tT^0$ is injective. 
Since $T^0=\tT^0\circ \mu_Y^0$, 
the map $T^0$ is injective. 
For $g_0\in C_0(\s{E}{0}{rg})$ with $T^0(g_0)=T^0_{\rs{rg}}(g_0)$, 
define $h\in C_0(\omega(Y))\subset C_0(E_Y^0)$ 
by $h(\omega(v))=g_0(v)$. 
We can take $f=g_0$ and $g=-g_0$ 
in the definition of $\tT^0(h)$. 
Hence we have $\tT^0(h)=T^0(g_0)-T^0_{\rs{rg}}(g_0)=0$, 
and so $h=0$.
This implies that $g_0\in C_0(\s{E}{0}{rg}\setminus Y)$. 
Thus we get $Y_T\supset Y$. 
Hence we have shown $Y_T=Y$ 
because the other inclusion is assumed. 

Conversely take an injective Toeplitz $E$-pair $T$ such that $Y_T=Y$. 
Take $h\in C_0(E_Y^0)$ with $\tT^0(h)=0$. 
Choose $f\in C_0(E^0)$ and $g\in C_0(\s{E}{0}{rg})$ 
such that $h(v)=f(v)+g(v)$ for $v\in E^0$ 
and $h(\omega(v))=f(v)$ for $v\in Y$. 
The condition $\tT^0(h)=T^0(f)+T^0_{\rs{rg}}(g)=0$ 
implies that $T^0(f)=T^0_{\rs{rg}}(-g)\in\varPhi\big(\cK(C_d(E^1))\big)$. 
Hence by Lemma \ref{YT}, 
we have $f\in C_0(\s{E}{0}{rg}\setminus Y_T)$ 
and $T^0(f)=T^0_{\rs{rg}}(f)$. 
Since $T^0$ is injective, 
so is $T^0_{\rs{rg}}$. 
Hence we have $f=-g$. 
Therefore we see $h(v)=f(v)+g(v)=0$ for $v\in E^0$ 
and $h(\omega(v))=f(v)=0$ for $v\in Y=Y_T$. 
Thus we have $h=0$. 
This shows that $\tT^0$ is injective. 
We are done. 
\end{proof}

Next, 
we define a linear map $\tT^1\colon C_{d_Y}(E_Y^1)\to C^*(T)$ 
such that $T^1=\tT^1\circ \mu_Y^1$. 
Recall that the open subset $\s{E}{1}{rg}$ of $E^1$ is defined 
by $\s{E}{1}{rg}=d^{-1}(\s{E}{0}{rg})$ 
and we showed that 
$$C_d(\s{E}{1}{rg})=\{\xi g\in C_d(E^1)\mid 
\xi\in C_d(E^1), g\in C_0(\s{E}{0}{rg})\}$$ 
in \cite[Lemma 1.12]{Ka1}. 

\begin{lemma}\label{Cd(F1)1}
The map $T_{\rs{rg}}^1\colon C_d(\s{E}{1}{rg})\to C^*(T)$ 
defined by $T_{\rs{rg}}^1(\xi g)=T^1(\xi)T_{\rs{rg}}^0(g)$ 
is a well-defined linear map satisfying that 
$$T_{\rs{rg}}^1(\eta_1)^*T_{\rs{rg}}^1(\eta_2)
=T_{\rs{rg}}^0(\ip{\eta_1}{\eta_2})$$
for $\eta_1,\eta_2\in C_d(\s{E}{1}{rg})$. 
\end{lemma}

\begin{proof}
For $i=1,2$, 
take $\eta_i=\xi_i g_i\in C_d(\s{E}{1}{rg})$ with 
$\xi_i\in C_d(E^1)$ and $g_i\in C_0(\s{E}{0}{rg})$. 
We see that 
\begin{align*}
\big(T^1(\xi_1)T_{\rs{rg}}^0(g_1)\big)^*
\big(T^1(\xi_2)T_{\rs{rg}}^0(g_2)\big)
&=T_{\rs{rg}}^0(\overline{g_1})T^0(\ip{\xi_1}{\xi_2})T_{\rs{rg}}^0(g_2)\\
&=T_{\rs{rg}}^0\big(\overline{g_1}\ip{\xi_1}{\xi_2}g_2\big)\\
&=T_{\rs{rg}}^0\big(\ip{\eta_1}{\eta_2}\big).
\end{align*}
This proves the last equality. 
For $\eta\in C_d(\s{E}{1}{rg})$, 
take $\xi_i\in C_d(E^1)$ and $g_i\in C_0(\s{E}{0}{rg})$
with $\eta=\xi_i g_i$ for $i=1,2$. 
By the computation above, 
we have 
$$\big(T^1(\xi_i)T_{\rs{rg}}^0(g_i)\big)^*
\big(T^1(\xi_j)T_{\rs{rg}}^0(g_j)\big)
=T_{\rs{rg}}^0\big(\ip{\eta}{\eta}\big)$$
for $i,j=1,2$. 
Hence we get $x^*x=0$ 
for 
$x=T^1(\xi_1)T_{\rs{rg}}^0(g_1)-T^1(\xi_2)T_{\rs{rg}}^0(g_2)\in C^*(T)$. 
Thus we have $T^1(\xi_1)T_{\rs{rg}}^0(g_1)=T^1(\xi_2)T_{\rs{rg}}^0(g_2)$. 
This proves the well-definedness of $T_{\rs{rg}}^1$. 
For $\eta_1,\eta_2\in C_d(\s{E}{1}{rg})$, 
we have 
$$\big(T_{\rs{rg}}^1(\eta_1+\eta_2)
-T_{\rs{rg}}^1(\eta_1)-T_{\rs{rg}}^1(\eta_2)\big)^*
\big(T_{\rs{rg}}^1(\eta_1+\eta_2)
-T_{\rs{rg}}^1(\eta_1)-T_{\rs{rg}}^1(\eta_2)\big)=0$$
by the computation above. 
This proves the linearity of $T_{\rs{rg}}^1$. 
\end{proof}

\begin{lemma}\label{T1=T1rg}
For $\eta\in C_d(\s{E}{1}{rg}\setminus d^{-1}(Y))\subset C_d(\s{E}{1}{rg})$, 
we have $T^1(\eta)=T_{\rs{rg}}^1(\eta)$. 
\end{lemma}

\begin{proof}
For $\eta\in C_d(\s{E}{1}{rg}\setminus d^{-1}(Y))$, 
we can find $\xi\in C_d(E^1)$ 
and $g\in C_0(\s{E}{0}{rg}\setminus Y)$ 
such that $\eta=\xi g$ by \cite[Lemma 1.12]{Ka1}. 
Hence we have 
$$T_{\rs{rg}}^1(\eta)=T^1(\xi)T_{\rs{rg}}^0(g)
=T^1(\xi)T^0(g)=T^1(\xi g)=T^1(\eta).$$
We are done. 
\end{proof}

By using Lemma \ref{T1=T1rg}, 
we can prove the following 
a proof of the following lemma 
may be given that is similar to the proof of Lemma \ref{DeftT0}. 

\begin{lemma}\label{DeftT1}
For $\zeta\in C_{d_Y}(E_Y^1)$, 
there exist $\xi\in C_d(E^1)$ 
and $\eta\in C_d(\s{E}{1}{rg})\subset C_d(E^1)$ 
such that 
$\zeta(e)=\xi(e)+\eta(e)$ for $e\in E^1$ 
and $\zeta(\omega(e))=\xi(e)$ for $e\in d^{-1}(Y)$. 
The element $T^1(\xi)+T^1_{\rs{rg}}(\eta)\in C^*(T)$ 
does not depend on the choices of $\xi$ and $\eta$ 
satisfying the above two equations. 
\end{lemma}

Hence we can define a linear map $\tT^1\colon C_{d_Y}(E_Y^1)\to C^*(T)$ 
by $\tT^1(\zeta)=T^1(\xi)+T^1_{\rs{rg}}(\eta)$ 
for $\zeta\in C_{d_Y}(E_Y^1)$, 
where $\xi\in C_d(E^1)$ and $\eta\in C_d(\s{E}{1}{rg})$ 
satisfy the two conditions in Lemma \ref{DeftT1}. 
It is easy to see that $T^1=\tT^1\circ \mu_Y^1$. 
We will prove that the pair $\tT=(\tT^0,\tT^1)$ 
is a Cuntz-Krieger $E_Y$-pair. 

\begin{lemma}\label{Cd(F1)1.5}
For $f\in C_0(E^0)$, $g\in C_0(\s{E}{0}{rg})$, 
$\xi\in C_d(E^1)$ and $\eta\in C_d(\s{E}{1}{rg})$, 
we have 
\begin{align*}
T^0(f)T_{\rs{rg}}^1(\eta)&=T_{\rs{rg}}^1(\pi_r(f)\eta),\quad 
T_{\rs{rg}}^0(g)T_{\rs{rg}}^1(\eta)=T_{\rs{rg}}^1(\pi_r(g)\eta),\\
\mbox{and }&\quad T^1(\xi)^*T_{\rs{rg}}^1(\eta)=T_{\rs{rg}}^0(\ip{\xi}{\eta}).
\end{align*}
\end{lemma}

\begin{proof}
Straightforward. 
\end{proof}

\begin{lemma}\label{Cd(F1)2}
We have $\tT^1(\zeta_1)^*\tT^1(\zeta_2)=\tT^0(\ip{\zeta_1}{\zeta_2})$ 
for $\zeta_1,\zeta_2\in C_{d_Y}(E_Y^1)$. 
\end{lemma}

\begin{proof}
For $i=1,2$, 
take $\xi_i\in C_d(E^1)$ and $\eta_i\in C_d(\s{E}{1}{rg})$ 
satisfying that $\zeta_i(e)=\xi_i(e)+\eta_i(e)$ for $e\in E^1$ 
and $\zeta_i(\omega(e))=\xi_i(e)$ for $e\in d^{-1}(Y)$. 
Set $f=\ip{\xi_1}{\xi_2}\in C_0(E^0)$ and 
$$g=\ip{\xi_1+\eta_1}{\xi_2+\eta_2}-\ip{\xi_1}{\xi_2}
=\ip{\eta_1}{\xi_2}+\ip{\xi_1}{\eta_2}+\ip{\eta_1}{\eta_2}
\in C_0(\s{E}{0}{rg}).$$
For $v\in E^0$, we have 
\begin{align*}
\ip{\zeta_1}{\zeta_2}(v)
&=\sum_{e\in (d_Y)^{-1}(v)}\overline{\zeta_1(e)}\zeta_2(e)
 =\sum_{e\in d^{-1}(v)}\overline{(\xi_1+\eta_1)(e)}(\xi_2+\eta_2)(e)\\
&=\ip{\xi_1+\eta_1}{\xi_2+\eta_2}(v)
 =f(v)+g(v). 
\end{align*}
Similarly for $v\in Y$, we have 
\begin{align*}
\ip{\zeta_1}{\zeta_2}(\omega(v))
&=\sum_{e\in (d_Y)^{-1}(\omega(v))}
  \overline{\zeta_1(e)}\zeta_2(e)
 =\sum_{e\in d^{-1}(v)\subset d^{-1}(Y)}
  \overline{\zeta_1(\omega(e))}\zeta_2(\omega(e))\\
&=\sum_{e\in d^{-1}(v)}\overline{\xi_1(e)}\xi_2(e)
 =\ip{\xi_1}{\xi_2}(v)
 =f(v).
\end{align*}
Hence we have $\tT^0(\ip{\zeta_1}{\zeta_2})=T^0(f)+T_{\rs{rg}}^0(g)$. 
By Lemma \ref{Cd(F1)1} and Lemma \ref{Cd(F1)1.5}, 
we have 
\begin{align*}
\tT^1(\zeta_1)^*\tT^1(\zeta_2)
&=\big(T^1(\xi_1)+T_{\rs{rg}}^1(\eta_1)\big)^*
  \big(T^1(\xi_2)+T_{\rs{rg}}^1(\eta_2)\big)\\
&=T^0(\ip{\xi_1}{\xi_2})
  +T_{\rs{rg}}^0(\ip{\eta_1}{\xi_2})
  +T_{\rs{rg}}^0(\ip{\xi_1}{\eta_2})
  +T_{\rs{rg}}^0(\ip{\eta_1}{\eta_2})\\
&=T^0(f)+T_{\rs{rg}}^0(g)\\
&=\tT^0(\ip{\zeta_1}{\zeta_2}).
\end{align*}
\end{proof}

\begin{lemma}\label{Cd(F1)3}
We have $\tT^0(h)\tT^1(\zeta)=\tT^1(\pi_{r_Y}(h)\zeta)$ 
for $h\in C_0(E_Y^0)$ and $\zeta\in C_{d_T}(E_Y^1)$. 
\end{lemma}

\begin{proof}
Take $f\in C_0(E^0)$ and $g\in C_0(\s{E}{0}{rg})$ 
such that $h(v)=f(v)+g(v)$ for $v\in E^0$ 
and $h(\omega(v))=f(v)$ for $v\in Y$. 
Take $\xi\in C_d(E^1)$ and $\eta\in C_d(\s{E}{1}{rg})$ 
such that $\zeta(e)=\xi(e)+\eta(e)$ for $e\in E^1$ 
and $\zeta(\omega(e))=\xi(e)$ for $e\in d^{-1}(Y)$. 
We have $T_{\rs{rg}}^0(g)T^1(\xi)=T^1(\pi_r(g)\xi)$. 
From this fact and Lemma \ref{Cd(F1)1.5}, 
we see that 
\begin{align*}
\tT^0(h)\tT^1(\zeta)
&=\big(T^0(f)+T_{\rs{rg}}^0(g)\big)
  \big(T^1(\xi)+T_{\rs{rg}}^1(\eta)\big)\\
&=T^1(\pi_r(f+g)\xi)+T_{\rs{rg}}^1(\pi_r(f+g)\eta)
\end{align*}
Let us set $\xi'=\pi_r(f+g)\xi\in C_d(E^1)$ 
and $\eta'=\pi_r(f+g)\eta\in C_d(\s{E}{1}{rg})$. 
For $e\in E^1$, 
we have 
$$(\pi_{r_Y}(h)\zeta)(e)=h(r(e))\zeta(e)
=\big(f(r(e))+g(r(e))\big)\big(\xi(e)+\eta(e)\big)
=\xi'(e)+\eta'(e).$$ 
For $e\in d^{-1}(Y)$, 
we have 
$$(\pi_{r_Y}(h)\zeta)(\omega(e))=h(r(e))\zeta(\omega(e))
=\big(f(r(e))+g(r(e))\big)\xi(e)
=\xi'(e).$$ 
Hence we get 
$$\tT^1(\pi_{r_Y}(h)\zeta)=T^1(\xi')+T_{\rs{rg}}^1(\eta')
=\tT^0(h)\tT^1(\zeta).$$ 
The proof is completed. 
\end{proof}

\begin{proposition}\label{CK}
The pair $\tT=(\tT^0,\tT^1)$ 
is a Cuntz-Krieger $E_Y$-pair 
satisfying the equation $T^i=\tT^i\circ \mu_Y^i$ for $i=0,1$ 
and $C^*(\tT)=C^*(T)$. 
\end{proposition}

\begin{proof}
By Lemma \ref{Cd(F1)2} and Lemma \ref{Cd(F1)3},  
$\tT$ is a Toeplitz $E_Y$-pair.
We will show that $\tT$ is a Cuntz-Krieger $E_Y$-pair.
To do so, 
it suffices to see that 
$\tT^0\big(C_0((E_Y^0)_{\rs{rg}})\big)\subset 
\widetilde{\varPhi}\big(\cK(C_{d_Y}(E_Y^1))\big)$ 
by Lemma \ref{YT}. 
Take $h\in C_0((E_Y^0)_{\rs{rg}})$. 
By Lemma \ref{F0a}, 
we have $(E_Y^0)_{\rs{rg}}=\s{E}{0}{rg}$. 
Hence there exists $g\in C_0(\s{E}{0}{rg})$ 
such that $h(v)=g(v)$ for $v\in E^0$. 
Thus we have $\tT^0(h)=T^0(0)+T^0_{\rs{rg}}(g)=\varPhi(\pi_r(g))$. 
Since $T^1(C_d(E^1))\subset \tT^1(C_{d_Y}(E_Y^1))$,
we have $\varPhi\big(\cK(C_d(E^1))\big)
\subset \widetilde{\varPhi}\big(\cK(C_{d_Y}(E_Y^1))\big)$.
Hence we get $\tT^0(h)\in\widetilde{\varPhi}
\big(\cK(C_{d_Y}(E_Y^1))\big)$ 
for every $h\in C_0((E_Y^0)_{\rs{rg}})$. 
Thus $\tT$ is a Cuntz-Krieger $E_Y$-pair. 
As we have already seen, 
the two equations $T^i=\tT^i\circ \mu_Y^i$ hold for $i=0,1$. 
This implies $\tT^0(C_0(E_Y^0))\supset T^0(C_0(E^0))$ and 
$\tT^1(C_{d_Y}(E_Y^1))\supset T^1(C_d(E^1))$. 
Hence we have $C^*(\tT)=C^*(T)$. 
\end{proof}

By Proposition \ref{CK}, 
we have a surjective $*$-ho\-mo\-mor\-phism $\rho_T\colon \cO(E_Y)\to C^*(T)$ 
such that $\tT^i=\rho_T\circ t_Y^i$ for $i=0,1$. 
We have 
$$\rho_T\circ T_Y^i=\rho_T\circ t_Y^i\circ \mu^i=\tT^i\circ \mu^i=T^i,$$
for $i=0,1$. 
We study for which Toeplitz pairs $T$ 
the surjections $\rho_T$ are injective. 

A {\em loop} is an element $e\in E^n$ with $n\geq 1$ 
such that $d(e)=r(e)$, 
and the vertex $d(e)=r(e)$ is called 
the {\em base point} of the loop $e$. 
A loop $e=(e_1,\ldots,e_n)$ is 
said to be {\em without entrances} 
if $r^{-1}(r(e_k))=\{e_{k}\}$ for $k=1,\ldots,n$.
Recall that a topological graph $E$ 
is said to be {\em topologically free} 
if the set of base points of loops without entrances 
has an empty interior (\cite[Definition 5.4]{Ka1}). 

\begin{proposition}\label{rhoisom2}
When $E$ is topologically free, 
the $*$-ho\-mo\-mor\-phism $\rho_T$ is an isomorphism 
if and only if the Toeplitz $E$-pair $T$ is injective 
and satisfies $Y_T=Y$. 
\end{proposition}

\begin{proof}
When $E$ is topologically free, 
the topological graph $E_Y$ is also topologically free 
because vertices in $\omega(Y)$ receive no edges. 
Hence $\rho_T$ is an isomorphism 
if and only if $\tT$ is injective by \cite[Theorem 5.12]{Ka1}. 
By Proposition \ref{tTinj}, 
$\tT$ is injective whenever 
$T$ is an injective Toeplitz $E$-pair such that $Y_T=Y$. 
This completes the proof. 
\end{proof}

\begin{lemma}\label{tTgauge}
The Cuntz-Krieger $E_Y$-pair $\tT$ admits a gauge action 
if and only if $T$ does. 
\end{lemma}

\begin{proof}
Since $\tT^0(C_0(E_Y^0))\supset T^0(C_0(E^0))$ and 
$\tT^1(C_{d_Y}(E_Y^1))\supset T^1(C_d(E^1))$, 
it is clear that if $\tT$ admits a gauge action 
then $T$ also does. 
Conversely, suppose that for each $z\in\T$, 
there exists an automorphism $\beta_z$ 
on $C^*(T)$ such that $\beta_z(T^0(f))=T^0(f)$ 
and $\beta_z(T^1(\xi))=zT^1(\xi)$ 
for $f\in C_0(E^0)$ and $\xi\in C_d(E^1)$. 
Then it is easy to see that $\beta_z(T_{\rs{rg}}^0(g))=T_{\rs{rg}}^0(g)$ 
and $\beta_z(T_{\rs{rg}}^1(\eta))=zT_{\rs{rg}}^1(\eta)$ 
for $g\in C_0(\s{E}{0}{rg})$ and $\eta\in C_d(\s{E}{1}{rg})$. 
Hence $\beta_z(\tT^0(h))=\tT^0(h)$ 
and $\beta_z(\tT^1(\zeta))=z\tT^1(\zeta)$ 
for $h\in C_0(E_Y^0)$ and $\zeta\in C_{d_Y}(E_Y^1)$. 
Thus the pair $\tT$ admits a gauge action. 
We are done. 
\end{proof}

\begin{proposition}\label{rhoisom}
The $*$-ho\-mo\-mor\-phism $\rho_T$ is an isomorphism 
if and only if the Toeplitz $E$-pair $T$ is injective, 
admits a gauge action, and satisfies $Y_T=Y$. 
\end{proposition}

\begin{proof}
This follows from Proposition \ref{GIUT} 
with the help of Lemma \ref{tTgauge} 
in a similar way to the proof of Proposition \ref{rhoisom2}. 
\end{proof}

\begin{corollary}\label{CT=OE}
Let $E$ be a topological graph. 
For an injective Toeplitz $E$-pair $T$ which admits a gauge action, 
the $C^*$-algebra $C^*(T)$ is isomorphic to $\cO(E_{Y_T})$. 
\end{corollary}

\begin{proposition}\label{CTY=OEY}
The Toeplitz $E$-pair $T_Y$ on $\cO(E_Y)$ 
in Proposition \ref{EY} satisfies $C^*(T_Y)=\cO(E_Y)$. 
\end{proposition}

\begin{proof}
Since we have $Y_{T_Y}=Y$ by Proposition \ref{EY}, 
we get a $*$-ho\-mo\-mor\-phism $\rho_{T_Y}\colon \cO(E_Y)\to C^*(T_Y)$ 
satisfying the equation $\rho_T\circ T_Y^i=T_Y^i$ for $i=0,1$. 
We have $\rho_{T_Y}(x)=x$ for $x\in C^*(T_Y)\subset \cO(E_Y)$. 
By Proposition \ref{rhoisom}, 
$\rho_{T_Y}$ is an isomorphism. 
Thus we have $C^*(T_Y)=\cO(E_Y)$. 
\end{proof}

Now we have the following proposition 
which implies that 
the $C^*$-algebra $\cO(E_Y)$ is the universal $C^*$-algebra 
generated by Toeplitz $E$-pairs $T$ satisfying $Y_T\subset Y$. 

\begin{proposition}\label{cOEY}
Let $E$ be a topological graph, 
and $Y$ be a closed subset of $\s{E}{0}{rg}$. 
For a Toeplitz $E$-pair $T$ with $Y_T\subset Y$, 
there exists a unique surjective $*$-ho\-mo\-mor\-phism 
$\rho_T\colon \cO(E_Y)\to C^*(T)$ 
satisfying $T^i=\rho_T\circ T_Y^i$ for $i=0,1$. 
\end{proposition}

\begin{proof}
We have already seen that there exists such a surjection $\rho_T$. 
The uniqueness follows from Proposition \ref{CTY=OEY}. 
\end{proof}

\begin{corollary}
For a topological graph $E$, 
we have $\cT(E)\cong\cO(E_{\s{E}{0}{rg}})$. 
\end{corollary}

\begin{proof}
Take $Y=\s{E}{0}{rg}$ in Proposition \ref{cOEY}. 
\end{proof}

As a consequence of the analysis above, 
we get the following proposition 
which means that 
$\cO(E)$ is the smallest $C^*$-algebra 
among $C^*$-algebras generated by injective Toeplitz $E$-pairs 
which admit gauge actions. 
Thus we get an alternative definition of $\cO(E)$ 
which does not use the space $\s{E}{0}{rg}$ 
or Cuntz-Krieger pairs. 

\begin{proposition}\label{OEminimum}
Let $E$ be a topological graph, 
and $T$ be an injective Toeplitz $E$-pair 
which admits a gauge action. 
Then there exists a surjective $*$-ho\-mo\-mor\-phism 
$\varPsi\colon C^*(T)\to\cO(E)$ such that 
$\varPsi\circ T^i=t^i$ for $i=0,1$. 
\end{proposition}

\begin{proof}
We have the isomorphism $C^*(T)\cong\cO(E_{Y_T})$ 
by Corollary \ref{CT=OE}, 
and we have the surjection $\cO(E_{Y_T})\to\cO(E)$ 
by Proposition \ref{cOEY}. 
It is routine to check that the composition 
$\varPsi\colon C^*(T)\to\cO(E)$ of the above two maps satisfies 
$\varPsi\circ T^i=t^i$ for $i=0,1$. 
\end{proof}

\begin{remark}
The surjection $\cO(E_{Y_T})\to\cO(E)$ 
in the proof of Proposition \ref{OEminimum} 
is induced by the regular factor map $m=(m^0,m^1)$ from $E$ to $E_{Y_T}$ 
defined by 
the embeddings $m^0\colon E^0\to F_T^0$ and $m^1\colon E^1\to F_T^1$. 
\end{remark}

We finish this section by generalizing 
Corollary \ref{CT=OE} to 
all Toeplitz pairs admitting gauge actions. 
Let us fix a topological graph $E=(E^0,E^1,d,r)$ 
and a Toeplitz $E$-pair $T$. 
We define a closed subset $X^0$ of $E^0$ by 
$\ker T^0=C_0(E^0\setminus X^0)$. 
We set $X^1=d^{-1}(X^0)$ which is a closed subset of $E^1$. 

\begin{proposition}\label{XT}
We have $r(X^1)\subset X^0$. 
\end{proposition}

\begin{proof}
Take $e\in X^1$. 
We can find $\xi\in C_d(E^1)$ with $\xi(e)=1$ and $\xi(e')=0$ 
for all $e'\in d^{-1}(d(e))\setminus\{e\}$ 
because $e$ is isolated in $d^{-1}(d(e))$.
For any $f\in\ker T^0$, 
we have 
$$T^0(\ip{\xi}{\pi_r(f)\xi})=T^1(\xi)T^0(f)T^1(\xi)=0$$
Hence $\ip{\xi}{\pi_r(f)\xi}\in\ker T^0$.
We get $\ip{\xi}{\pi_r(f)\xi}(d(e))=0$ 
for all $f\in \ker T^0$. 
Since 
\begin{align*}
\ip{\xi}{\pi_r(f)\xi}(d(e))
=\sum_{e'\in d^{-1}(d(e))}\overline{\xi(e')}f(r(e'))\xi(e')
=f(r(e)),
\end{align*}
we have $f(r(e))=0$ for all $f\in\ker T^0$. 
This implies $r(e)\in X^0$. 
Thus $r(X^1)\subset X^0$. 
\end{proof}

By Proposition \ref{XT}, 
the quadruple $X=(X^0,X^1,d,r)$ is a topological graph. 
For $\xi\in C_d(E^1\setminus X^1)$, 
we can find $\eta\in C_d(E^1)$ and $f\in C_0(E^0\setminus X^0)$ 
with $\xi=\eta f$. 
Hence we have $T^1(\xi)=T^1(\eta)T^0(f)=0$. 
Thus the maps $T^0\colon C_0(E^0)\to C^*(T)$ 
and $T^1\colon C_d(E^1)\to C^*(T)$ factor through
$\dot{T}^0\colon C_0(X^0)\to C^*(T)$ 
and $\dot{T}^1\colon C_d(X^1)\to C^*(T)$. 
It is easy to see the following. 

\begin{lemma}
The pair $\dot{T}=(\dot{T}^0,\dot{T}^1)$ 
is an injective Toeplitz $X$-pair 
with $C^*(\dot{T})=C^*(T)$. 
The pair $\dot{T}$ admits a gauge action 
if and only if $T$ does. 
\end{lemma}

We define a closed subset $Y$ of $\s{X}{0}{rg}$ by 
$$C_0(\s{X}{0}{rg}\setminus Y)
=\{f\in C_0(\s{X}{0}{rg})\mid \dot{T}^0(f)=\dot{\varPhi}(\pi_r(f))\}.$$
where the $*$-ho\-mo\-mor\-phism 
$\dot{\varPhi}\colon \cK(C_d(X^1))\to C^*(T)$ is defined 
from the pair $\dot{T}$. 
The following proposition easily follows from Corollary \ref{CT=OE}. 

\begin{proposition}\label{CT=OF}
Let $E$ be a topological graph. 
For a Toeplitz $E$-pair $T$ 
which admits a gauge action, 
the $C^*$-algebra $C^*(T)$ 
is isomorphic to $\cO(X_Y)$ 
where the topological graph $X_Y$ 
is obtained from the topological graph $X$ 
by attaching extra vertices isomorphic to $Y$ 
and extra edges isomorphic to $d^{-1}(Y)$ as above. 
\end{proposition}

\begin{remark}
For a Toeplitz $E$-pair $T$ 
which does not admit a gauge action, 
we just get an injective Cuntz-Krieger $X_Y$-pair $T'$ on $C^*(T)$ 
such that $C^*(T')=C^*(T)$. 
\end{remark}

\begin{remark}
In general, 
the Toeplitz $X$-pair $\dot{T}$ 
defined above 
may not be a Cuntz-Krieger pair 
even when $T$ is a Cuntz-Krieger $E$-pair. 
This phenomena will be studied in the analysis of 
the ideal structures of $\cO(E)$ in \cite{Ka2}. 
\end{remark}

\section{Projective systems of topological graphs}\label{SecProj}

In this section, 
we define projective systems of topological graphs
and their projective limits, 
and study how these relate to $C^*$-algebras $\cT(E)$ and $\cO(E)$. 

\begin{definition}
A {\em projective system of topological graphs} 
over a directed set $\Lambda$ 
consists of a set of topological graphs 
$E_\lambda=(E_\lambda^0,E_\lambda^1,d_\lambda,r_\lambda)$ 
for $\lambda\in\Lambda$ 
and a set of factor maps $m_{\lambda,\lambda'}\colon E_{\lambda'}\to E_\lambda$ 
for $\lambda,\lambda'\in\Lambda$ with $\lambda\preceq\lambda'$ 
satisfying the equations 
$m_{\lambda,\lambda}=\id_{E_\lambda}$ for $\lambda\in\Lambda$
and 
$m_{\lambda,\lambda'}\circ m_{\lambda',\lambda''}=m_{\lambda,\lambda''}$
for $\lambda\preceq\lambda'\preceq\lambda''$. 
\end{definition}

Let us take a projective system of topological graphs 
$(\{E_\lambda\}_{\lambda\in\Lambda}, 
\{m_{\lambda,\lambda'}\}_{\lambda\preceq\lambda'})$ 
and fix it.
For $i=0,1$, 
we define a compact set $\widetilde{E}^i$ 
by the projective limit of the compact sets $\widetilde{E}_\lambda^i$ 
by the maps $m_{\lambda,\lambda'}^i$. 
For each $\lambda\in\Lambda$ and $i=0,1$, 
we denote by $m_\lambda^i$ 
the natural continuous map from $\widetilde{E}^i$ 
to $\widetilde{E}_\lambda^i$. 
For $i=0,1$,
we denote by $\infty$ the element $x\in\widetilde{E}^i$ 
with $m_\lambda^i(x)=\infty$ for all $\lambda\in\Lambda$, 
and set $E^i=\widetilde{E}^i\setminus\{\infty\}$. 
The set $E^i$ is a locally compact space 
whose one-point compactification is $\widetilde{E}^i$. 
For any $v\in E^0$, 
there exists $\lambda_0\in\Lambda$ 
such that $m_{\lambda_0}^0(v)\in E_{\lambda_0}^0$.
Then we have $m_\lambda^0(v)\in E_\lambda^0$ 
for any $\lambda\succeq\lambda_0$.
The net $\{m_\lambda^0(v)\}_{\lambda\succeq\lambda_0}$ satisfies the equation  
$m_{\lambda,\lambda'}^0(m_{\lambda'}^0(v))=m_\lambda^0(v)$ 
for $\lambda'\succeq\lambda\succeq\lambda_0$. 
Conversely, for $\lambda_0\in\Lambda$, 
a net $\{v_\lambda\}_{\lambda\succeq {\lambda_0}}$ 
with $v_\lambda\in E_\lambda^0$ 
satisfying the equation $m_{\lambda,\lambda'}^0(v_{\lambda'})=v_\lambda$ 
for $\lambda'\succeq\lambda\succeq\lambda_0$,
gives an element in $E^0$.
Thus elements in $E^0$ are represented by such nets.
Elements in $E^1$ are represented similarly.
We define maps $d,r\colon E^1\to E^0$ by 
$d(e)=\{d_\lambda(e_\lambda)\}_{\lambda\succeq\lambda_0}$ 
and $r(e)=\{r_\lambda(e_\lambda)\}_{\lambda\succeq\lambda_0}$ 
for $e=\{e_\lambda\}_{\lambda\succeq\lambda_0}\in E^1$.
This is well-defined because $m_{\lambda,\lambda'}$'s are factor maps. 
To prove that 
the quadruple $E=(E^0,E^1,d,r)$ is 
a topological graph, 
we need the following lemma. 

\begin{lemma}\label{homeo}
Let $m=(m^0,m^1)$ be a factor map 
from a topological graph $E=(E^0,E^1,d_E,r_E)$ 
to a topological graph $F=(F^0,F^1,d_F,r_F)$. 
If the restriction of $d_E$ to 
an open set $U\subset E^1$ is a homeomorphism 
onto the open set $d_E(U)\subset E^0$, 
then the restriction of $d_F$ to $(m^1)^{-1}(U)\subset F^1$ 
is a homeomorphism onto $(m^0)^{-1}(d_E(U))\subset F^0$. 
\end{lemma}

\begin{proof}
Take such an open subset $U\subset E^1$ 
and set $V=d_E(U)$. 
Since a bijective local homeomorphism is a homeomorphism, 
it suffices to see that 
the restriction of $d_F$ to $(m^1)^{-1}(U)\subset F^1$ 
is a bijection onto $(m^0)^{-1}(V)\subset F^0$. 
For $e\in (m^1)^{-1}(U)$, 
we have $m^0(d_F(e))=d_E(m^1(e))\in V$
by the condition (i) in Definition \ref{DefFactor}. 
Hence $d_F((m^1)^{-1}(U))\subset (m^0)^{-1}(V)$.
Take $v\in (m^0)^{-1}(V)$. 
Since $m^0(v)\in V$, 
there exists unique $e'\in U$ such that $d_E(e')=m^0(v)$.
By the condition (ii) in Definition \ref{DefFactor}, 
there exists unique $e\in F^1$ 
such that $m^1(e)=e'$ and $d_F(e)=v$.
This $e\in F^1$ is the unique element in $(m^1)^{-1}(U)$ with $d_F(e)=v$. 
Hence the restriction of $d_F$ to $(m^1)^{-1}(U)$ 
is a bijection onto $(m^0)^{-1}(V)$. 
We are done. 
\end{proof}

\begin{proposition}\label{projlim}
The quadruple $E=(E^0,E^1,d,r)$ defined 
from a projective system of topological graphs 
$(\{E_\lambda\}_{\lambda\in\Lambda}, 
\{m_{\lambda,\lambda'}\}_{\lambda\preceq\lambda'})$ 
as above is a topological graph. 
\end{proposition}

\begin{proof}
For $\lambda_0\in\Lambda$ and an open subset $V$ of $E^0_{\lambda_0}$, 
the set 
$$r^{-1}\big((m_{\lambda_0}^0)^{-1}(V)\big)
=\bigcup_{\lambda\succeq\lambda_0}
(m_{\lambda}^1)^{-1}\big(r_{\lambda}^{-1}(V)\big)$$ 
is an open subset of $E^1$. 
Since the family 
$$\{(m_{\lambda_0}^0)^{-1}(V)\mid \lambda_0\in\Lambda,\ 
V \mbox{ is an open subset of } E^0_{\lambda_0}\}$$
is a basis of $E^0$, 
we see that $r\colon E^1\to E^0$ is continuous. 
We will show that $d$ is locally homeomorphic. 
Take $e=\{e_\lambda\}_{\lambda\succeq\lambda_0}\in E^1$. 
Take a neighborhood $U_{\lambda_0}$ of $e_{\lambda_0}\in E_{\lambda_0}^1$ 
such that the restriction of $d_{\lambda_0}$ to $U_{\lambda_0}$ 
is a homeomorphism onto $d_{\lambda_0}(U_{\lambda_0})$. 
Set $V_{\lambda_0}=d_{\lambda_0}(U_{\lambda_0})$ 
which is a neighborhood 
of $d_{\lambda_0}(e_{\lambda_0})\in E_{\lambda_0}^0$. 
For $\lambda\succeq\lambda_0$, 
we set $U_\lambda=(m_{\lambda_0,\lambda}^1)^{-1}(U_{\lambda_0})$
and $V_\lambda=(m_{\lambda_0,\lambda}^0)^{-1}(V_{\lambda_0})$.
By Lemma \ref{homeo}, 
the restriction of $d_\lambda$ to $U_\lambda$ 
is a homeomorphism onto $V_\lambda$. 
Set $U=(m_{\lambda_0}^1)^{-1}(U_{\lambda_0})$ 
and $V=(m_{\lambda_0}^0)^{-1}(V_{\lambda_0})$. 
We see that $U$ is a neighborhood of $e\in E^1$
and $V$ is a neighborhood of $d(e)\in E^0$ 
because $m_{\lambda_0}^0(d(e))=d_{\lambda_0}(e_{\lambda_0})$. 
We will show that the restriction of $d$ to $U$ 
is a homeomorphism onto $V$. 
Take $\lambda\succeq\lambda_0$ 
and open subsets $V'\subset V_\lambda$ and $U'\subset U_\lambda$ 
with $V'=d_\lambda(U')$. 
For $e'\in (m_\lambda^1)^{-1}(U')\subset U$, 
we have 
$$m_{\lambda}^0(d(e'))=d_{\lambda}(m_{\lambda}^1(e'))\in V'.$$ 
Hence $d((m_\lambda^1)^{-1}(U'))\subset (m_\lambda^0)^{-1}(V')$. 
Take $v'\in (m_\lambda^0)^{-1}(V')\subset V$. 
Since the restriction of $d_{\lambda}$ to $U_{\lambda}$ 
is a bijection onto $V_{\lambda}$ and 
$m_{\lambda}^0(v')\in V'\subset V_{\lambda}$, 
there exists unique $e'_{\lambda}\in U_{\lambda}$
with $d_{\lambda}(e'_{\lambda})=m_{\lambda}^0(v')$. 
We have $e'_{\lambda}\in U'$. 
For each $\lambda'\succeq\lambda$ 
we have $m_{\lambda'}^0(v')\in V_{\lambda'}$. 
Hence there exists unique $e'_{\lambda'}\in U_{\lambda'}$
with $d_{\lambda'}(e'_{\lambda'})=m_{\lambda'}^0(v')$. 
By the uniqueness of $e'_{\lambda'}$, 
we have $m_{\lambda',\lambda''}^1(e'_{\lambda''})=e'_{\lambda'}$ 
for $\lambda''\succeq\lambda'\succeq\lambda$.
Thus $e'=\{e'_{\lambda'}\}_{\lambda'\succeq\lambda}
\in (m_\lambda^1)^{-1}(U')\subset U$ is 
the unique element in $U$ with $d(e')=v'$.
Hence the restriction of $d$ to $(m_\lambda^1)^{-1}(U')$ 
is a bijection onto $(m_\lambda^0)^{-1}(V')$. 
Since the family 
$$\{(m_\lambda^1)^{-1}(U')\mid \lambda\succeq\lambda_0,\ 
U' \mbox{ is an open subset of } U_\lambda\}$$
is a basis of $U$, 
and the family 
$$\{(m_\lambda^0)^{-1}(V')\mid \lambda\succeq\lambda_0,\ 
V' \mbox{ is an open subset of } V_\lambda\}$$
is a basis of $V$, 
we see that the restriction of $d$ to $U$ is a homeomorphism onto $V$. 
Thus $d\colon E^1\to E^0$ is a local homeomorphism. 
\end{proof}

\begin{definition}
The topological graph $E$ in Proposition \ref{projlim} 
is called the {\em projective limit} of the projective system 
$(\{E_\lambda\}_{\lambda\in\Lambda}, 
\{m_{\lambda,\lambda'}\}_{\lambda\preceq\lambda'})$, 
and denoted by 
$$\varprojlim \big(\{E_\lambda\}_{\lambda\in\Lambda}, 
\{m_{\lambda,\lambda'}\}_{\lambda\preceq\lambda'}\big),$$ 
or simply by $\varprojlim E_\lambda$. 
\end{definition}

\begin{proposition}\label{factormap}
For each $\lambda_0\in\Lambda$,
the pair $m_{\lambda_0}=(m_{\lambda_0}^0,m_{\lambda_0}^1)$ is 
a factor map from $E=\varprojlim E_\lambda$ to $E_{\lambda_0}$.
\end{proposition}

\begin{proof}
By the definitions of $d$ and $r$,
if $e\in E^1$ satisfies $m_{\lambda_0}^1(e)\in E_{\lambda_0}^1$,
then we have $d_{\lambda_0}(m_{\lambda_0}^1(e))=m_{\lambda_0}^0(d(e))$ and 
$r_{\lambda_0}(m_{\lambda_0}^1(e))=m_{\lambda_0}^0(r(e))$. 
Take $e_{\lambda_0}\in E_{\lambda_0}^1$ and $v\in E^0$ 
such that $d_{\lambda_0}(e_{\lambda_0})=m_{\lambda_0}^0(v)$.
For each $\lambda\succeq\lambda_0$,
we have 
$d_{\lambda_0}(e_{\lambda_0})=m_{\lambda_0,\lambda}^0(m_{\lambda}^0(v))$.
Since $m_{\lambda_0,\lambda}$ is a factor map, 
there exists unique $e_{\lambda}\in E_{\lambda}^1$ 
with $m_{\lambda_0,\lambda}^1(e_{\lambda})=e_{\lambda_0}$ 
and $d_{\lambda}(e_{\lambda})=m_{\lambda}^0(v)$.
The uniqueness implies that 
$m_{\lambda,\lambda'}^1(e_{\lambda'})=e_{\lambda}$ 
for $\lambda'\succeq\lambda\succeq{\lambda_0}$.
Hence $e=\{e_{\lambda}\}_{\lambda\succeq\lambda_0}$ 
is the unique element in $E^1$ satisfying 
the equations $m_{{\lambda_0}}^1(e)=e_{{\lambda_0}}$ and $d(e)=v$.
We are done. 
\end{proof}

Denote by $T_\lambda=(T_\lambda^0,T_\lambda^1)$ 
the universal Toeplitz $E_\lambda$-pair in $\cT(E_\lambda)$. 
For $\lambda\preceq\lambda'$, 
the factor map $m_{\lambda,\lambda'}$ gives 
a $*$-ho\-mo\-mor\-phism 
$\mu_{\lambda,\lambda'}^0\colon C_0(E_\lambda^0)\to C_0(E_{\lambda'}^0)$, 
a linear map 
$\mu_{\lambda,\lambda'}^1\colon C_{d_\lambda}(E_\lambda^1)\to 
C_{d_{\lambda'}}(E_{\lambda'}^1)$ 
and a $*$-ho\-mo\-mor\-phism 
$\mu_{\lambda,\lambda'}\colon \cT(E_{\lambda})\to\cT(E_{\lambda'})$
such that 
$\mu_{\lambda,\lambda'}\circ T_{\lambda}^i
=T_{\lambda'}^i\circ\mu_{\lambda,\lambda'}^i$ 
for $i=0,1$ 
by Proposition \ref{functT}. 
For $\lambda\preceq\lambda'\preceq\lambda''$, 
we have 
$\mu_{\lambda',\lambda''}^i\circ\mu_{\lambda,\lambda'}^i=
\mu_{\lambda,\lambda''}^i$ for $i=0,1$ and 
$\mu_{\lambda',\lambda''}\circ\mu_{\lambda,\lambda'}=
\mu_{\lambda,\lambda''}$ by Proposition \ref{comp}. 
For each $\lambda\in\Lambda$, 
the factor map $m_\lambda$ gives 
a $*$-ho\-mo\-mor\-phism 
$\mu_{\lambda}^0\colon C_0(E_{\lambda}^0)\to C_0(E^0)$, 
a linear map $\mu_{\lambda}^1\colon C_{d_\lambda}(E_{\lambda}^1)\to C_d(E^1)$ 
and a $*$-ho\-mo\-mor\-phism 
$\mu_{\lambda}\colon \cT(E_{\lambda})\to\cT(E)$ 
by Proposition \ref{factormap}. 
Since we have $\mu_{\lambda'}\circ\mu_{\lambda,\lambda'}=\mu_{\lambda}$ 
for $\lambda\preceq\lambda'$, 
$\{\mu_{\lambda}\}_{\lambda\in\Lambda}$ 
gives us a $*$-ho\-mo\-mor\-phism 
$\varinjlim \cT(E_{\lambda})\to\cT(E)$. 
We can naturally consider 
$C_0(E^0)=\varinjlim C_0(E_{\lambda}^0)$ and 
$C_d(E^1)=\varinjlim C_{d_\lambda}(E_{\lambda}^1)$. 
There exist a $*$-ho\-mo\-mor\-phism 
$T^0\colon C_0(E^0)\to \varinjlim \cT(E_{\lambda})$ 
and a linear map $T^1\colon C_d(E^1)\to \varinjlim \cT(E_{\lambda})$ 
such that $T^i\circ\mu_{\lambda}^i=\mu_{\lambda}\circ T_\lambda^i$ 
for $i=0,1$ and $\lambda\in\Lambda$. 
One can easily see that the pair $T=(T^0,T^1)$ is a Toeplitz $E$-pair. 
By the universality of $\cT(E)$, 
there exists a $*$-ho\-mo\-mor\-phism 
$\cT(E)\to\varinjlim \cT(E_{\lambda})$. 
It is easy to verify that the two maps 
$\varinjlim \cT(E_{\lambda})\to\cT(E)$ and 
$\cT(E)\to\varinjlim \cT(E_{\lambda})$ are 
the inverses of each others. 
Thus we get the following.

\begin{proposition}\label{Tind}
For a projective system 
$(\{E_\lambda\}_{\lambda\in\Lambda}, 
\{m_{\lambda,\lambda'}\}_{\lambda\preceq\lambda'})$, 
we have 
$$\cT\big(\varprojlim E_{\lambda}\big)\cong\varinjlim \cT(E_{\lambda}).$$ 
\end{proposition}

We define the regurality of projective limits 
and seek an analogous result of Proposition \ref{Tind} 
for $\cO(E)$. 

\begin{definition}
A projective system 
$(\{E_\lambda\}_{\lambda\in\Lambda},
\{m_{\lambda,\lambda'}\}_{\lambda\preceq\lambda'})$ 
of topological graphs 
is said to be {\em regular} if $m_{\lambda,\lambda'}$ is regular 
for all $\lambda,\lambda'\in\Lambda$ with $\lambda\preceq\lambda'$. 
\end{definition}

Take a regular projective system 
$(\{E_\lambda\}_{\lambda\in\Lambda}, 
\{m_{\lambda,\lambda'}\}_{\lambda\preceq\lambda'})$ 
of topological graphs, 
and denote by $E=(E^0,E^1,d,r)$ 
its projective limit $\varprojlim E_\lambda$. 
Let $m_{\lambda_0}=(m_{\lambda_0}^0,m_{\lambda_0}^1)$ be 
the natural factor map from $E=\varprojlim E_\lambda$ 
to $E_{\lambda_0}$ for each $\lambda_0\in\Lambda$. 

\begin{proposition}\label{factormapreg}
For each $\lambda_0\in\Lambda$,
the factor map $m_{\lambda_0}$ is regular.
\end{proposition}

\begin{proof}
Take $v\in E^0$ with $m_{\lambda_0}^0(v)\in (E_{\lambda_0}^0)_{\rs{rg}}$,
and we will show that $r^{-1}(v)$ is a non-empty subset 
of $(m_{\lambda_0}^1)^{-1}(E_{\lambda_0}^1)$.
Let $e$ be an element in $r^{-1}(v)$.
Take $\lambda\in\Lambda$ 
with $m_{\lambda}^1(e)\in E_{\lambda}^1$.
We may assume $\lambda\succeq\lambda_0$.
Since the factor map $m_{\lambda_0,\lambda}$ is regular and 
$$m_{\lambda_0,\lambda}^0(r_{\lambda}(m_{\lambda}^1(e)))
=m_{\lambda_0,\lambda}^0(m_{\lambda}^0(r(e)))
=m_{\lambda_0}^0(v)\in (E_{\lambda_0}^0)_{\rs{rg}},$$
we have $m_{\lambda}^1(e)\in (m_{\lambda_0,\lambda}^1)^{-1}(E_{\lambda_0}^1)$.
Thus we have $e\in (m_{\lambda_0}^1)^{-1}(E_{\lambda_0}^1)$. 
Hence $r^{-1}(v)\subset (m_{\lambda_0}^1)^{-1}(E_{\lambda_0}^1)$. 
We will show that $r^{-1}(v)\neq\emptyset$. 

Set $\Lambda_0=\{\lambda\in\Lambda\mid \lambda\succeq\lambda_0\}$. 
For each $\lambda\in\Lambda_0$, 
we define $G_\lambda\subset E_{\lambda}^1$ 
by $G_\lambda=r_{\lambda}^{-1}(m_{\lambda}^0(v))$. 
Lemma \ref{regular} implies 
$m_{\lambda}^0(v)\in (E_{\lambda}^0)_{\rs{rg}}$ 
for $\lambda\in\Lambda_0$ 
because $m_{\lambda_0,\lambda}^0(m_{\lambda}^0(v))=m_{\lambda_0}^0(v)
\in (E_{\lambda_0}^0)_{\rs{rg}}$. 
Hence $G_\lambda$ is a non-empty compact set. 
We define a compact set $G$ 
by $G=\prod_{\lambda\in\Lambda_0}G_\lambda$. 
For $\lambda,\lambda'\in \Lambda_0$ with $\lambda\preceq\lambda'$, 
we define a subset 
$$F_{\lambda,\lambda'}
=\big\{\{e_\lambda\}_{\lambda\in \Lambda_0}\in G\ \big|\ 
m_{\lambda,\lambda'}(e_{\lambda'})=e_\lambda\big\}$$
It is clear that $F_{\lambda,\lambda'}$ 
is a closed subset of the compact set $G$ 
for $\lambda,\lambda'\in \Lambda_0$ with $\lambda\preceq\lambda'$. 
We will show that 
$\bigcap_{k=1}^n F_{\lambda_k,\lambda_k'}\neq \emptyset$ 
for $\lambda_k,\lambda_k'\in \Lambda_0$ 
with $\lambda_k\preceq\lambda_k'$ 
for $k=1,2,\ldots,n$. 
We can find $\lambda_0'\in \Lambda_0$ 
with $\lambda_k'\preceq\lambda_0'$ for $k=1,2,\ldots,n$. 
Since $G_{\lambda_0'}\neq \emptyset$, 
we can take $e_{\lambda_0'}\in G_{\lambda_0'}$. 
By the regurality of $m_{\lambda,\lambda_0'}$, 
we have $m_{\lambda,\lambda_0'}(e_{\lambda_0'})\in G_{\lambda}$ 
for $\lambda\in \Lambda_0$ with $\lambda\preceq\lambda_0'$. 
Thus we can find $\{e_{\lambda}\}_{\lambda\in\Lambda_0}\in G$ 
such that $e_{\lambda}=m_{\lambda,\lambda_0'}(e_{\lambda_0'})$ 
with $\lambda\in \Lambda_0$ with $\lambda\preceq\lambda_0'$. 
For $k=1,2,\ldots,n$, 
we have 
$$m_{\lambda_k,\lambda_k'}(e_{\lambda_k'})
=m_{\lambda_k,\lambda_k'}(m_{\lambda_k',\lambda_0'}(e_{\lambda_0'}))
=m_{\lambda_k,\lambda_0'}(e_{\lambda_0'})
=e_{\lambda_k}.$$
Hence we have shown that 
$$\{e_{\lambda}\}_{\lambda\in\Lambda_0}
\in \bigcap_{k=1}^n F_{\lambda_k,\lambda_k'}\neq \emptyset$$ 
for $\lambda_k,\lambda_k'\in \Lambda_0$ 
with $\lambda_k\preceq\lambda_k'$ 
for $k=1,2,\ldots,n$. 
Since $G$ is compact, 
we can find 
$$\{e_{\lambda}\}_{\lambda\in\Lambda_0}\in 
\bigcap_{\lambda\preceq\lambda'} F_{\lambda,\lambda'}.$$
Since we have $m_{\lambda,\lambda'}(e_{\lambda'})=e_\lambda$ 
for $\lambda,\lambda'\in \Lambda_0$ with $\lambda\preceq\lambda'$. 
we get an element 
$e=\{e_{\lambda}\}_{\lambda\in\Lambda_0}\in E^1$ 
which satisfies $r(e)=v$. 
Hence $r^{-1}(v)\neq\emptyset$. 
The proof is completed. 
\end{proof}

As in the case of $\cT(E)$, 
we get the following commutative diagram ($\lambda\preceq\lambda'$): 
\begin{center}
\mbox{\xymatrix{
C_0(E_\lambda^0) \ar_{\mu_{\lambda,\lambda'}^0}[r] 
                 \ar@/^0.8pc/^{\mu_{\lambda}^0}[rrrr] 
                 \ar_{t_{\lambda}^0}[d] &
C_0(E_{\lambda'}^0) \ar_{\mu_{\lambda'}^0}[rrr]  
                    \ar_{t_{\lambda'}^0}[d] &&&
C_0(E^0) \ar_{t^0}[d] \\
\cO(E_\lambda) \ar_{\mu_{\lambda,\lambda'}}[r] 
               \ar@/^0.8pc/^{\mu_{\lambda}}[rrrr] &
\cO(E_{\lambda'}) \ar_{\mu_{\lambda'}}[rrr] &&& 
\cO(E) \\
C_{d_\lambda}(E_\lambda^1) \ar_{\mu_{\lambda,\lambda'}^1}[r] 
                           \ar@/^0.8pc/^{\mu_{\lambda}^1}[rrrr] 
                           \ar^{t_{\lambda}^1}[u] &
C_{d_{\lambda'}}(E_{\lambda'}^1) \ar_{\mu_{\lambda'}^1}[rrr] 
                                 \ar^{t_{\lambda'}^1}[u] &&&
C_d(E^1) \ar^{t^1}[u]
}}
\end{center}
We denote the natural map 
by $\nu_\lambda\colon \cO(E_\lambda)\to\varinjlim\cO(E_\lambda)$ 
for each $\lambda\in\Lambda$. 
For each $\lambda\in\Lambda$, 
we define maps $T^0\colon C_0(E^0)\to \varinjlim\cO(E_\lambda)$
and $T^1\colon C_d(E^1)\to \varinjlim\cO(E_\lambda)$ so that 
$T^i\circ\mu_\lambda^i=\nu_\lambda\circ t_\lambda^i$ for $i=0,1$. 
By the universality of inductive limits, 
there exists a $*$-ho\-mo\-mor\-phism 
$\nu\colon \varinjlim\cO(E_\lambda)\to \cO(E)$
with $\nu\circ\nu_\lambda=\mu_\lambda$ for $\lambda\in \Lambda$. 
We have $\nu\circ T^i=t^i$ for $i=0,1$. 
Hence the image of $\nu$ contains 
$t^0(C_0(E^0))$ and $t^1(C_d(E^1))$, 
and so $\nu$ is surjective. 
One might expect that $\nu$ is an isomorphism, 
as it is in the case of $\cT(E)$. 
However this is not the case as the following example shows. 

\begin{example}\label{ex1}
Let $F=(F^0,F^1,d_F,r_F)$ be a discrete graph 
given by 
$$
\begin{array}{cc}
F^0=\{v,v',w\},&
F^1=\{e_k\}_{k\in\N},\\
d_F(e_k)=\left\{\begin{array}{ll}
v& (k=0)\\
v'& (k\geq 1)
\end{array}\right. ,&
r_F(e_k)=w\ (k\in\N).
\end{array}
\qquad
\raisebox{1cm}{\xymatrix{
&\bullet\ v' \ar@<-1.8ex>[d]_{e_1} \ar@<-1.0ex>[d]^{\cdots}\\
v\ \bullet \ar[r]_{e_0} & \bullet\ w}}
$$
We have $\s{F}{0}{sce}=\{v,v'\}$ and $\s{F}{0}{inf}=\{w\}$.
Hence $\s{F}{0}{rg}=\emptyset$. 
We see that $\cO(F)\cong\M_2\oplus\widetilde{\K}$ 
where $\widetilde{\K}$ means the unitization of $\K$. 
Define a regular factor map $m=(m^0,m^1)$ from $F$ to itself by 
$$m^0(v)=v, m^0(w)=w, m^0(v')=\infty,$$
$$m^1(e_0)=e_0, \ m^1(e_k)=\infty\  (k\geq 1).$$
Under the isomorphism $\cO(F)\cong\M_2\oplus\widetilde{\K}$, 
the $*$-ho\-mo\-mor\-phism $\mu\colon \cO(F)\to\cO(F)$ induced by $m$ 
is expressed as 
$$\M_2\oplus\widetilde{\K}\ni (x,y+\lambda)\mapsto 
(x,\lambda)\in \M_2\oplus\widetilde{\K}$$
for $x\in\M_2$, $y\in\K$ and $\lambda\in\C$. 
For $k\in\N$, define $E_k=F$ and $m_{k,k+1}=m$. 
Then $\{E_k\}_{k\in\N}$ and $\{m_{k,k+1}\}_{k\in\N}$ 
give a regular projective system. 
Its projective limit $E=(E^0,E^1,d,r)$ is a discrete graph 
such that $E^0=\{v,w\}, E^1=\{e_0\}, d(e_0)=v, r(e_0)=w$. 
We have $\cO(E)\cong\M_2$. 
On the other hand, we have $\varinjlim \cO(E_k)\cong \M_2\oplus\C$. 
Hence the surjection $\varinjlim \cO(E_k)\to \cO(E)$ 
is not an isomorphism. 
Note that $\varinjlim \cO(E_k)\cong\cT(E)$. 
\end{example}

We define an open set $O\subset E^0$ by 
$$O=\bigcup_{\lambda\in\Lambda}
(m_{\lambda}^0)^{-1}((E_{\lambda}^0)_{\rs{rg}}).$$
From Lemma \ref{regular} and Proposition \ref{factormapreg}, 
we see that 
$(m_{\lambda}^0)^{-1}((E_{\lambda}^0)_{\rs{rg}})\subset\s{E}{0}{rg}$ 
for every $\lambda\in\Lambda$. 
Hence $O$ is an open subset of $\s{E}{0}{rg}$. 

\begin{lemma}\label{C0O}
We have $\{f\in C_0(\s{E}{0}{rg})\mid T^0(f)=\varPhi(\pi_r(f))\}=C_0(O)$. 
\end{lemma}

\begin{proof}
Take $\lambda\in\Lambda$ and $g\in C_0((E_{\lambda}^0)_{\rs{rg}})$. 
We have $\mu_{\lambda}^0(g)\in C_0(O)\subset C_0(\s{E}{0}{rg})$. 
We define 
$\psi_{\lambda}\colon \cK(C_{d_{\lambda}}(E_{\lambda}^1))\to\cK(C_d(E^1))$
by $\psi_{\lambda}(\theta_{\xi,\eta})
=\theta_{\mu_\lambda^1(\xi),\mu_\lambda^1(\eta)}$
for $\xi,\eta\in C_{d_{\lambda}}(E_{\lambda}^1)$.
Then we have 
$\varPhi\circ\psi_{\lambda}=\nu_{\lambda}\circ\varphi_{\lambda}$.
By Lemma \ref{funct0}, we have
\begin{align*}
T^0(\mu_{\lambda}^0(g))
&=\nu_{\lambda}(t_{\lambda}^0(g)) 
 =\nu_{\lambda}(\varphi_{\lambda}(\pi_{r_{\lambda}}(g)))\\
&=\varPhi(\psi_{\lambda}(\pi_{r_{\lambda}}(g))) 
 =\varPhi(\pi_r(\mu_{\lambda}^0(g))). 
\end{align*}
Since 
$$C_0(O)=\overline{
\bigcup_{\lambda\in\Lambda}\mu_{\lambda}^0
\big(C_0((E_{\lambda}^0)_{\rs{rg}})\big)},$$ 
we have $T^0(f)=\varPhi(\pi_r(f))$ for 
all $f\in C_0(O)$.

To derive a contradiction, 
assume that there exists $f\in C_0(\s{E}{0}{rg})$
such that $T^0(f)=\varPhi(\pi_r(f))$ and 
$f\notin C_0(O)$.
There exists $v\notin O$ 
with $|f(v)|=\e>0$. 
We can find $\lambda_0\in\Lambda$ and $g_0\in C_0(E_\lambda^0)$
such that $\|\nu_{\lambda_0}(t_{\lambda_0}^0(g_0))-T^0(f)\|<\e/3$.
Since 
$$\varPhi\big(\cK(C_d(E^1))\big)
=\overline{\bigcup_{\lambda\in\Lambda}\nu_{\lambda}
 \big(\varphi_{\lambda}\big(\cK(C_{d_{\lambda}}(E_{\lambda}^0))\big)\big)},$$
we can find $\lambda_1\in\Lambda$ 
and $x_1\in\cK(C_{d_{\lambda_1}}(E_{\lambda_1}^0))$
such that 
$\|\nu_{\lambda_1}(\varphi_{\lambda_1}(x_1))-\varPhi(\pi_r(f))\|<\e/3$.
Take $\lambda_2\in\Lambda$ so that $\lambda_2\succeq\lambda_0$ and 
$\lambda_2\succeq\lambda_1$, and set 
$$a=\mu_{\lambda_0,\lambda_2}\big(t_{\lambda_0}^0(g_0)\big)-
\mu_{\lambda_1,\lambda_2}\big(\varphi_{\lambda_1}(x_1)\big)
\in\cO(E_{\lambda_2}).$$
Then we have 
\begin{align*}
\|\nu_{\lambda_2}(a)\|
&=\big\|\nu_{\lambda_0}\big(t_{\lambda_0}^0(g_0)\big)-
    \nu_{\lambda_1}\big(\varphi_{\lambda_1}(x_1)\big)\big\|\\
&\leq\big\|\nu_{\lambda_0}\big(t_{\lambda_0}^0(g_0)\big)-T^0(f)\big\|
 +\big\|\varPhi\big(\pi_r(f)\big)
        -\nu_{\lambda_1}\big(\varphi_{\lambda_1}(x_1)\big)\big\|\\
&<2\e/3.
\end{align*}
By the definition of the inductive limit of $C^*$-algebras, 
there exists $\lambda_3\succeq\lambda_2$ 
such that $\|\mu_{\lambda_2,\lambda_3}(a)\|<2\e/3$. 
Take $\lambda\in\Lambda$ such that $\lambda\succeq\lambda_3$ and 
$m_\lambda^0(v)\in E_{\lambda}^0$. 
Set $g=\mu_{\lambda_0,\lambda}^0(g_0)\in C_0(E_{\lambda}^0)$ 
and $x=\psi_{\lambda_1,\lambda}(x_1)\in\cK(C_{d_{\lambda}}(E_{\lambda}^1))$, 
where $\psi_{\lambda_1,\lambda}\colon \cK(C_{d_{\lambda_1}}(E_{\lambda_1}^1))\to
\cK(C_{d_{\lambda}}(E_{\lambda}^1))$ 
is defined similarly as $\psi_{\lambda}$. 
Then we have 
\begin{align*}
\|t_{\lambda}^0(g)-\varphi(x)\|
&=\big\|\mu_{\lambda_0,\lambda}\big(t_{\lambda_0}^0(g_0)\big)
  -\mu_{\lambda_1,\lambda}\big(\varphi_{\lambda_1}(x_1)\big)\big\|\\
&=\big\|\mu_{\lambda_2,\lambda}\big(
   \mu_{\lambda_0,\lambda_2}(t_{\lambda_0}^0(g_0))
  -\mu_{\lambda_1,\lambda_2}(\varphi_{\lambda_1}(x_1))\big)\big\|\\
&=\|\mu_{\lambda_2,\lambda}(a)\|
 =\big\|\mu_{\lambda_3,\lambda}\big(\mu_{\lambda_2,\lambda_3}(a)\big)\big\|
 \leq\|\mu_{\lambda_2,\lambda_3}(a)\|
 <2\e/3.
\end{align*}
Since $v\notin O$, 
we have $m_\lambda^0(v)\in (E_{\lambda}^0)_{\rs{sg}}$. 
Since 
$$\|\mu_{\lambda}^0(g)-f\|=\|T^0(\mu_{\lambda_0}^0(g_0)-f)\|
 =\|\nu_{\lambda_0}(t_{\lambda_0}^0(g_0))-T^0(f)\|<\e/3,$$
we have $|(\mu_{\lambda}^0(g)-f)(v)|<\e/3$. 
Hence we have $|g(m_\lambda^0(v))|>2\e/3$ because $|f(v)|=\e$.
Since $\F^1/\G^1\cong C_0((E_{\lambda}^0)_{\rs{sg}})$,
we have 
$$\inf_{x'\in\cK(C_{d_{\lambda}}(E_{\lambda}^1))}
\|t_{\lambda}^0(g)-\varphi_{\lambda}(x')\|
=\sup_{v'\in (E_{\lambda}^0)_{\rs{sg}}}|g(v')|.$$
However, we have $\|t_{\lambda}^0(g)-\varphi_{\lambda}(x)\|<2\e/3$ 
and $|g(m_\lambda^0(v))|>2\e/3$.
This is a contradiction.
Thus we have shown that $T^0(f)=\varPhi(\pi_r(f))$ if and only if 
$f\in C_0(O)$.
\end{proof}

We define $Y=\s{E}{0}{rg}\setminus O$ 
which is a closed subset of $\s{E}{0}{rg}$. 
In Example \ref{ex1}, we have $O=\emptyset$ and $\s{E}{0}{rg}=\{w\}$, 
hence $Y=\{w\}$. 
It is easy to see the following. 

\begin{proposition}\label{O=Erg}
The inductive limit $\varinjlim\cO(E_\lambda)$ 
is isomorphic to $\cO(E_Y)$, 
and the surjection $\nu\colon \varinjlim\cO(E_\lambda)\to \cO(E)$ 
is an isomorphism if and only if $Y=\emptyset$. 
\end{proposition}

Let us say that a projective system 
$(\{E_\lambda\}_{\lambda\in\Lambda}, 
\{m_{\lambda,\lambda'}\}_{\lambda\preceq\lambda'})$ 
is surjective 
when $m_{\lambda,\lambda'}^0$ is surjective 
for every $\lambda\preceq\lambda'$. 
By Proposition \ref{funct1}, 
we get an injective inductive system of $C^*$-algebras 
from a surjective regular projective system. 

\begin{lemma}\label{limE0r}
If a regular projective system 
$(\{E_\lambda\}_{\lambda\in\Lambda}, 
\{m_{\lambda,\lambda'}\}_{\lambda\preceq\lambda'})$ 
is surjective, 
then we have $Y=\emptyset$. 
\end{lemma}

\begin{proof}
Suppose that $m_{\lambda,\lambda'}^0$ is surjective 
for every $\lambda\preceq\lambda'$.
Note that $m_{\lambda,\lambda'}^1$ is surjective 
for every $\lambda\preceq\lambda'$
and that $m_{\lambda}^0,m_{\lambda}^1$ is surjective 
for every $\lambda\in\Lambda$.
To prove that $Y=\emptyset$, 
it suffices to see that 
for $v\in\s{E}{0}{rg}$, 
there exists $\lambda\in\Lambda$ 
such that $m_{\lambda}^0(v)\in (E_\lambda^0)_{\rs{rg}}$.
Take $v\in\s{E}{0}{rg}$. 
By Lemma \ref{E0r}, 
there exists a compact neighborhood $V$ of $v$ 
such that $r^{-1}(V)$ is compact
and $r(r^{-1}(V))=V$. 
Since $V$ is a neighborhood of $v$, 
there exist $\lambda_0\in\Lambda$ 
and a neighborhood $V'$ of $m_{\lambda_0}^0(v)$ 
such that $(m_{\lambda_0}^0)^{-1}(V')\subset V$.
For $\lambda\succeq\lambda_0$, 
we have $(m_{\lambda_0,\lambda}^0)^{-1}(V')\subset m_{\lambda}^0(V)$
because $m_{\lambda}^0$ is surjective. 
Hence $m_{\lambda}^0(V)$ is a neighborhood of $m_{\lambda}^0(v)$ 
for $\lambda\succeq\lambda_0$. 
Since $\widetilde{E}^0\setminus V$ is a neighborhood of 
$\infty\in\widetilde{E}^0$, 
there exists $\lambda_1\in\Lambda$ such that 
$(m_{\lambda_1}^0)^{-1}(V'')\subset \widetilde{E}^0\setminus V$ 
for some neighborhood $V''$ of $\infty\in\widetilde{E}_{\lambda_1}^0$. 
Hence $\infty\notin m_{\lambda_1}(V)$.
This implies that $m_{\lambda}^0(V)\subset E_{\lambda}^0$ 
for every $\lambda\succeq\lambda_1$. 
Similarly there exists $\lambda_2\in\Lambda$ 
such that $m_{\lambda}^1(r^{-1}(V))\subset E_{\lambda}^1$ 
for every $\lambda\succeq\lambda_2$ 
because $r^{-1}(V)\subset E^1$ is compact. 
Take $\lambda\in\Lambda$ with $\lambda\succeq\lambda_i$ for $i=0,1,2$. 
Set $V_\lambda=m_{\lambda}^0(V)$ and $U_\lambda=m_{\lambda}^1(r^{-1}(V))$. 
Then $V_\lambda$ is a compact neighborhood of $m_{\lambda}^0(v)$ 
in $E_\lambda^0$ and $U_\lambda$ is a compact subset of $E_\lambda^1$. 
We get $U_\lambda=r_\lambda^{-1}(V_\lambda)$ 
and $r_\lambda(U_\lambda)=V_\lambda$ 
because $m_{\lambda}^0$ and $m_{\lambda}^1$ are surjective
and $r(r^{-1}(V))=V$. 
Hence we have $m_{\lambda}^0(v)\in (E_\lambda^0)_{\rs{rg}}$ 
by Lemma \ref{E0r}. 
Thus we have shown that $Y=\emptyset$. 
\end{proof}

By Proposition \ref{O=Erg} and Lemma \ref{limE0r}, 
we get the following which is satisfactory for application. 

\begin{proposition}\label{projisom}
For a surjective regular projective system 
$(\{E_\lambda\}_{\lambda\in\Lambda},
\{m_{\lambda,\lambda'}\}_{\lambda\preceq\lambda'})$, 
we have $\cO\big(\varprojlim E_\lambda\big)\cong\varinjlim\cO(E_\lambda)$. 
\end{proposition}

\section{Subalgebras of $\cO(E)$}\label{Secsubalg}

Let us take a topological graph $E=(E^0,E^1,d,r)$ and fix it.
In this section, we study subalgebras of $\cO(E)$. 

\begin{definition}
A {\em subgraph} of the topological graph $E=(E^0,E^1,d,r)$ 
is a quadruple $F=(F^0,F^1,d_F,r_F)$ 
where $F^0\subset E^0$ and $F^1\subset E^1$ 
are open subsets such that $d(F^1),r(F^1)\subset F^0$, 
and $d_F,r_F$ are the restrictions of $d,r$ to $F^1$. 
\end{definition}

Take a subgraph $F=(F^0,F^1,d_F,r_F)$ 
of the topological graph $E=(E^0,E^1,d,r)$. 
We can identify $C_{d_F}(F^1)$ 
with $C_d(E^1)\cap C_0(F^1)$. 
The $C^*$-algebra generated by $t^0(C_0(F^0))$ 
and $t^1(C_{d_F}(F^1))$ is different from $\cO(F)$ in general. 
We will explore what the difference is. 
Let us define a $*$-ho\-mo\-mor\-phism $T^0\colon C_0(F^0)\to\cO(E)$ 
and a linear map $T^1\colon C_{d_F}(F^1)\to\cO(E)$ 
by the restrictions of $t^0$ and $t^1$, respectively. 
It is clear that the pair $T=(T^0,T^1)$ is 
an injective Toeplitz $F$-pair. 
In general, it is not a Cuntz-Krieger $F$-pair. 

\begin{lemma}\label{lemFY}
For $f\in C_0(F^0)$, 
we have $T^0(f)=\varPhi(\pi_{r_F}(f))$ 
if and only if 
$f\in C_0(\s{F}{0}{rg}\setminus \overline{r(E^1\setminus F^1)})$. 
\end{lemma}

\begin{proof}
First note that the $*$-ho\-mo\-mor\-phism 
$\varPhi\colon \cK(C_{d_F}(F^1))\to \cO(E)$ 
obtained by the Toeplitz $F$-pair $T$ 
is the restriction of the map 
$\varphi\colon \cK(C_d(E^1))\to \cO(E)$ 
to $\cK(C_{d_F}(F^1))\subset \cK(C_d(E^1))$. 
We first show that 
$f\in C_0(F^0)$ satisfies $T^0(f)=\varPhi(\pi_{r_F}(f))$ 
if and only if 
$f\in C_0(F^0\cap \s{E}{0}{rg}\setminus \overline{r(E^1\setminus F^1)})$. 
Let us take $f\in C_0(F^0)$ satisfies $T^0(f)=\varPhi(\pi_{r_F}(f))$. 
Since we have $t^0(f)\in \varphi(\cK(C_d(E^1)))$, 
we get $f\in C_0(\s{E}{0}{rg})$ and 
$$\varPhi(\pi_{r_F}(f))=T^0(f)=t^0(f)=\varphi(\pi_{r}(f)).$$
The latter condition implies that $f\circ r\in C_0(F^1)$ 
and so $f\in C_0(E^0\setminus \overline{r(E^1\setminus F^1)})$. 
Therefore 
$f\in C_0(F^0\cap \s{E}{0}{rg}\setminus \overline{r(E^1\setminus F^1)})$. 
Conversely if $f\in C_0(F^0)$ satisfies 
$f\in C_0(F^0\cap \s{E}{0}{rg}\setminus \overline{r(E^1\setminus F^1)})$, 
then we get 
$$T^0(f)=t^0(f)=\varphi(\pi_r(f))=\varPhi(\pi_{r_F}(f)).$$ 
Thus we have shown that 
$f\in C_0(F^0)$ satisfies $T^0(f)=\varPhi(\pi_{r_F}(f))$ 
if and only if 
$f\in C_0\big(F^0\cap \s{E}{0}{rg}\setminus 
\overline{r(E^1\setminus F^1)}\big)$. 
The proof completes once we show 
$$F^0\cap\s{E}{0}{rg}\setminus \overline{r(E^1\setminus F^1)}
=\s{F}{0}{rg}\setminus \overline{r(E^1\setminus F^1)}.$$
Take $v\in \s{F}{0}{rg}\setminus \overline{r(E^1\setminus F^1)}$. 
By Proposition \ref{E0r}
we can find a neighborhood $V\subset F^0$ of $v\in F^0$ 
satisfying the conditions 
that $(r_F)^{-1}(V)\subset F^1$ is compact 
and $r_F((r_F)^{-1}(V))=V$. 
By replacing it by a smaller set if necessary, 
we may assume that $V\cap r(E^1\setminus F^1)=\emptyset$. 
Since $F^0$ is an open subset of $E^0$, 
$V$ is a neighborhood of $v\in E^0$. 
Since $V\cap r(E^1\setminus F^1)=\emptyset$, 
we have $r^{-1}(V)\subset F^1$. 
Hence $r^{-1}(V)=(r_F)^{-1}(V)$ is compact 
and satisfies $r(r^{-1}(V))=V$. 
By Proposition \ref{E0r}, we have $v\in \s{E}{0}{rg}$. 
Thus we get 
$v\in F^0\cap\s{E}{0}{rg}\setminus \overline{r(E^1\setminus F^1)}$. 
Conversely, take 
$v\in F^0\cap\s{E}{0}{rg}\setminus \overline{r(E^1\setminus F^1)}$. 
Then we can find a neighborhood $V$ of $v\in E^0$ 
such that $r^{-1}(V)$ is compact and $r(r^{-1}(V))=V$. 
By replacing it by a smaller set if necessary, 
we may assume $V\subset F^0$ and 
$V\cap r(E^1\setminus F^1)=\emptyset$. 
Then $V\subset F^0$ is a neighborhood of $v\in F^0$ 
such that $(r_F)^{-1}(V)\subset F^1$ is compact 
and $r_F((r_F)^{-1}(V))=V$. 
Thus we have $v\in \s{F}{0}{rg}\setminus \overline{r(E^1\setminus F^1)}$. 
Therefore we get 
$$F^0\cap\s{E}{0}{rg}\setminus \overline{r(E^1\setminus F^1)}
=\s{F}{0}{rg}\setminus \overline{r(E^1\setminus F^1)}.$$
This completes the proof. 
\end{proof}

\begin{proposition}\label{FY}
For a subgraph $F=(F^0,F^1,d_F,r_F)$ 
of the topological graph $E$, 
the $C^*$-subalgebra of $\cO(E)$ 
generated by $t^0(C_0(F^0))$ and $t^1(C_{d_F}(F^1))$ 
is isomorphic to $\cO(F_Y)$ where 
$Y$ is the closed subset of $\s{F}{0}{rg}$ defined 
by $Y=\s{F}{0}{rg}\cap \overline{r(E^1\setminus F^1)}$. 
\end{proposition}

\begin{proof}
The injective Toeplitz $F$-pair $T=(T^0,T^1)$ 
admits a gauge action and 
$C^*(T)$ is the $C^*$-subalgebra of $\cO(E)$ 
generated by $t^0(C_0(F^0))$ and $t^1(C_{d_F}(F^1))$. 
By Lemma \ref{lemFY}, 
we have $Y_{T}=Y$. 
Hence by Corollary \ref{CT=OE}, 
the $C^*$-subalgebra of $\cO(E)$ 
generated by $t^0(C_0(F^0))$ and $t^1(C_{d_F}(F^1))$ 
is isomorphic to $\cO(F_Y)$.
\end{proof}

If a subgraph $F$ of $E$ 
satisfies $\s{F}{0}{rg}\cap r(E^1\setminus F^1)=\emptyset$, 
then the $C^*$-subalgebra of $\cO(E)$ 
generated by $t^0(C_0(F^0))$ and $t^1(C_{d_F}(F^1))$ 
is isomorphic to $\cO(F)$ by Proposition \ref{FY}. 
The following is one useful construction of such subgraphs. 

For an open subset $V$ of $E^0$, 
we define $F_V^0=V\cup d(r^{-1}(V))\subset E^0$ 
and $F_V^1=r^{-1}(V)\subset E^1$. 
We see that $F_V^0,F_V^1$ are open subsets, 
and that $d(F_V^1),r(F_V^1)\subset F_V^0$. 
Hence $F_V=(F_V^0,F_V^1,d_V,r_V)$ 
is a subgraph of the topological graph $E$, 
where $d_V,r_V$ are the restrictions of $d,r$ to $F_V^1$. 
We denote by $A_V\subset\cO(E)$ 
the $C^*$-subalgebra of $\cO(E)$ generated by $t^0(C_0(F_V^0))$ 
and $t^1(C_{d_V}(F_V^1))$. 
We will show that $A_V\cong\cO(F_V)$. 

\begin{lemma}\label{lemAV}
We have $(F_V^0)_{\rs{rg}}=V\cap \s{E}{0}{rg}$. 
\end{lemma}

\begin{proof}
Take $v\in (F_V^0)_{\rs{rg}}$. 
By Lemma \ref{E0r},
there exists a neighborhood $W$ of $v\in F_V^0$ 
such that $r_V^{-1}(W)\subset F_V^1$ is compact 
and $r_V(r_V^{-1}(W))=W$. 
We can choose such a $W$ with $W\subset V$
because $v\in (F_V^0)_{\rs{rg}}\subset r(F_V^1)\subset V$. 
Hence we have $r_V^{-1}(W)=r^{-1}(W)$. 
This implies that $v\in \s{E}{0}{rg}$ by Lemma \ref{E0r}. 
Thus $(F_V^0)_{\rs{rg}}\subset V\cap \s{E}{0}{rg}$. 

Conversely take $v\in V\cap \s{E}{0}{rg}$. 
Then we can find a neighborhood $W$ of $v\in E^0$ 
such that $W\subset V$, $r^{-1}(W)\subset E^1$ is compact 
and $r(r^{-1}(W))=W$ 
by Lemma \ref{E0r}. 
This implies that $v\in (F_V^0)_{\rs{rg}}$. 
Hence $(F_V^0)_{\rs{rg}}\supset V\cap \s{E}{0}{rg}$. 
Therefore we have $(F_V^0)_{\rs{rg}}=V\cap \s{E}{0}{rg}$. 
\end{proof}

\begin{proposition}\label{AV=OFV}
For an open subset $V$ of $E^0$, 
we have $A_V\cong\cO(F_V)$. 
\end{proposition}

\begin{proof}
By Proposition \ref{FY}, 
it suffices to check 
$(F_V^0)_{\rs{rg}}\cap r(E^1\setminus F^1)=\emptyset$, 
which easily follows from Lemma \ref{lemAV}. 
\end{proof}

\begin{proposition}\label{indlim}
If an increasing family of open subsets 
$V_1\subset V_2\subset\cdots\subset V_n\subset \cdots$ of $E^0$ 
satisfies $\bigcup_{k=1}^\infty V_k=E^0$, 
then we have $\cO(E)=\overline{\bigcup_{k=1}^\infty A_{V_k}}$ 
and $\cO(E)\cong\lim_{k\to\infty}\cO(F_{V_k})$. 
\end{proposition}

\begin{proof}
It is easy to see that 
we have $A_{V_1}\subset A_{V_2}$ 
for two open subsets $V_1,V_2$ of $E^0$ satisfying $V_1\subset V_2$. 
Hence $\overline{\bigcup_{k=1}^\infty A_{V_k}}$ 
is a $C^*$-subalgebra of $\cO(E)$ which contains 
$t^0(C_0(E^0))$ and $t^1(C_d(E^1))$ 
because $\bigcup_{k=1}^\infty V_k=E^0$. 
Hence we have 
$\overline{\bigcup_{k=1}^\infty A_{V_k}}=\cO(E)$. 
By Proposition \ref{AV=OFV}, 
we get $\cO(E)\cong\lim_{k\to\infty}\cO(F_{V_k})$. 
\end{proof}

\begin{remark}\label{Remsubalg}
For two open subsets $V_1,V_2$ of $E^0$ satisfying $V_1\subset V_2$, 
the injective map $\cO(F_{V_1})\to\cO(F_{V_2})$ 
induced by the inclusion $A_{V_1}\subset A_{V_2}$ 
is the same map 
as the one induced by a regular factor map $m=(m^0,m^1)$ 
from $F_{V_2}$ to $F_{V_1}$ 
defined by 
$m^i|_{F_{V_1}^i}=\id_{F_{V_1}^i}$ 
and $m^i(F_{V_2}^i\setminus F_{V_1}^i)=\infty$ 
for $i=0,1$. 
Thus the latter statement of Proposition \ref{indlim} is 
the special case of Proposition \ref{projisom}. 
\end{remark}

We determine conditions on an open subset $V$ 
that imply $A_V$ is hereditary or full. 

\begin{lemma}\label{lemhere}
For an open subset $V$ of $E^0$ with $d(r^{-1}(V))\subset V$, 
we have $F_V^n=(r^n)^{-1}(V)$ for $n\in\N$. 
\end{lemma}

\begin{proof}
For $n=0$, we have $F_V^0=V\cup d(r^{-1}(V))=V$. 
For $n=1$, by definition $F_V^1=r^{-1}(V)$. 
Let $n$ be an integer greater than $1$. 
Since $F_V^1=r^{-1}(V)$, 
we have $F_V^n\subset (r^n)^{-1}(V)$. 
Take $e=(e_1,e_2,\ldots,e_n)\in (r^n)^{-1}(V)$. 
From $r(e_1)=r^n(e)\in V$, 
we get $e_1\in r^{-1}(V)=F_V^1$. 
Then we have $r(e_2)=d(e_1)\in d(r^{-1}(V))\subset V$. 
Hence we get $e_2\in r^{-1}(V)=F_V^1$. 
Recursively, we have $e_i\in r^{-1}(V)=F_V^1$ 
for $i=1,2,\ldots,n$. 
Hence we get $e\in F_V^n$. 
Therefore we have $F_V^n=(r^n)^{-1}(V)$ for $n\in\N$. 
\end{proof}

\begin{proposition}\label{herfull}
If an open subset $V$ of $E^0$ satisfies $d(r^{-1}(V))\subset V$, 
then $A_V$ is a hereditary subalgebra of $\cO(E)$. 
If in addition each $v\in E^0\setminus V$ is regular 
and satisfies $d^n((r^n)^{-1}(v))\subset V$ for some $n\in\N$, 
then the hereditary subalgebra $A_V$ is full in $\cO(E)$.
Hence in this case, 
$\cO(F_V)$ is strongly Morita equivalent to $\cO(E)$. 
\end{proposition}

\begin{proof}
Let $V$ be an open subset of $E^0$ 
such that $d(r^{-1}(V))\subset V$. 
The linear span of elements 
of the form $t^n(\xi)t^m(\eta)^*$ 
is dense in $\cO(E)$ 
(see \cite[Section 2]{Ka1} for the proof of this fact 
and the definition of 
the linear map $t^n\colon C_d(E^n)\to \cO(E)$).
Therefore, once we get $t^0(f)t^n(\xi)\in A_V$ 
for arbitrary $f\in C_0(V)$, $\xi\in C_d(E^n)$ and $n\in\N$. 
we know that $A_V$ is a hereditary subalgebra 
generated by $t^0(C_0(V))$. 
Take $f\in C_0(V)$ and $\xi\in C_d(E^n)$ for $n\in\N$. 
For $e\notin F_V^n$, 
we have $(\pi_{r^n}(f)\xi)(e)=f(r^n(e))\xi(e)=0$ 
because $r^n(e)\notin V$ by Lemma \ref{lemhere}. 
Hence 
we have $\pi_{r^n}(f)\xi\in C_{d_V}(F_V^n)$. 
Thus we get 
$$t^0(f)t^n(\xi)=t^n(\pi_{r^n}(f)\xi)\in 
t^n(C_{d_V}(F_V^n))\subset A_V.$$ 
This proves the first part. 

Now further assume that each $v\in E^0\setminus V$ is regular 
and satisfies $d^n((r^n)^{-1}(v))\subset V$ for some $n\in\N$. 
We set $V_0=V$, 
and define $V_n\subset E^0$ for $n=1,2,\ldots$ by 
$$V_n=\{v\in E^0\mid d^n((r^n)^{-1}(v))\subset V\}.$$
Since $V$ satisfies $d(r^{-1}(V))\subset V$, 
we have $V_n\subset V_{n+1}$ for $n\in\N$. 
By the assumption, we have $\bigcup_{n=0}^\infty V_n=E^0$. 
For $n\in\N$, we have 
$$V_{n+1}=\{v\in E^0\mid r^{-1}(v)\subset d^{-1}(V_{n})\}.$$
If we set 
$$V_{n+1}'=V_{n+1}\cap \s{E}{0}{rg}
=\{v\in \s{E}{0}{rg}\mid r^{-1}(v)\subset d^{-1}(V_{n})\},$$
then we have $V_{n+1}=V\cup V_{n+1}'$ 
because $\s{E}{0}{sg}\subset V$. 
By \cite[Lemma 1.21]{Ka1}, 
if $V_{n}$ is open then $V_{n+1}'$ is open.
Hence we can show that $V_n$ is open recursively. 
Let $I$ be an ideal generated by $A_V$. 
We will show $t^0(C_0(V_n))\subset I$ 
by induction with respect to $n\in\N$. 
For $n=0$, 
we have $t^0(C_0(V_0))=t^0(C_0(V))\subset A_V\subset I$. 
Assume we have $t^0(C_0(V_n))\subset I$ for $n\in\N$. 
Take $f\in C_0(V_{n+1}')$. 
Since $V_{n+1}'\subset \s{E}{0}{rg}$, 
we have $t^0(f)=\varphi(\pi_r(f))$. 
By $r^{-1}(V_{n+1}')\subset d^{-1}(V_n)$, 
we see $\pi_r(f)\in C_0(d^{-1}(V_n))\subset \cK(C_d(d^{-1}(V_n)))$. 
Since $t^0(C_0(V_n))\subset I$, 
we have $t^1(C_d(d^{-1}(V_n)))\subset I$. 
Hence $\varphi\big(\cK(C_d(d^{-1}(V_n)))\big)\subset I$. 
This shows that $t^0(f)\in I$. 
Hence we have $t^0(C_0(V_{n+1}))\subset I$ 
because $V_{n+1}=V\cup V_{n+1}'$. 
Thus we have shown that 
$t^0(C_0(V_n))\subset I$ for all $n\in\N$. 
Since $\bigcup_{n=0}^\infty V_n=E^0$, 
we get $t^0(C_0(E^0))\subset I$. 
Since $\cO(E)$ is generated by $t^0(C_0(E^0))$ as a hereditary subalgebra, 
we have $I=\cO(E)$. 
Thus $A_V$ is a full hereditary subalgebra. 

The last part follows from Proposition \ref{AV=OFV}.
\end{proof}

\begin{remark}
We will see in \cite[Remark 6.2]{Ka2} that 
for an open subset $V$ of $E^0$ 
with $d(r^{-1}(V))\subset V$, 
the condition that 
each $v\in E^0\setminus V$ is regular 
and satisfies $d^n((r^n)^{-1}(v))\subset V$ for some $n\in\N$ 
is not only sufficient but also necessary 
for $A_V$ to be full. 
\end{remark}

\section{Strong Morita equivalence}\label{SecMorita}

In Proposition \ref{herfull}, 
we found subgraphs $F$ of $E$ 
such that $\cO(F)$ is strongly Morita equivalent to $\cO(E)$. 
In this section, 
we give a construction of a topological graph $F$ 
which contains a given topological graph $E$ 
such that $\cO(F)$ is strongly Morita equivalent to $\cO(E)$. 

Let $E=(E^0,E^1,d,r)$ be a topological graph, 
and $N$ be a positive integer or $\infty$. 
For $k=1,2,\ldots,N$, 
take locally compact spaces $X_k^0,X_k^1$, 
local homeomorphisms $d_k\colon X_k^1\to X_{k-1}^0$, 
and proper continuous surjections $r_k\colon X_k^1\to X_k^0$, 
where $X_0^0=E^0$. 
Set 
$$F^0=E^0\amalg X_1^0\amalg\cdots\amalg X_N^0,\quad 
F^1=E^1\amalg X_1^1\amalg\cdots\amalg X_N^1$$ 
and define $d_F,r_F\colon F^1\to F^0$ 
from $d,r$ and $d_k,r_k$ for $k=1,2,\ldots,N$. 
Then $F=(F^0,F^1,d_F,r_F)$ is a topological graph. 

\begin{lemma}\label{fuco}
We have $\s{F}{0}{rg}=\s{E}{0}{rg}\amalg X_1^0\amalg\cdots\amalg X_N^0$. 
\end{lemma}

\begin{proof}
Clearly, $\s{F}{0}{rg}\cap E^0=\s{E}{0}{rg}$.
For $k=1,2,\ldots,N$, 
we have $X_k^0\subset \s{F}{0}{rg}$
because $r_k$ is surjective and proper. 
Hence we have 
$\s{F}{0}{rg}=\s{E}{0}{rg}\amalg X_1^0\amalg\cdots\amalg X_N^0$. 
\end{proof}

We have 
\begin{align*}
C_0(F^0)&=C_0(E^0)\oplus C_0(X_1^0)\oplus\cdots\oplus C_0(X_N^0),\\ 
C_{d_F}(F^0)&=C_d(E^1)\oplus C_{d_1}(X_1^1)
\oplus\cdots\oplus C_{d_N}(X_N^1). 
\end{align*}
Hence there exist natural inclusions
$\mu^0\colon C_0(E^0)\to C_0(F^0)$ 
and $\mu^1\colon C_d(E^1)\to C_{d_F}(F^1)$. 
Let $t_E=(t_E^0,t_E^1)$ and $t_F=(t_F^0,t_F^1)$ 
be the universal Cuntz Krieger $E$-pair in $\cO(E)$ 
and the universal Cuntz Krieger $F$-pair in $\cO(F)$, 
respectively. 

\begin{proposition}\label{embed}
There exists an injective $*$-ho\-mo\-mor\-phism 
$\mu\colon \cO(E)\to\cO(F)$ 
such that $t_F^i\circ\mu^i=\mu\circ t_E^i$ for $i=0,1$. 
\end{proposition}

\begin{proof}
It is routine to check that 
the pair $T=(t_F^0\circ\mu^0,t_F^1\circ\mu^1)$ is a Toeplitz $E$-pair. 
By Lemma \ref{fuco}, 
$T$ is a Cuntz-Krieger $E$-pair. 
Hence there exists a $*$-ho\-mo\-mor\-phism $\mu\colon \cO(E)\to\cO(F)$ 
such that $t_F^i\circ\mu^i=\mu\circ t_E^i$ for $i=0,1$. 
Clearly the pair $T$ is injective 
and admits a gauge action. 
Hence $\mu$ is injective by Proposition \ref{GIUT}, 
\end{proof}

\begin{remark}
As in Remark \ref{Remsubalg}, 
we see that 
the injection $\mu\colon \cO(E)\to\cO(F)$ in Proposition \ref{embed} 
is induced by a regular factor map from $F$ to $E$. 
\end{remark}

\begin{proposition}\label{fulcor}
The image $\mu(\cO(E))$ is a full corner of $\cO(F)$. 
Hence $\cO(F)$ is strongly Morita equivalent to $\cO(E)$. 
\end{proposition}

\begin{proof}
Since $t_F^0\colon C_0(F^0)\to \cO(F)$ is non-degenerate 
by \cite[Lemma 1.20]{Ka1}, 
it extends the map $C_b(F^0)=\cM(C_0(F^0))\to\cM(\cO(F))$ 
where $\cM(\cdot)$ means a multiplier algebra 
(see \cite[Proposition 2.1]{La}). 
Let $p\in \cM(\cO(F))$ be 
the image of the characteristic function of $E^0\subset F^0$ 
under this map. 
We will show that $p\cO(F)p=\mu(\cO(E))$. 
Since $\cO(F)$ is the linear span of elements of the form 
$t_F^n(\xi)t_F^m(\eta)^*$ for $\xi\in C_{d_F}(F^n)$, 
$\eta\in C_{d_F}(F^m)$ and $n,m\in\N$, 
the corner $p\cO(F)p$ is the linear span of elements of the form 
$pt_F^n(\xi)t_F^m(\eta)^*p$ 
(see \cite[Section 2]{Ka1}). 
For $\xi\in C_{d_F}(F^n)$, 
we have $pt_F^n(\xi)=\mu(t_E^n(\xi_0))$ 
where $\xi_0\in C_d(E^n)$ 
is the restriction of $\xi$ to $E^n\subset F^n$, 
because for $e\in F^n$, $r_F^n(e)\in E^0$ implies $e\in E^n$. 
Thus we have $p\cO(F)p=\mu(\cO(E))$. 

Let $I$ be an ideal of $\cO(F)$ generated by $p\cO(F)p$.
We will prove that $I=\cO(F)$ 
in a fashion that is similar to the proof of Proposition \ref{herfull}. 
By the former part, 
we have $t_F^0(C_0(E^0)),t_F^1(C_d(E^1))\subset I$. 
We will prove $t_F^0(C_0(X_1^0)),t_F^1(C_{d_1}(X_1^1))\subset I$. 
Since $d_1(X_1^1)\subset E^0$, 
we see that $t_F^1(\xi)p=t_F^1(\xi)$ 
for every $\xi\in C_{d_1}(X_1^1)\subset C_{d_F}(F^1)$. 
Therefore we get $t_F^1(\xi)\in I$ 
for every $\xi\in C_{d_1}(X_1^1)$. 
This implies that $\varphi_F\big(\cK(C_{d_1}(X_1^1))\big)\subset I$. 
For $f\in C_0(X_1^0)$, 
we have $\pi_{r_F}(f)\in \cK(C_{d_1}(X_1^1))$, 
and so $t_F^0(f)=\varphi_F(\pi_{r_F}(f))\in I$ 
because $X_1^0\subset\s{F}{0}{rg}$. 
Thus we get $t_F^0(C_0(X_1^0)),t_F^1(C_{d_1}(X_1^1))\subset I$. 
The same argument shows that 
$t_F^0(C_0(X_k^0)),t_F^1(C_{d_k}(X_k^1))\subset I$ 
for $k=2,\ldots,N$ by induction. 
Hence we get $t_F^0(C_0(F^0)),t_F^1(C_{d_F}(F^1))\subset I$. 
Since $\cO(F)$ is generated by the images of $t_F^0$ and $t_F^1$, 
we have $I=\cO(F)$. 
Therefore $\mu(\cO(E))$ is a full corner of $\cO(F)$. 
The last part follows from the injectivity of $\mu$. 
\end{proof}

Now we specialize the above discussion. 
Let $E=(E^0,E^1,d,r)$ be a topological graph 
and $N$ be a positive integer. 
For $k=1,2,\ldots,N$, 
set $X_k^0,X_k^1\cong E^0$, 
and define $d_k\colon X_k^1\to X_{k-1}^0$ 
and $r_k\colon X_k^1\to X_k^0$ by the identity map of $E^0$ 
where $X_0^0=E^0$. 
Set 
$$E_N^0=E^0\amalg X_1^0\amalg\cdots\amalg X_{N}^0,\quad 
E_N^1=E^1\amalg X_1^1\amalg\cdots\amalg X_{N}^1$$ 
and define $d_N,r_N\colon E_N^1\to E_N^0$ 
from $d,r$ and $d_k,r_k$ ($k=1,2,\ldots,N$). 
Then $E_N=(E_N^0,E_N^1,d_N,r_N)$ is a topological graph. 
By Proposition \ref{fulcor}, 
the $C^*$-algebra $\cO(E_N)$ 
is strongly Morita equivalent to $\cO(E)$. 
We can say more. 

\begin{proposition}\label{otimesMn}
We have $\cO(E_N)\cong \cO(E)\otimes\M_{N+1}$.
\end{proposition}

\begin{proof}
In the proof, 
we use the natural identification 
\begin{align*}
C_0(E_N^0)&=C_0(E^0)\oplus C_0(X_1^0)\oplus\cdots\oplus C_0(X_{N}^0),\\ 
C_{d_N}(E_N^1)&=C_d(E^1)\oplus 
C_{d_1}(X_1^1)\oplus\cdots\oplus C_{d_{N}}(X_{N}^1). 
\end{align*}
We consider $C_0(E^0), C_0(X_k^0)$ as subalgebras of $C_0(E_N^0)$, 
and $C_d(E^1), C_{d_k}(X_k^1)$ as subspaces of $C_{d_1}(X_1^1)$ 
for $k=1,2,\ldots,N$. 
We will identify $C_0(X_k^0)$ and $C_{d_k}(X_k^1)$ with $C_0(E^0)$ 
for each $k$. 
Then the map $\theta_{g,g'}\mapsto g\overline{g'}$ 
for $g,g'\in C_{d_k}(X_k^1)\cong C_0(E^0)$ 
gives an isomorphism from $\cK(C_{d_k}(X_k^1))$ to $C_0(E^0)$. 
Since $r_k\colon X_k^1\to X_k^0$ is defined by the identity map of $E^0$, 
the map $\pi_{r_k}\colon C_0(X_k^0)\to \cK(C_{d_k}(X_k^1))$ 
is the inverse of the above isomorphism 
modulo the identification $C_0(X_k^0)\cong C_0(E^0)$. 
Let us denote by $\{u_{k,l}\}_{0\leq k,l\leq N}$ 
the matrix units of $\M_{N+1}$. 
We define two maps $T^0\colon C_0(E_N^0)\to \cO(E)\otimes\M_{N+1}$ and 
$T^1\colon C_{d_N}(E_N^1)\to \cO(E)\otimes\M_{N+1}$ by 
\begin{align*}
T^0(f_0,f_1,\ldots,f_{N})
&=\sum_{k=0}^{N} t^0(f_k)\otimes u_{k,k},\\
T^1(\xi,g_1,\ldots,g_{N})
&=t^1(\xi)\otimes u_{0,0}+\sum_{k=1}^{N}t^0(g_k)\otimes u_{k,k-1}.
\end{align*}
We will show that $T=(T^0,T^1)$ is an injective Cuntz-Krieger $E_N$-pair.
Take $\eta=(\xi,g_1,\ldots,g_{N}),
\eta'=(\xi',g_1',\ldots,g_{N}')\in C_{d_N}(E_N^1)$.
We see that 
$$\ip{\eta}{\eta'}
=\big(\ip{\xi}{\xi'}+\overline{g_1}g_1',
\overline{g_2}g_2',\ldots,\overline{g_{N}}g_{N}',0\big)
\in C_0(E_N^0).$$
We have 
\begin{align*}
T^1&(\eta)^*T^1(\eta')\\
&=\bigg(t^1(\xi)\otimes u_{0,0}
  +\sum_{k=1}^{N}t^0(g_k)\otimes u_{k,k-1}\bigg)^*
  \bigg(t^1(\xi')\otimes u_{0,0}
  +\sum_{k=1}^{N}t^0(g_k')\otimes u_{k,k-1}\bigg)\\
&=\bigg(t^1(\xi)^*\!\otimes u_{0,0}
  +\sum_{k=1}^{N}t^0(g_k)^*\!\otimes u_{k-1,k}\bigg)
  \bigg(t^1(\xi')\otimes u_{0,0}
  +\sum_{k=1}^{N}t^0(g_k')\otimes u_{k,k-1}\bigg)\\
&=\big(t^1(\xi)^*t^1(\xi')\big)\otimes u_{0,0}
  +\sum_{k=1}^{N}\big(t^0(g_k)^*t^0(g_k')\big)\otimes u_{k-1,k-1}\\
&=t^0(\ip{\xi}{\xi'})\otimes u_{0,0}
  +\sum_{k=1}^{N}t^0(\overline{g_k}g_k')\otimes u_{k-1,k-1}\\
&=T^0(\ip{\eta}{\eta'}).
\end{align*}
Take $f=(f_0,f_1,\ldots,f_{N})\in C_0(E_N^0)$ and 
$\eta=(\xi,g_1,\ldots,g_{N})\in C_{d_N}(E_N^1)$. 
We see that 
$\pi_{r_N}(f)\eta=(\pi_r(f_0)\xi,f_1g_1,\ldots,f_{N} g_{N})$.
We have 
\begin{align*}
T^0(f)T^1(\eta)
&=\bigg(\sum_{k=0}^{N} t^0(f_k)\otimes u_{k,k}\bigg)
  \bigg(t^1(\xi)\otimes u_{0,0}
  +\sum_{k=1}^{N}t^0(g_k)\otimes u_{k,k-1}\bigg)\\
&=\big(t^0(f_0)t^1(\xi)\big)\otimes u_{0,0}
  +\sum_{k=1}^{N}\big(t^0(f_k)t^0(g_k)\big)\otimes u_{k,k-1}\\
&=t^1(\pi_r(f_0)\xi)\otimes u_{0,0}
  +\sum_{k=1}^{N}t^0(f_kg_k)\otimes u_{k,k-1}\\
&=T^1(\pi_{r_N}(f)\eta).
\end{align*}
Thus $T=(T^0,T^1)$ is an injective Toeplitz $E_N$-pair. 
We will prove that $T$ is a Cuntz-Krieger $E_N$-pair. 
Take $f=(f_0,f_1,\ldots,f_{N})\in C_0((E_N^0)_{\rs{rg}})$. 
By Lemma \ref{fuco}, 
we have $f_0\in C_0(\s{E}{0}{rg})$. 
Hence we see that 
\begin{align*}
\pi_{r_N}(f_0,f_1,\ldots,f_{N})
&=\big(\pi_{r}(f_0),\pi_{r_1}(f_1),\ldots,\pi_{r_{N}}(f_{N})\big)\\
&\in\cK(C_d(E^1))\oplus \cK(C_{d_1}(X_1^1))
\oplus\cdots\oplus \cK(C_{d_{N}}(X_{N}^1)).
\end{align*}
We compute 
$\varPhi(\pi_{r_N}(f_0,f_1,\ldots,f_{N}))\in \cO(E)\otimes\M_{N+1}$ 
where $\varPhi\colon \cK(C_{d_N}(E_N^1))\to \cO(E)\otimes\M_{N+1}$ 
is the map induced by the Toeplitz $E_N$-pair $T=(T^0,T^1)$. 
For $\xi,\xi'\in C_d(E^1)\subset C_{d_N}(E_N^1)$, 
we have 
$$\varPhi(\theta_{\xi,\xi'})=T^1(\xi)T^1(\xi')^*
=(t^1(\xi)\otimes u_{0,0})(t^1(\xi')\otimes u_{0,0})^*
=\varphi(\theta_{\xi,\xi'})\otimes u_{0,0}.$$
Hence we get 
$$\varPhi(\pi_{r}(f_0))=\varphi(\pi_{r}(f_0))\otimes u_{0,0}
=t^0(f_0)\otimes u_{0,0}$$ 
because $f_0\in C_0(\s{E}{0}{rg})$. 
For $k=1,2,\ldots,N$, 
we get 
$$\varPhi(\pi_{r_k}(f_k))=t^0(f_{k})\otimes u_{k,k},$$ 
by the remark in the beginning of this proof and the computation 
\begin{align*}
\varPhi(\theta_{g_k,g_k'})
&=T^1(g_k)T^1(g_k')^*\\
&=(t^0(g_k)\otimes u_{k,k-1})(t^0(g_k')\otimes u_{k,k-1})^*\\
&=t^0(g_k\overline{g_k'})\otimes u_{k,k},
\end{align*}
for $g_k,g_k'\in C_{d_k}(X_k^1)\subset C_{d_N}(E_N^1)$. 
Therefore we have
\begin{align*}
\varPhi(\pi_{r_N}(f_0,f_1,\ldots,f_{N}))
&=t^0(f_0)\otimes u_{0,0}+t^0(f_1)\otimes u_{1,1}+\cdots
 +t^0(f_{N})\otimes u_{N,N}\\
&=T^0(f_0,f_1,\ldots,f_{N}).
\end{align*}
This shows that $T$ is a Cuntz-Krieger $E_N$-pair.

Next we will show that the $C^*$-algebra $C^*(T)$ 
generated by the images of $T^0$ and $T^1$ is $\cO(E)\otimes\M_{N+1}$. 
Since $T^0(f)=t^0(f)\otimes u_{0,0}$ 
and $T^1(\xi)=t^1(\xi)\otimes u_{0,0}$ 
for $f\in C_0(E^0)\subset C_0(E_N^0)$ 
and $\xi\in C_d(E^1)\subset C_{d_N}(E_N^1)$, 
we have $\cO(E)\otimes u_{0,0}\subset C^*(T)$. 
For $f\in C_0(X_1^0)\subset C_0(E_N^0)$, 
we have $T^0(f)=t^0(f)\otimes u_{1,0}$. 
Hence $t^0(f)x\otimes u_{1,0}\subset C^*(T)$ 
for all $f\in C_0(E^0)$ and all $x\in\cO(E)$. 
Since we have $t^0(C_0(E^0))\cO(E)=\cO(E)$
(see \cite[Proposition 2.5]{Ka1} and the remark before it), 
we get $\cO(E)\otimes u_{1,0}\subset C^*(T)$. 
Recursively, 
we have $\cO(E)\otimes u_{k,0}\subset C^*(T)$ 
for $k=2,\ldots,N$. 
Since $\cO(E)\otimes\M_{N+1}$ 
is generated by $\cO(E)\otimes u_{k,0}$ for $k=0,1,\ldots,N$, 
we have $C^*(T)=\cO(E)\otimes\M_{N+1}$. 

Finally we will find a gauge action for the pair $T$. 
For each $z\in\T$, define a unitary 
$u_z\in\M_{N+1}$ by 
$u_z=\sum_{k=0}^{N} z^ku_{k,k}$, 
and an automorphism $\Ad u_z$ on $\M_{N+1}$ 
by $\Ad u_z(x)=u_zxu_z^*$ for $x\in\M_{N+1}$.
Let $\beta$ be the gauge action on $\cO(E)$. 
The automorphism $\beta'_z=\beta_z\otimes\Ad u_z$ 
of $\cO(E)\otimes\M_{N+1}$ satisfies the equations  
$\beta'_z(T^0(f))=T^0(f)$ and $\beta'_z(T^1(\eta))=zT^1(\eta)$ 
for $f\in C_0(E_N^0)$ and $\eta\in C_{d_N}(E_N^1)$. 
Now $T$ gives an isomorphism $\cO(E_N)\cong \cO(E)\otimes\M_{N+1}$ 
with the help of Proposition \ref{GIUT}. 
\end{proof}

The discussion above works for $N=\infty$.
Namely, if we define a topological graph 
$E_\infty=(E_\infty^0,E_\infty^1,d_\infty,r_\infty)$
by 
\begin{align*}
E_\infty^0&=E^0\amalg X_1^0\amalg X_2^0\amalg\cdots,& 
E_\infty^1&=E^1\amalg X_1^1\amalg X_2^1\amalg\cdots \quad 
(X_k^0, X_k^1\cong E^0),
\end{align*}
\begin{align*}
d_\infty&\colon E^1\stackrel{d}{\longrightarrow} E^0,\quad 
d_\infty \colon X_k^1\stackrel{\id}{\longrightarrow}X_{k-1}^0&
&\mbox{where }X_0^0=E^0\subset E_\infty^0, \\
r_\infty&\colon E^1\stackrel{r}{\longrightarrow} E^0,\quad 
r_\infty \colon X_k^1\stackrel{\id}{\longrightarrow}X_k^0&
&\mbox{for }k=1,2,\ldots,
\end{align*}
then the proof of Proposition \ref{otimesMn}, 
with appropriate simple modifications, 
proves the following proposition. 

\begin{proposition}
We have $\cO(E_\infty)\cong\cO(E)\otimes\K$.
\end{proposition}

\begin{remark}
This Proposition generalizes \cite[Theorem 4.2]{T}. 
\end{remark}

\begin{remark}
For a topological graph $E$ and a positive integer $N$, 
there are many ways to construct topological graphs $E'$ 
such that the associated $C^*$-algebras $\cO(E')$ 
are isomorphic to $\cO(E)\otimes\M_{N+1}$. 
Besides the topological graph $E_N$ defined above, 
we give another example 
$\bar{E}_N=(\bar{E}_N^0,\bar{E}_N^1,\bar{d}_N,\bar{r}_N)$. 
We set $\bar{E}_N^0=E_N^0$, $\bar{E}_N^1=E_N^1$, 
$\bar{r}_N=r_N$ 
and define $\bar{d}_N\colon \bar{E}_N^1\to \bar{E}_N^0$ by 
$$\bar{d}_N\colon E^1\stackrel{d}{\longrightarrow} E^0,\quad 
\bar{d}_N\colon X_k^1\stackrel{\id}{\longrightarrow}E^0\quad 
\mbox{ for }k=1,2,\ldots,N.$$
We can prove that $\cO(\bar{E}_N)\cong\cO(E)\otimes\M_{N+1}$ 
by using the injective Cuntz-Krieger $\bar{E}_N$-pair 
$\bar{T}=(\bar{T}^0,\bar{T}^1)$ where $\bar{T}^0$ is the same map as $T^0$ 
in Proposition \ref{otimesMn}, and 
$$\bar{T}^1(\xi,g_1,\ldots,g_{N})
=t^1(\xi)\otimes u_{0,0}+\sum_{k=1}^{N}t^0(g_k)\otimes u_{k,0},$$
for $(\xi,g_1,\ldots,g_{N})\in C_{\bar{d}_N}(\bar{E}_N^1)$, 
and also using the automorphism $\beta_z\otimes \Ad u_z'$ 
where $u_z'=u_{0,0}+z\sum_{k=1}^{N} u_{k,k}\in\M_{N+1}$ for $z\in\T$. 
This construction also works for $N=\infty$. 
\end{remark}

\section{Other operations}\label{SecOther}

\begin{proposition}
For a topological graph $E=(E^0,E^1,d,r)$, 
the $C^*$-algebra $\cO(E)$ is unital if and only if $E^0$ is compact. 
\end{proposition}

\begin{proof}
We have that 
$\cO(E)$ is unital if and only if $t^0(C_0(E^0))$ is unital 
because the hereditary subalgebra generated by $t^0(C_0(E^0))$ 
is $\cO(E)$ (see \cite[Proposition 2.5]{Ka1}).
Since $t^0$ is injective, 
$t^0(C_0(E^0))$ is unital if and only if $E^0$ is compact.
We are done.
\end{proof}

\begin{definition}
Let $E=(E^0,E^1,d,r)$ be a topological graph. 
The topological graph $\widetilde{E}=(\widetilde{E}^0,E^1,d,r)$ 
is called the {\rm one-point compactification} of $E$
where $\widetilde{E}^0=E^0\cup\{\infty\}$ 
is the one-point compactification of $E^0$. 
\end{definition}

\begin{lemma}\label{tErg}
For the one-point compactification $\widetilde{E}$ 
of a topological graph $E=(E^0,E^1,d,r)$, 
we have $\s{\widetilde{E}}{0}{rg}=\s{E}{0}{rg}$. 
\end{lemma}

\begin{proof}
It is clear that $\s{\widetilde{E}}{0}{fin}\cap E^0=\s{E}{0}{fin}$ and 
$\s{\widetilde{E}}{0}{sce}\cap E^0=\s{E}{0}{sce}$. 
Therefore $\s{\widetilde{E}}{0}{rg}\cap E^0=\s{E}{0}{rg}$. 
Since $r^{-1}(\infty)=\emptyset$, 
we have $\infty\notin\s{\widetilde{E}}{0}{rg}$ 
by Lemma \ref{E0r}. 
Hence we have $\s{\widetilde{E}}{0}{rg}=\s{E}{0}{rg}$. 
\end{proof}

In Lemma \ref{tErg}, we see that 
$\infty\in\s{\widetilde{E}}{0}{sg}=\s{\widetilde{E}}{0}{inf}
\cup\overline{\s{\widetilde{E}}{0}{sce}}$. 
Note that there exist topological graphs $E$ 
with $\infty\notin\s{\widetilde{E}}{0}{inf}$ 
(for example, in the case that $E^1$ is compact)
as well as ones with $\infty\notin\overline{\s{\widetilde{E}}{0}{sce}}$
(for example, in the case that $E^1=E^0$ and $r=\id$). 

\begin{proposition}
Let $E$ be a topological graph, 
and $\widetilde{E}$ be its one-point compactification. 
Then $\cO\big(\widetilde{E}\big)$ is isomorphic to 
the unitization ${\cO(E)}^{\sim}$ of the $C^*$-algebra $\cO(E)$. 
\end{proposition}

\begin{proof}
Define a $*$-ho\-mo\-mor\-phism 
$\tilde{t}^0\colon C\big(\widetilde{E}^0\big)\to {\cO(E)}^{\sim}$ 
by $\tilde{t}^0(f)=t^0(f-f(\infty))+f(\infty)$. 
Then it is easy to see that $\tilde{t}=(\tilde{t}^0,t^1)$ 
is an injective Toeplitz $\widetilde{E}$-pair 
which admits a gauge action and satisfies $C^*(\tilde{t})={\cO(E)}^{\sim}$. 
By Lemma \ref{tErg}, the pair $\tilde{t}$ is 
a Cuntz-Krieger $\widetilde{E}$-pair. 
Hence by Proposition \ref{GIUT}, 
we have $\cO\big(\widetilde{E}\big)\cong{\cO(E)}^{\sim}$. 
\end{proof}

For a discrete graph $E=(E^0,E^1,d,r)$ 
with infinitely many vertices, 
its one-point compactification $\widetilde{E}$ is no longer discrete. 

\begin{definition}
We define a {\em disjoint union} $E\amalg F$ of 
two topological graphs $E$ and $F$ by
$(E\amalg F)^0=E^0\amalg F^0$, $(E\amalg F)^1=E^1\amalg F^1$ and 
$d,r\colon (E\amalg F)^1\to (E\amalg F)^0$ are natural ones. 
The {\em disjoint union} $\coprod_{\lambda\in\Lambda}E_\lambda$ of 
a family of topological graphs $\{E_\lambda\}_{\lambda\in\Lambda}$
is defined similarly. 
\end{definition}

It is easy to see the following.

\begin{proposition}\label{sum}
For two topological graphs $E$ and $F$, 
we have $\cO(E\amalg F)=\cO(E)\oplus \cO(F)$. 
We also have $\cO(\coprod_{\lambda\in\Lambda}E_\lambda)=
\bigoplus_{\lambda\in\Lambda}\cO(E_\lambda)$ 
for a family of topological graphs 
$\{E_\lambda\}_{\lambda\in\Lambda}$. 
\end{proposition}

Let $E=(E^0,E^1,d_E,r_E)$ be a topological graph, 
and $X$ be a locally compact space. 
We define a topological graph $E\times X$ as follows. 
We set $(E\times X)^0=E^0\times X$, $(E\times X)^1=E^1\times X$, 
and define $d,r\colon (E\times X)^1\to (E\times X)^0$ by 
$d((e,x))=(d_E(e),x)$ and $r((e,x))=(r_E(e),x)$ for $(e,x)\in (E\times X)^1$. 
It is easy to see $\s{(E\times X)}{0}{rg}=\s{E}{0}{rg}\times X$. 

\begin{proposition}\label{otimesCX}
We have $\cO(E\times X)\cong\cO(E)\otimes C_0(X)$. 
\end{proposition}

\begin{proof}
First note that we can identify  
$C_0((E\times X)^0)=C_0(E^0)\otimes C_0(X)$. 
We can also see that $C_{d}((E\times X)^1)$ is isomorphic to 
the completion of the algebraic 
tensor product $C_d(E^1)\odot C_0(X)$. 
We define a $*$-ho\-mo\-mor\-phism $T^0\colon C_0((E\times X)^0)\to \cO(E)\otimes C_0(X)$ 
and a linear map $T^1\colon C_{d}((E\times X)^1)\to \cO(E)\otimes C_0(X)$ 
by $T^0(f\otimes g)=t^0(f)\otimes g$ and 
$T^1(\xi\otimes g)=t^1(\xi)\otimes g$ 
for $f\in C_0(E^0)$, $\xi\in C_d(E^1)$ and $g\in C_0(X)$. 
It is routine to check that the pair $T=(T^0,T^1)$ 
is an injective Cuntz-Krieger $E\times X$-pair admitting a gauge action. 
It is also easy to see that $C^*(T)=\cO(E)\otimes C_0(X)$. 
Hence by Proposition \ref{GIUT}, 
we have an isomorphism from $\cO(E\times X)$ to $\cO(E)\otimes C_0(X)$. 
\end{proof}

\section{Examples 1}\label{SecExample1}

Thanks to the study above, 
we can show that 
the class of $C^*$-algebras arising from topological graphs 
contains all AF-algebras and many AH-algebras. 

Let us take an AF-algebra $A$ and write $A=\varinjlim (A_n,\mu_n)$ 
where $A_n$ is a finite dimensional $C^*$-algebra 
and $\mu_n\colon A_n\to A_{n+1}$ is an injective $*$-ho\-mo\-mor\-phism 
for $n\in\N$. 
We write $A_n=\bigoplus_{i=1}^{i_n}A_{n}^{(i)}$ 
where $A_{n}^{(i)}\cong\M_{k_{n}^{(i)}}$ 
for a positive integer $k_{n}^{(i)}$. 
For each $n\in\N$, 
a map $\mu_n\colon A_n\to A_{n+1}$ 
is characterized (up to unitary equivalence) 
by an $\N$-valued rectangular matrix 
$(\sigma_n^{(i,j)})_{1\leq j\leq i_n, 1\leq i\leq i_{n+1}}$ 
where $\sigma_n^{(i,j)}$ is the multiplicity of 
the map $\mu_n^{(i,j)}\colon A_{n}^{(j)}\to A_{n+1}^{(i)}$ 
which is obtained by restricting $\mu_n$. 
Note that we have 
$\sum_{j=1}^{i_n}\sigma_n^{(i,j)}k_{n}^{(j)}\leq k_{n+1}^{(i)}$ 
for each $n\in\N$ and $1\leq i\leq i_{n+1}$. 

For each $n\in\N$, 
we define a topological graph $E_n=(E_n^0,E_n^1,d_n,r_n)$ 
as follows:
\begin{align*}
&E_n^0=\{v_n^{(i,k)}\mid 1\leq i\leq i_n,\ 1\leq k\leq k_{n}^{(i)}\},\\
&E_n^1=\{e_n^{(i,k)}\mid 1\leq i\leq i_n,\ 1\leq k\leq k_{n}^{(i)}-1\},\\
&d_n(e_n^{(i,k)})=v_n^{(i,k)},\quad r_n(e_n^{(i,k)})=v_n^{(i,k+1)}.
\end{align*}
We see that 
$\cO(E_n)\cong\bigoplus_{i=1}^{i_n}\M_{k_{n}^{(i)}}\cong A_n$ 
by Proposition \ref{otimesMn} and Proposition \ref{sum}. 
We define two maps $m_n^0\colon \widetilde{E}_{n+1}^0\to \widetilde{E}_n^0$
and $m_n^1\colon \widetilde{E}_{n+1}^1\to \widetilde{E}_n^1$ as follows. 
Take $i\in\{1,\ldots,i_{n+1}\}$ and $k\in\{1,\ldots, k_{n+1}^{(i)}\}$. 
If $k>\sum_{j=1}^{i_n}\sigma_n^{(i,j)}k_{n}^{(j)}$, 
then define $m_n^0(v_{n+1}^{(i,k)})=\infty$ 
and $m_n^1(e_{n+1}^{(i,k)})=\infty$. 
Otherwise, we can find $j\in\{1,\ldots,i_n\}$ 
such that $k'=k-\sum_{j'=1}^{j-1}\sigma_n^{(i,j')}k_{n}^{(j')}$ satisfies 
that $1\leq k'\leq \sigma_n^{(i,j)}k_{n}^{(j)}$. 
Take $l\in\{1,\ldots, k_{n}^{(j)}\}$ 
with $l\equiv k'\bmod k_{n}^{(j)}$. 
We define $m_n^0(v_{n+1}^{(i,k)})=v_n^{(j,l)}$ 
and 
$$m_n^1(e_{n+1}^{(i,k)})=\left\{\begin{array}{ll}
e_n^{(j,l)}& \mbox{if }1\leq l\leq k_{n}^{(j)}-1\\
\infty& \mbox{if }l=k_{n}^{(j)}. 
\end{array}\right.$$
Then $m_n=(m_n^0,m_n^1)$ is a regular factor map from $E_{n+1}$ to $E_n$ 
and the $*$-ho\-mo\-mor\-phism $\cO(E_n)\to\cO(E_{n+1})$ induced by $m_n$ 
is the same map as the injection $\mu_n$. 
Hence if we denote by $E=(E^0,E^1,d,r)$ 
the projective limit of the projective system $\{E_n\}$ and $\{m_n\}$, 
then we have $A\cong\cO(E)$ by Proposition \ref{projisom}.
Note that $E^0$ is a totally disconnected space, 
and $d,r\colon E^1\to E^0$ are homeomorphisms 
onto open subsets of $E^0$. 
Thus this is an example of a crossed product 
by partial homeomorphisms explained in Subsection \ref{SecPH}, 
and this construction is the same as in \cite{Ex2}. 

Note that 
the class of graph algebras contains all AF-algebras 
up to strong Morita equivalence, 
but does not contain many AF-algebras 
such as simple unital infinite dimensional AF-algebras 
(e.g.\ UHF-algebras). 

Next we see that many AH-algebras can be obtained as 
$C^*$-algebras of topological graphs. 
For a topological graph $E=(E^0,E^1,d,r)$ with 
$E^1=\emptyset$, we have $\cO(E)\cong C_0(E^0)$. 
Thus the class of our algebras contains 
all commutative $C^*$-algebras. 
Combining this fact with 
Proposition \ref{otimesMn} and Proposition \ref{sum}, 
we see that a $C^*$-algebra $A$ of the form 
$A=\bigoplus_{k=1}^K C_0(X_k)\otimes\M_{n_k}$, 
where $X_k$ is a locally compact space and $n_k$ is a positive integer, 
is obtained as a $C^*$-algebra $\cO(E)$ arising 
from a topological graph $E$ 
(we can also use Proposition \ref{otimesCX}). 

Let us take $C^*$-algebras $A,B$ of the form 
$$A=\bigoplus_{k=1}^K C_0(X_k)\otimes\M_{n_k},\quad 
B=\bigoplus_{l=1}^L C_0(Y_l)\otimes\M_{m_l},$$ 
and choose topological graphs $E,F$ 
such that $\cO(E)\cong A$ and $\cO(F)\cong B$ as above. 
Not all $*$-ho\-mo\-mor\-phisms from $A$ to $B$ 
come from regular factor maps from $F$ to $E$. 
However, every diagonal $*$-ho\-mo\-mor\-phism from $A$ to $B$ 
comes from a regular factor map from $F$ to $E$, 
where a $*$-ho\-mo\-mor\-phism $\mu\colon A\to B$ 
is called diagonal 
if for each $k\in\{1,\ldots,K\}$ and $l\in\{1,\ldots,L\}$, 
the restriction 
$\mu_{l,k}\colon C_0(X_k)\otimes\M_{n_k}\to C_0(Y_l)\otimes\M_{m_l}$ 
of $\mu$ is of the form 
$$\mu_{l,k}(f)
=\diag\big\{0,\ldots,0,f\circ m_{l,k}^{(1)},\ldots,
f\circ m_{l,k}^{(\sigma_{l,k})},0,\ldots,0\big\}
\in C_0(Y_l)\otimes\M_{m_l}$$
for $f\in C_0(X_k)\otimes\M_{n_k}$ 
where the $m_{l,k}^{(i)}$ are continuous maps 
from $\widetilde{Y}_l$ to $\widetilde{X}_k$ 
preserving $\infty$ for $i=1,2,\ldots,\sigma_{l,k}$. 
By Proposition \ref{projisom}, 
our class includes all $C^*$-algebras which are given 
by inductive limits of $C^*$-algebras of the form 
$\bigoplus_{k=1}^K C_0(X_k)\otimes\M_{n_k}$ 
with diagonal $*$-ho\-mo\-mor\-phisms. 
In particular, 
all simple real rank zero AT-algebras 
and all Goodearl algebras 
appear as $C^*$-algebras of topological graphs
(see, \cite[Theorem 4.7.5]{Li} and \cite[Example 3.1.7]{RS}). 
The $C^*$-algebras ${\mathcal A}_T$ 
of totally ordered, compact metrizable sets $T$ 
defined in \cite{Ro} 
are also in our class. 
In particular, 
the example ${\mathcal A}_{[0,1]}$ of a purely infinite AH-algebra 
constructed in \cite{Ro} 
is obtained as the $C^*$-algebra of a topological graph. 

Besides AF-algebras and AH-algebras, 
many nuclear $C^*$-algebras satisfying 
the Universal Coefficient Theorem 
appear as $C^*$-algebras of topological graphs, 
for example purely infinite $C^*$-algebras (see \cite{Ka3})
and stabely projectionless $C^*$-algebras. 
In \cite{Ka5}, 
we study the $C^*$-algebras of topological graphs 
arising from constant maps 
in order to analyze the $C^*$-algebras generated by scaling elements.

\section{Examples 2}\label{SecExample2}

Our construction of $C^*$-algebras from topological graphs 
is motivated by graph algebras. 
Graph algebras are one of the generalization of 
Cuntz-Krieger algebras. 
In this section, 
we observe that our construction encompasses 
other generalizations of Cuntz-Krieger algebras. 

\subsection{Exel-Laca algebras}

In \cite{EL}, 
R. Exel and M. Laca presented a method for constructing $C^*$-algebras, 
now called Exel-Laca algebras, 
from an infinite matrix with entries in $\{0,1\}$. 
They introduced these algebras 
in order to extend the work of Cuntz and Krieger 
who focused primarily on finite $\{0,1\}$-valued matrices. 
In \cite[Subsection 3.5]{Sc}, 
Schweizer observed how to present an Exel-Laca algebra 
as the $C^*$-algebra of a topological graph. 

\begin{remark}
Schweizer \cite{Sc} called a topological graph a continuous diagram. 
His presentation of an Exel-Laca algebra, 
given by a matrix $A$, in terms of a topological graph 
was made under the assumption that 
$A$ has no columns that are identically zero. 
This assumption can be removed if, 
in the notation of \cite[Subsection 3.5]{Sc}, 
one adds the characteristic function $\delta_{i}$ to $B$ 
for each $i\in\G$ 
such that the $i^{\text{th}}$ column of $A$ is identically zero. 
\end{remark}

\subsection{Matsumoto algebras}

In \cite{M1}, 
K. Matsumoto introduced a method for constructing $C^*$-algebras, 
now called {\em Matsumoto algebras}, from subshifts. 
When the subshift is a topological Markov shift, 
then his construction coincides with the construction of Cuntz and Krieger. 
In \cite{M2}, 
he generalized subshifts by introducing the notion of a $\lambda$-graph 
and he showed how one can attach a $C^*$-algebra to one of these, 
generalizing the Matsumoto algebras from \cite{M1} 
(see \cite[Corollary 4.5]{M2}). 
His construction used topological graphs, 
which he called continuous graphs. 
Thus, Matsumoto algebras and $\lambda$-graph algebras 
are all instances $C^*$-algebras associated to topological graphs. 

\begin{remark}
In \cite{CM}, 
an alternate construction of $C^*$-algebras from subshifts 
is presented. 
These are based on $\lambda$-graphs and so, ultimately, 
may be viewed as coming from topological graphs. 
\end{remark}

\section{Examples 3}\label{SecExample3}

The class of $C^*$-algebras of topological graphs 
contains the ones of graph algebras and of homeomorphism $C^*$-algebras. 
Graph algebras generalizes Cuntz-Krieger algebras, 
and we study two other such classes in the previous section. 
In this section, 
we study three classes of $C^*$-algebras 
generalizing homeomorphism $C^*$-algebras. 

\subsection{Crossed products by partial homeomorphisms}\label{SecPH}

\begin{definition}
Let $X$ be a locally compact space. 
A {\em partial homeomorphism} is a homeomorphism $\sigma$ 
from an open subset $U$ of $X$ to another open subset $V$ of $X$. 
\end{definition}

If a partial homeomorphism $\sigma\colon U\to V$ on $X$ is given, 
we can define a $*$-ho\-mo\-mor\-phism $\theta\colon C_0(U)\to C_0(V)$ 
by $\theta(f)=f\circ\sigma^{-1}$. 
The triple $(\theta,C_0(V),C_0(U))$ is called 
a {\em partial automorphism} of $C_0(X)$ in \cite{Ex1}. 
R. Exel associated a $C^*$-algebra $C_0(X)\rtimes_{\theta}\Z$ 
with the partial automorphism $(\theta,C_0(V),C_0(U))$ 
\cite[Definition 3.7]{Ex1}. 
Instead of giving a definition of 
the $C^*$-algebra $C_0(X)\rtimes_{\theta}\Z$, 
we give its universal property 
(see \cite[Definition 2.4]{AEE}). 

\begin{proposition}[{\cite[Example 3.2]{AEE}}]\label{cpuniv}
The $C^*$-algebra $C_0(X)\rtimes_{\theta}\Z$ 
is generated by the images of 
a $*$-ho\-mo\-mor\-phism 
$\rho^0\colon C_0(X)\to C_0(X)\rtimes_{\theta}\Z$ 
and a linear map 
$\rho^1\colon C_0(V)\to C_0(X)\rtimes_{\theta}\Z$ 
satisfying 
\benu
\item $\rho^0(f)\rho^1(g)=\rho^1(fg)$, 
\item $\rho^1(g)\rho^0(f)=\rho^1(\theta(\theta^{-1}(g)f))$, 
\item $\rho^1(g)\rho^1(h)^*=\rho^0(g\overline{h})$, 
\item $\rho^1(g)^*\rho^1(h)=\rho^0(\theta^{-1}(\overline{g}h))$\qquad\qquad
$(f\in C_0(X)$, $g,h\in C_0(V))$. 
\eenu
Moreover $C_0(X)\rtimes_{\theta}\Z$ is universal 
among such $C^*$-algebras. 
\end{proposition}

\begin{remark}\label{two}
The similar computation done after Definition \ref{DefTpl} 
shows that the conditions (i) and (ii) above are 
automatically satisfied from the conditions (iii) and (iv), 
respectively. 
\end{remark}

From a partial homeomorphism $\sigma\colon U\to V$ on $X$, 
we can define a topological graph $E=(E^0,E^1,d,r)$ 
by $E^0=X$, $E^1=U$, $r=\sigma$, and $d$ is a natural embedding. 
We have $C_d(E^1)=C_0(U)$ with the natural inner product 
and the natural right action. 
For $f\in C_0(X)=C_0(E^0)$ and $g\in C_0(U)=C_d(E^1)$, 
we have $\pi_r(f)g=(f\circ\sigma)g\in C_0(U)$. 

\begin{lemma}\label{cpph}
We have $\s{E}{0}{rg}=V$, 
and $\pi_r(f)=\theta_{g,h}$ for $f\in C_0(V)$ 
where $g,h\in C_0(U)$ satisfy $\theta^{-1}(f)=g\overline{h}$. 
\end{lemma}

\begin{proof}
By Lemma \ref{E0r}, 
$\s{E}{0}{rg}$ is the largest open subset of $E^0$ 
satisfying the property 
that the restriction of $r$ to $r^{-1}(\s{E}{0}{rg})$ 
is a proper surjection onto $\s{E}{0}{rg}$. 
Hence we have $\s{E}{0}{rg}=V$ 
because $r\colon E^1\to E^0$ is a homeomorphism onto $V\subset E^0$. 

For $f\in C_0(V)$, 
we have $\pi_r(f)g'=(f\circ\sigma)g'=\theta^{-1}(f)g'$
for $g'\in C_0(U)$. 
We also have 
$\theta_{g,h}g'=g\overline{h}g'$ 
for $g,h,g'\in C_0(U)$. 
Now the latter part is easy to see. 
\end{proof}

\begin{proposition}\label{OE=pcp}
There exists a natural isomorphism 
$\cO(E)\cong C_0(X)\rtimes_{\theta}\Z$. 
\end{proposition}

\begin{proof}
From a Cuntz Krieger $E$-pair $T=(T^0,T^1)$
on a $C^*$-algebra, 
we get a $*$-ho\-mo\-mor\-phism 
$\rho^0=T^0\colon C_0(X)\to A$ 
and a linear map $\rho^1=T^1\circ\theta^{-1}\colon C_0(V)\to A$. 
We will show that $\rho^0$ and $\rho^1$ 
satisfies four conditions in Proposition \ref{cpuniv}. 
By Remark \ref{two}, it suffices to see (iii) and (iv). 
Take $f\in C_0(X)$, $g,h\in C_0(V)$. 
For (iii), we have 
\begin{align*}
\rho^1(g)\rho^1(h)^*
&=T^1(\theta^{-1}(g))T^1(\theta^{-1}(h))^*
=\varPhi(\theta_{\theta^{-1}(g),\theta^{-1}(h)})\\
&=\varPhi(\pi_r(g\overline{h}))
=T^0(g\overline{h})=\rho^0(g\overline{h}),
\end{align*}
by Lemma \ref{cpph}. 
For (iv), we have 
$$\rho^1(g)^*\rho^1(h)=T^1(\theta^{-1}(g))^*T^1(\theta^{-1}(h))
=T^0(\theta^{-1}(\overline{g}h))=\rho^0(\theta^{-1}(\overline{g}h)).$$
We can similarly prove that 
$T^0=\rho^0$ and $T^1=\rho^1\circ\theta$ define 
a Cuntz-Krieger $E$-pair 
when two maps $\rho^0,\rho^1$ satisfy 
four conditions in Proposition \ref{cpuniv}.
Hence there exists a natural isomorphism 
$\cO(E)\cong C_0(X)\rtimes_{\theta}\Z$. 
\end{proof}

\begin{remark}
Using Lemma \ref{cpph}, 
Proposition \ref{OE=pcp} follows from \cite[Proposition 2.23]{MS}. 
\end{remark}

\subsection{\boldmath{$C^*$}-algebras associated with branched coverings}
\label{SecCov}

In \cite{DM}, 
V. Deaconu and P. S. Muhly defined a $C^*$-algebra $C^*(X,\sigma)$ 
from a branched covering $\sigma\colon X\to X$. 
They define a $C^*$-correspondence over $C_0(X)$ 
by taking a completion of $C_c(X\setminus S)$ 
where $S$ is a branch set of $\sigma$. 
This $C^*$-correspondence is the same 
as $C_d(E^1)$ obtained from a topological graph $E=(E^0,E^1,d,r)$ 
where $E^0=X$, $E^1=X\setminus S$, $d=\sigma$, 
and $r$ is a natural embedding. 
They showed that $C^*(X,\sigma)$ 
is isomorphic to the augmented Cuntz-Pimsner algebra of 
the $C^*$-correspondence $C_d(E^1)$ over $C_0(E^0)$ 
(\cite[Theorem 3.2]{DM}). 
Hence by \cite[Proposition 3.9]{Ka1}, 
we see that $C^*(X,\sigma)$ is 
isomorphic to $\cO(E)$. 

\begin{remark}
In \cite{DM}, 
a map $\sigma\colon X\to X$ was assumed to be surjective, 
but the proof of \cite[Theorem 3.2]{DM} goes well 
without assuming this (cf.\ \cite[Theorem 7]{DKM}). 
\end{remark}

\subsection{\boldmath{$C^*$}-algebras associated with singly generated topological systems}

In \cite{Re2}, J. Renault introduces the following notion. 

\begin{definition}
A {\em singly generated dynamical system} (SGDS) is 
a pair $(X,\sigma)$ where $X$ is a locally compact topological space 
and $\sigma$ is a local homeomorphism from an open subset 
$\dom(\sigma)$ of $X$ onto an open subset $\mbox{ran}(\sigma)$ of $X$. 
\end{definition}

J. Renault constructed a groupoid $G(X,\sigma)$ 
from an SGDS $(X,\sigma)$ by 
$$G(X,\sigma)=\{(x,m-n,y)\mid 
m,n\in\N, x\in\dom(\sigma^m), y\in\dom(\sigma^n), \sigma^m(x)=\sigma^n(y)\}.$$
The groupoid $G(X,\sigma)$ has a topology whose basic open sets are 
in the form 
$${\mathcal U}(U_0;m_0,m_1;U_1)=
\{(x,m_0-m_1,y)\mid (x,y)\in U_0\times U_1, 
\sigma^{m_0}(x)=\sigma^{m_1}(y)\},$$
where $U_i$ is an open subset of 
$\dom(\sigma^{m_i})$ on which $\sigma^{m_i}$ is injective for $i=0,1$. 
Note that ${\mathcal U}(U_0;m_0,m_1;U_1)$ is homeomorphic to 
$\sigma^{m_0}(U_0)\cap \sigma^{m_1}(U_1)\subset X$. 
By this topology, $G(X,\sigma)$ is a locally compact groupoid. 

We should remark that 
in \cite{Re2} a topological space $X$ in an SGDS $(X,\sigma)$ 
was not assumed to be locally compact, 
or even Hausdorff, 
but eventually $X$ was assumed to be locally compact 
(hence Hausdorff) and second countable 
in order to apply the construction in \cite{Re1} 
to the groupoid $G(X,\sigma)$. 
We do not assume that $X$ is second countable here 
because we do not need this assumption 
and we can apply the construction in \cite{Re1} 
even though $X$ is not second countable. 

In \cite{Re2}, J. Renault defined 
the $C^*$-algebra $C^*(X,\sigma)$ of an SGDS $(X,\sigma)$ 
to be the $C^*$-algebra of 
the locally compact groupoid $G(X,\sigma)$. 
In other words, 
the $C^*$-algebra $C^*(X,\sigma)$ 
is the norm closure of the $*$-algebra $C_c(G(X,\sigma))$ 
whose operations are defined by 
\begin{align*}
fg(x,m-n,y)&=\sum_{z,l}f(x,m-l,z)g(z,l-n,y)\\
f^*(x,m-n,y)&=\overline{f(y,n-m,x)}
\end{align*}
for $f,g\in C_c(G(X,\sigma))$ 
with respect to a certain norm 
(for the detail, see \cite{Re1}). 

From an SGDS $(X,\sigma)$, 
we have a topological graph $E=(E^0,E^1,d,r)$ 
by setting $E^0=X$, $E^1=\dom(\sigma)$, $d=\sigma$, 
and $r$ is a natural embedding. 
Since $r$ is a natural embedding, 
we have $\s{E}{0}{rg}=\dom(\sigma)$. 

\begin{proposition}\label{SGDS=>OE}
For an SGDS $(X,\sigma)$, 
the $C^*$-algebra $C^*(X,\sigma)$ is naturally isomorphic to 
$\cO(E)$. 
\end{proposition}

\begin{proof}
We can and will identify the open set 
$$\{(x,0,x)\in G(X,\sigma) \mid x\in X\}$$
in $G(X,\sigma)$ with $X$. 
It is routine to check that 
the embedding $C_c(X)\to C_c(G(X,\sigma))$ is 
a $*$-ho\-mo\-mor\-phism. 
Thus we get an injective $*$-ho\-mo\-mor\-phism 
$T^0\colon C_0(X)\to C^*(X,\sigma)$. 
The open set 
$$\{(x,1,\sigma(x))\in G(X,\sigma) \mid x\in \dom(\sigma)\}$$
is homeomorphic to $E^1=\dom(\sigma)$, 
and the embedding 
$T^1\colon C_c(E^1)\to C_c(G(X,\sigma))$ 
satisfies $T^1(\xi)^* T^1(\eta)=T^0(\ip{\xi}{\eta})$ 
for $\xi,\eta\in C_c(E^1)$. 
Thus we get a linear map $T^1\colon C_d(E^1)\to C^*(X,\sigma)$. 
It is not difficult to see that $T=(T^0,T^1)$ 
is an injective Toeplitz $E$-pair. 
We will show that $T$ is a Cuntz-Krieger $E$-pair. 

Let $U$ be an open subset of $E^1=\dom(\sigma)$ 
on which $d=\sigma$ is injective. 
Take $\xi,\eta\in C_c(U)\subset C_d(E^1)$ 
and set $f=\xi\overline{\eta}\in C_c(U)\subset C_0(E^0)$. 
We have $\pi_r(f)=\theta_{\xi,\eta}$ 
in a similar way to the proof of Lemma \ref{cpph}. 
We also have $T^0(f)=T^1(\xi)T^1(\eta)^*$ 
by straightforward computation. 
Thus we get $T^0(f)=\varPhi(\pi_r(f))$ 
for all $f\in C_c(U)$ and all $U\subset \dom(\sigma)=\s{E}{0}{rg}$ . 
This shows that $T^0(f)=\varPhi(\pi_r(f))$ 
for all $f\in C_0(\s{E}{0}{rg})$. 
Thus $T$ is a Cuntz-Krieger $E$-pair. 
Hence there exists a $*$-ho\-mo\-mor\-phism 
$\rho\colon \cO(E)\to C^*(X,\sigma)$. 
Since the cocycle $G(X,\sigma)\ni (x,k,y)\mapsto k\in \Z$ 
gives an action $\T\curvearrowright C^*(X,\sigma)$ 
which is a gauge action of $T$, 
the map $\rho$ is injective by Proposition \ref{GIUT}. 

The proof ends once we show that $\rho$ is surjective. 
To do so, 
it suffices to see that 
for all $(x_0,k,x_1)\in G(X,\sigma)$, 
there exists a neighborhood $W$ of $(x_0,k,x_1)$ 
such that $C_c(W)\subset C^*(X,\sigma)$ is 
in the image of $\rho$. 
Take $(x_0,k,x_1)\in G(X,\sigma)$. 
Then there exist $m_0,m_1\in\N$ 
such that 
$m_0-m_1=k$, $x_i\in \dom(\sigma^{m_i})$ for $i=0,1$ 
and $\sigma^{m_0}(x_0)=\sigma^{m_1}(x_1)$. 
For each $i=0,1$, 
take a neighborhood $U_i\subset\dom(\sigma^{m_i})$ of $x_i$ 
on which $\sigma^{m_i}$ is injective. 
Set $W={\mathcal U}(U_0;m_0,m_1;U_1)$ 
which is a neighborhood of $(x_0,k,x_1)\in G(X,\sigma)$. 
Let us set $W'=\sigma^{m_0}(U_0)\cap \sigma^{m_1}(U_1)$. 
Then $W\ni (x,k,y)\mapsto \sigma^{m_0}(x)\in W'$ 
is a homeomorphism. 
Take $f\in C_c(W)$. 
We have $f'\in C_c(W')$ such that 
$f'(\sigma^{m_0}(x))=f((x,k,y))$ for all $(x,k,y)\in W$. 
Let $X$ be the support of $f'$. 
There exist $X_0\subset U_0$ and $Y_0\subset U_1$ 
such that $\sigma^{m_0}(X_0)=\sigma^{m_1}(Y_0)=X$. 
We set $X_n=\sigma^{n}(X_0)$ for $n=1,\ldots,m_0-1$, 
and $Y_n=\sigma^{n}(Y_0)$ for $n=1,\ldots,m_1-1$. 
For $n=0,\ldots,m_0-1$, 
choose $\xi_n\in C_c(\dom(\sigma))$ 
so that $\xi_n(x)=1$ for $x\in X_n$. 
Similarly for each $n=0,\ldots,m_1-1$, 
we choose $\eta_n\in C_c(\dom(\sigma))$ 
so that $\eta_n(x)=1$ for $x\in Y_n$. 
Then it is not difficult to check 
$$f=T^1(\xi_0)\cdots T^1(\xi_{m_0-1})T^0(f')
T^1(\eta_{m_1-1})^*\cdots T^1(\eta_{0})^*\in \rho(\cO(E)).$$
Thus $C_c(W)\subset \rho(\cO(E))$. 
This completes the proof. 
\end{proof}

By Proposition \ref{SGDS=>OE}, 
all $C^*$-algebras of SGDS's are obtained as 
$C^*$-algebras of topological graphs. 
Conversely 
we will see in \cite{Ka4} that from a topological graph 
we can construct an SGDS 
so that they define the same $C^*$-algebra. 
Thus the class of $C^*$-algebras of topological graphs 
coincides with the one of SGDS's.

\begin{proposition}
The SGDS $(X,\sigma)$ is essentially free if and only if 
the topological graph $E$ is topologically free. 
\end{proposition}

\begin{proof}
Since every vertices of $E$ receives at most one edge, 
every loop has no entrances. 
Thus $E$ is topologically free if and only if 
the set of base points of loops has an empty interior. 
By Baire's theorem, 
this is equivalent to say that for every positive integer $n$ 
the set of base points of loops with length $n$ 
has an empty interior (see \cite[Proposition 6.10]{Ka2} for the detail). 
The point $x\in E^0=X$ is a base point of a loop with length $n$ 
if and only if $x\in \mbox{dom}(\sigma^n)$ 
and $\sigma^n(x)=x$. 
Thus we have shown that 
the topological graph $E$ is topologically free 
if and only if the set 
$$\{x\in \mbox{dom}(\sigma^n) \mid \sigma^n(x)=x\}$$
has an empty interior for all positive integer $n$. 
This is equivalent to the essential freeness of the SGDS $(X,\sigma)$ 
defined in \cite[Definition 2.5]{Re2}
\end{proof}

When a local homeomorphism $\sigma\colon \dom(\sigma)\to\mbox{ran}(\sigma)$ 
is a partial homeomorphism, 
the $C^*$-algebra $C^*(X,\sigma)$ is naturally isomorphic to 
the $C^*$-algebra considered in Subsection \ref{SecPH}, 
and when $\sigma$ is 
obtained by restricting a branched covering $\sigma\colon X\to X$ 
to the nonsingular set, 
$C^*(X,\sigma)$ is naturally isomorphic to 
the $C^*$-algebra considered in Subsection \ref{SecCov}.


\begin{thebibliography}{AEE}

\bibitem[AEE]{AEE}
Abadie, B.; Eilers, S.; Exel, R. {\it Morita equivalence for crossed products by Hilbert $C\sp *$-bimodules.} Trans. Amer. Math. Soc. {\bf 350} (1998), no. 8, 3043--3054.

\bibitem[CM]{CM}
Carlsen, T. M.; Matsumoto, K. {\it Some remarks on the $C^*$-algebras assosiated with subshifts.} Preprint. 

\bibitem[D]{D}
Deaconu, V. {\it Continuous graphs and C*-algebras.} Operator theoretical methods, 137--149, Theta Found., Bucharest, 2000.

\bibitem[DKM]{DKM}
Deaconu, V.; Kumjian, A.; Muhly, P. {\it Cohomology of topological graphs and Cuntz-Pimsner algebras.} J. Operator Theory {\bf 46} (2001), no. 2, 251--264.

\bibitem[DM]{DM}
Deaconu, V.; Muhly, P. S. {\it $C\sp *$-algebras associated with branched coverings.} Proc. Amer. Math. Soc. {\bf 129} (2001), no. 4, 1077--1086.

\bibitem[E1]{Ex1}
Exel, R. {\it Circle actions on $C\sp *$-algebras, partial automorphisms, and a generalized Pimsner-Voiculescu exact sequence.} J. Funct. Anal. {\bf 122} (1994), no. 2, 361--401.

\bibitem[E2]{Ex2}
Exel, R. {\it Approximately finite $C\sp *$-algebras and partial automorphisms.} Math. Scand. {\rm 77} (1995), no. 2, 281--288.

\bibitem[EL]{EL}
Exel, R.; Laca, M. {\it Cuntz-Krieger algebras for infinite matrices.} J. Reine Angew. Math. {\bf 512} (1999), 119--172.

\bibitem[K1]{Ka1}
Katsura, T. {\it A class of $C^*$-algebras generalizing both graph algebras and homeomorphism $C^*$-algebras I, fundamental results.} 
Trans. Amer. Math. Soc. {\bf 356} (2004), no. 11, 4287-4322.

\bibitem[K2]{Ka6}
Katsura, T. {\it A construction of $C^*$-algebras from $C^*$-correspondences.} Advances in Quantum Dynamics, 173--182, Contemp. Math., {\bf 335}, Amer. Math. Soc., Providence, RI, 2003.

\bibitem[K3]{Ka7}
Katsura, T. {\it Ideal structure of $C^*$-algebras associated with $C^*$-correspondences.} Preprint 2003, math.OA/0309294. 

\bibitem[K4]{Ka5}
Katsura, T. {\it $C^*$-algebras generated by scaling elements.} 
To appear in J. Operator Theory. 

\bibitem[K5]{Ka2}
Katsura, T. {\it A class of $C^*$-algebras generalizing 
both graph algebras and homeomorphism $C^*$-algebras III, 
ideal structures.} Preprint 2004, math.OA/0408190. 

\bibitem[K6]{Ka3}
Katsura, T. {\it A class of $C^*$-algebras generalizing 
both graph algebras and homeomorphism $C^*$-algebras IV, 
pure infiniteness.} IPreprint 2005, math.OA/0509343. 

\bibitem[K7]{Ka4}
Katsura, T. {\it Topological graphs and singly generated dynamical systems.} In preparation. 

\bibitem[La]{La}
Lance, E. C. {\it Hilbert $C\sp *$-modules. A toolkit for operator algebraists.} London Mathematical Society Lecture Note Series, {\bf 210}. Cambridge University Press, Cambridge, 1995.

\bibitem[Li]{Li}
Lin, H. {\it An introduction to the classification of amenable $C\sp *$-algebras.} World Scientific Publishing Co., Inc., River Edge, NJ, 2001. 

\bibitem[M1]{M1}
Matsumoto, K. {\it On $C\sp *$-algebras associated with subshifts.} Internat. J. Math. {\bf 8} (1997), no. 3, 357--374.

\bibitem[M2]{M2}
Matsumoto, K. {\it $C\sp *$-algebras associated with presentations of subshifts.} Doc. Math. {\bf 7} (2002), 1--30.

\bibitem[MS]{MS}
Muhly, P. S.; Solel, B. {\it Tensor algebras, induced representations, and the Wold decomposition.} Canad. J. Math. {\bf 51} (1999), no. 4, 850--880.

\bibitem[P]{P}
Pimsner, M. V. {\it A class of $C\sp *$-algebras generalizing both Cuntz-Krieger algebras and crossed products by ${Z}$.} Free probability theory, 189--212, Fields Inst. Commun., {\bf 12}, Amer. Math. Soc., Providence, RI, 1997. 

\bibitem[Re1]{Re1}
Renault, J. {\it A groupoid approach to $C\sp{*} $-algebras.} Lecture Notes in Mathematics, 793. Springer, Berlin, 1980.

\bibitem[Re2]{Re2}
Renault, J. {\it Cuntz-like algebras.} Operator theoretical methods, 371--386, Theta Found., Bucharest, 2000.

\bibitem[R\o]{Ro}
R\o rdam, M. {\it A purely infinite AH-algebra and an application to AF-embeddability.} Preprint.

\bibitem[RS]{RS}
R\o rdam, M.; St\o rmer, E. {\it Classification of nuclear $C\sp *$-algebras. Entropy in operator algebras.} Encyclopaedia of Mathematical Sciences, {\bf 126}. Operator Algebras and Non-commutative Geometry, {\bf 7}. Springer-Verlag, Berlin, 2002.

\bibitem[S]{Sc}
Schweizer, J. {\it Crossed products by $C\sp *$-correspondences and Cuntz-Pimsner algebras.} $C\sp *$-algebras, 203--226, Springer, Berlin, 2000.

\bibitem[T]{T}
Tomforde, M. {\it Stability of $C\sp *$-algebras associated to graphs.} Proc. Amer. Math. Soc. {\bf 132} (2004), no. 6, 1787--1795.

\end{thebibliography}
\end{document}